\documentclass[reqno,10pt]{amsart}
\usepackage{amssymb}
\usepackage{amsmath}
\usepackage{amsthm}
\usepackage[usenames]{color}
\usepackage{graphicx}
\usepackage{cite}
\allowdisplaybreaks[4]

\usepackage[colorlinks,linkcolor=blue,anchorcolor=red,citecolor=blue]{hyperref}
\usepackage[margin=1in]{geometry} 
\usepackage{marginnote}

\usepackage{color}

\newtheorem{thm}{Theorem}[section]

\newtheorem{lem}[thm]{Lemma}

\newtheorem{rem}{Remark}[section]

\numberwithin{equation}{section}

\newcommand{\R}{\mathbb{R}}

\newcommand{\ka}{\kappa}

\newtheorem{defn}{Definition}[section]

\newcommand{\vertiii}[1]{{\left\vert\kern-0.25ex\left\vert\kern-0.25ex\left\vert #1
		\right\vert\kern-0.25ex\right\vert\kern-0.25ex\right\vert}}

\begin{document}

\title[Solutions to the Fokker-Planck equation with time-dependent force and source]{Solutions to the Fokker-Planck equation with time-dependent force and source in $\mathbb R^3$}

%
\author[R.-J. Duan]{Renjun Duan}
\address[RJD]{Department of Mathematics, The Chinese University of Hong Kong,  Shatin, Hong Kong}
\email{rjduan@math.cuhk.edu.hk}

\author[F. Li]{Fucai Li}
\address[FL]{School of Mathematics, Nanjing University, Nanjing 210093, P. R. China}
\email{fli@nju.edu.cn}

\author[J. Ni] {Jinkai Ni}
\address[JKN]{Department of Applied Mathematics, The Hong Kong Polytechnic University, Hong Kong}
\email{jinkaini123@gmail.com}

 \author[Z. Zhang]{Zhu Zhang}
\address[ZZ]{Department of Applied Mathematics, The Hong Kong Polytechnic University, Hong Kong}
\email{zhuama.zhang@polyu.edu.hk}

\begin{abstract}
This paper concerns the global solutions of the kinetic Fokker-Planck (KFP) equation in $\mathbb R^3$ driven by a given external force   $E=E(t,x)$ and an inhomogeneous source term $S_0=S_0(t,x,v)$.
First, when $E$ and $S_0$ are sufficiently small in a
suitable hybrid homogeneous Besov-Sobolev framework without any time-decay assumption, we establish the global well-posedness of the Cauchy problem in the hybrid homogeneous space $L_v^2(\dot B^{ {1}/{2}}_{2,\infty}(\mathbb R^3_x)  \cap\dot H^N(\mathbb R^3_x))$ with $N\geq 3$. We  further  obtain the uniform stability of the Cauchy problem via semi-group theory and a refined energy method. Second, when both $E$ and $S_0$ are time-periodic of the period $T$, we show the existence of time-periodic solutions with the same period $T$ and their asymptotic stability as $t\to \infty$. Moreover, we derive the time-decay estimate of the perturbation around time-periodic solutions for initial data in $L^p$-norm with $ {3}/{2}< p \leq 2$. Finally, when both $E$ and $S_0$ are independent of $t$, we  establish the existence of stationary solutions and their asymptotic stability. Our results provide a unified treatment of time-periodic dynamics and stationary states for the KFP equation with small force and source terms.
\end{abstract}

	\date{\today}
	

	\maketitle
	\thispagestyle{empty}
	

\section{Introduction}

\subsection{Problem and background}
Kinetic Fokker-Planck (KFP) equation appears in a broad range of applications, including plasma physics
\cite{bit,lib}, stellar dynamics \cite{chav,cha}, polymeric fluids \cite{doi}, neurophysics \cite{Frank-2005}, and biomedicine \cite{BBJM-esaim-2005}. In the present paper, we study a KFP equation with a general prescribed external force   $E = E(t,x)$ and a general inhomogeneous source   $S_0 = S_0(t,x,v)$ in the whole space $\mathbb R^3$. More precisely, we consider the following equation:
\begin{align}\label{G1.1}  
\partial_t F + v \cdot \nabla_x F + E \cdot \nabla_v F = \nabla_v \cdot (\nabla_v F + v F) + S_0.  
\end{align}  
Here, the unknown function $F = F(t,x,v) \geq 0$ denotes a probability density distribution function for time $t \geq 0$, position $x \in \mathbb{R}^3$, and velocity $v \in \mathbb{R}^3$. The term $\nabla_v \cdot (\nabla_v F + v F)$ is the classical Fokker–Planck operator (cf. \cite{Villani-MAMS-2009}).

A natural reference case to \eqref{G1.1} is obtained by setting $E = S_0 = 0$  which reads as  the   linear KFP equation  
\begin{align}\label{2015AA}  
\partial_t F + v \cdot \nabla_x F = \nabla_v \cdot (\nabla_v F + v F).  
\end{align}  
Equation \eqref{2015AA} above is well-known in kinetic theory and is usually used as a fundamental model to motivate more complex systems. In the literature, there have been many results on well-posedness and large time behavior for \eqref{2015AA} using different approaches: hypocoercivity, semigroup,  functional-analytic framework, spectral method, 
 and variational method, see   \cite{Villani-MAMS-2009, AAMN-APDE-2024, D-2011-Nonlinearity,GMM-MSM-2017,LY-AA-2015,BS-sima25,dms,do-li,cm},  among others.

We now turn to \eqref{G1.1} with  an external force but without source, namely $S_0 = 0$ but $E \neq 0$. In this direction, many works have focused on KFP-type equations with structured force fields. Degond \cite{Degond-1986} proved the global existence of smooth solutions to the KFP equation in $\mathbb{R}^d$ ($d = 1, 2$) with self-consistent  electronic field. Bouchut and Dolbeault \cite{Bouchut-Dolbeault-1995} studied the long-time asymptotics of the KFP equation and of the Vlasov-Poisson-Fokker-Planck (VPFP) system. In the weakly self-consistent setting, Bolley, Guillin, and Malrieu \cite{Bolley-Guillin-Malrieu-2010} obtained trend-to-equilibrium and particle approximation results, while Bayraktar, Feng, and Li \cite{Bayraktar-Feng-Li-2024} proved exponential entropy dissipation under suitable Hessian-type assumptions.  Recently,  Bedrossian,    Zhao, and   Zi \cite{bzz} studied the landau damping, collisionless limit, and stability threshold for the VPFP system with weak collisions.   
For KFP equations with confining or general force fields, Cao \cite{Ccq-CMS-2019,Ccq-JEE-2021} studied the existence and uniqueness of positive normalized equilibria and convergence to equilibrium. 
Mischler and Mouhot  \cite{mis-mou} developed an  abstract quantitative method to obtain the decay rates for KFP equation with a   confinement potential in  various metrics, including Lebesgue
norms, negative Sobolev norms, and the Monge-Kantorovich-Wasserstein distance. In  \cite{DFT-2010-CMP,Choi-2016-nonlinearity}, the authors studied some related source-free kinetic models with additional force structures—such as flocking or alignment forces. The above results show that the force term plays a decisive role in the global behavior of KFP-type equations, but they usually rely on structural assumptions on the force field, such as confinement, self-consistency, potential structure, or alignment mechanisms.
The force term $E$ in  \eqref{G1.1} can be also be generated
by electromagnetic field or fluid velocity, and there are many contributions on those systems, 
for example, see  \cite{YY-SIMA-2010,chae,BG} on Vlasov-Maxwell-Fokker-Planck system and 
\cite{CDM-krm-2011,DL-2013-KRM,LMW-SIMA-2017, 
LWW-ARMA-2022,LLY-JSP-2022,FQW-M3AS-2024} on the coupled Vlasov-Fokker-Planck-Euler/Navier-Stokes system. 
We also mention that there are some results on stationary  solutions and periodic solutions to KFP equations  under some structural assumptions \cite{ah,do,gsz,DL-2015-AMSSB,dre}.

The remaining and most relevant case for the present paper is the simultaneous presence of a prescribed Vlasov-type force and an extrnal source, namely $ E\neq0$ and  $S_0\not\equiv0$. This situation is not a direct superposition of the two separate cases above. The source term contributes an inhomogeneous Duhamel term, while the force
term contains the velocity derivative $E\cdot\nabla_vF$ and therefore interacts with the hypocoercive structure of the KFP operator.  To the best of our knowledge, the global well-posedness, time-periodic dynamics, and stationary states for \eqref{G1.1} in this simultaneous prescribed-force and
prescribed-source setting have not been systematically studied before.

The time-periodic problem is another main motivation of this work. Time-periodic phenomena is a  basic phenomena in many physical process. For example,  time periodic problems for fluid dynamical equations have been extensively studied. We refer to, for example, \cite{Serrin-1959-ARMA,Mp-1991-Nonlinearity,ft-1,ft-2} for the incompressible Navier–Stokes
equations, and \cite{feia,feib,feic,Ts,Valli-1983,Deguchi-2025} for the compressible Navier–Stokes equations. In kinetic level, time-periodic solutions have been studied for several models driven either by time-periodic source terms or by time-periodic forces; see, for instance,
\cite{DL-2015-AMSSB,DLN-2026,DN-2026,ukai,ukaiy}. These works suggest that a natural way to construct time-periodic solutions is to first establish a robust {\it global-in-time} stability theory for the associated Cauchy problem and then apply the well-known Serrin's method.

In the present paper, we follow this strategy for the KFP equation \eqref{G1.1}. We first establish the global well-posedness of the Cauchy problem near the global Maxwellian and prove the uniform stability of the
corresponding global solutions by combining semigroup estimates with a refined energy method. When both $E$ and $S_0$ are $T$-periodic in time, we use Serrin's method to prove the existence and asymptotic stability of time-periodic solutions.  When proving the existence and stability of time-periodic solutions and stationary solutions to the KFP equation \eqref{G1.1} with the period $T > 0$, it is given that the external force $E(t,x)$ and the source term $S_0(t,x,v)$ are both time-periodic. That is, for any $t\in\mathbb R$ and $x\in\mathbb R^3$, $E(t,x)=E(t + T,x)$ and $S_0(t,x,v)=S_0(t+T,x,v)$. In this case, the KFP equation generated by the inhomogeneous term $S_0(t,x,v)$ and the external force $E(t,x)$ is described as the solution to the following time-periodic problem:
 \begin{equation}\label{G1.2}
\left\{\begin{aligned}
& \partial_t F_{T}+ v\cdot \nabla_{x} F_{T}+E\cdot\nabla_{v} F_{T}= \nabla_{v}\cdot(\nabla_{v}F_{T}+vF_{T})+S_0\quad \text{in}\quad \mathbb R\times \mathbb R_{x}^3\times \mathbb R_{v}^3,\\
&F_{T}(t,x,v)=F_{T}(t+T,x,v)\quad \text{for}\quad (t,x,v)\in \mathbb R\times \mathbb R_x^3\times \mathbb R_v^3.  
 \end{aligned}
 \right.
\end{equation} 
Usually, as in \cite{DL-2015-AMSSB,DLN-2026,DN-2026,ukai,ukaiy}, such a solution will be sought in a neighborhood of a global Maxwellian $M$:
\begin{align*}
M=M(v)= {(2\pi)^{-\frac{3}{2}}}\exp\Big\{ {-\frac{|v|^2}{2}} \Big\}.    
\end{align*}
Note that $M$ has a zero bulk velocity, a unit density, and a unit temperature. In this paper, 
in the same spirit as the above cited works, we shall consider the time-periodic problem \eqref{G1.2} around the global Maxwellian $M$.
For this purpose, we introduce the perturbation $f = f(t,x,v)$ by taking
\begin{align*}
F = M+\sqrt{M}f.    
\end{align*}
Then, the unknown $f=f(t,x,v)$ fulfills
\begin{align}\label{G1.3}
\partial_t f + v\cdot\nabla_x f - \mathcal{L}f = E\cdot \Big(v\sqrt{M} -  \nabla_{v} f+\frac{1}{2}  v f\Big) + S,   
\end{align}
for any $(t,x,v)\in \mathbb R\times \mathbb R^3_{x}\times \mathbb R^3_{v}$, where   the inhomogeneous term $S(t,x,v)$ is given by
\begin{align*}
S=S(t,x,v)=\frac{S_0}{\sqrt{M}},   
\end{align*}
and $\mathcal{L}$ is the linearized Fokker-Planck operator defined by
\begin{align*}
\mathcal{L}f=\frac{1}{\sqrt{M}} \nabla_{v} \cdot \Big[ M \nabla_{v}  \Big(  \frac{f}{\sqrt{M}}   \Big) \Big].    
\end{align*}

\subsection{Main results}

Before treating the time-periodic problem \eqref{G1.2}, we first establish the global well-posedness and uniform stability of the associated Cauchy problem in the vicinity of the global Maxwellian $M$. 
These results provide the crucial dynamical elements necessary for implementing Serrin's approach \cite{Serrin-1959-ARMA} to construct time-periodic solutions. To achieve the end, we consider the following Cauchy problem for the perturbation $f(t,x,v)$:
 \begin{equation} \label{G1.4}
\left\{\begin{aligned}
& \partial_t f + v\cdot\nabla_x f - \mathcal{L}f =   E\cdot \Big(v\sqrt{M} -  \nabla_{v} f+\frac{1}{2}  v f\Big)+ S, \quad (t,x,v)\in \mathbb R^+\times \mathbb R^3_{x}\times \mathbb R^3_{v},\\
 &f(0,x,v)=f_0(x,v):=\frac{F_{0}(x,v)-M}{\sqrt{M}},\quad \quad\,   (x,v)\in    \mathbb R^3_{x}\times \mathbb R^3_{v}.
 \end{aligned}
 \right.
\end{equation}

 The first main result of this paper is stated as follows.
\begin{thm}[Global well-posedness] \label{Th1}
Let $F_0(x,v)=M+\sqrt{M}f_0(x,v)\geq 0$. Assume that $f_0\in L_v^2(\dot B_{2,\infty}^{\frac{1}{2}}\cap\dot H^N)$ for an integer $N\geq 3$. There exists a constant $\delta_0>0$ such that if the inhomogeneous source term $S(t,x,v) = \frac{S_0}{\sqrt{M}}$, the external force term $E(t,x)$, and the initial datum $f_0$ satisfy
\begin{align}\label{G1.5}
\|f_0\|_{L_v^2(\dot B_{2,\infty}^{\frac{1}{2}}\cap \dot H^N)}+\|S (t)\|_{\mathcal{C} ([0,\infty);L_v^2(\dot B_{2,\infty}^{-\frac{3}{2}}\cap \dot H^N) )}+\|E (t)\|_{\mathcal{C} ([0,\infty); \dot B_{2,\infty}^{-\frac{3}{2}}\cap \dot H^N )}\leq \delta_0,
\end{align}
then the Cauchy problem \eqref{G1.4} admits a unique global solution  $f(t,x,v)\in \mathcal{C}([0,\infty);L_v^2(\dot B_{2,\infty}^{\frac{1}{2}}\cap \dot H^N))$, 
 such that $F(t,x,v)=M+\sqrt{M}f(t,x,v)\geq 0$ and
\begin{align}\label{G1.6}
\|f(t)\|_{L_v^2(\dot B_{2,\infty}^{\frac{1}{2}}\cap \dot H^N)}
\leq C_0\Big( \|f_0\|_{L_v^2(\dot B_{2,\infty}^{\frac{1}{2}}\cap \dot H^N)}
+\|S (t)\|_{\mathcal{C} ([0,\infty);L_v^2(\dot B_{2,\infty}^{-\frac{3}{2}}\cap \dot H^N) )}
+\|E (t)\|_{\mathcal{C} ([0,\infty); \dot B_{2,\infty}^{-\frac{3}{2}}\cap \dot H^N )} \Big)  , 
\end{align}
for any $t\geq 0$.  Here $C_0$ is a positive constant independent of $t$. 
\end{thm}

\begin{rem}
For the Cauchy problem \eqref{G1.4}, we use the hybrid homogeneous space $L_v^2\bigl(\dot B^{\frac12}_{2,\infty}\cap \dot H^N\bigr)$ with $N\geq3$,
instead of the usual Sobolev space $L_v^2(H^N)$. This low-frequency framework
is adapted to the semigroup estimates and allows us to handle the prescribed
force $E(t,x)$ and source term $S_0(t,x,v)$. In contrast with the
compressible Navier-Stokes-Fourier system \cite{Deguchi-2025}, where one takes
$N\geq4$, the condition $N\geq3$ is sufficient here because the nonlinear terms in \eqref{G1.4}$_1$ are of zeroth order with respect to $x$.
\end{rem}

\begin{rem}
A key difference from the Boltzmann equation studied  in \cite{DN-2026} lies in the dissipative structure of the Fokker-Planck operator. Since this operator is linear and provides direct velocity dissipation, we can close the estimates without imposing extra velocity-weighted conditions on the initial data. This feature distinguishes the present analysis on KFP equation  from  the corresponding Boltzmann equation case.
\end{rem}

\begin{rem}
A significant difference from the Boltzmann equation treated in \cite{DN-2026,DLN-2026} is that our low-frequency analysis mainly depends on the spectral analysis of the associated semi-group. In this work, we introduce, seemingly for the first time in this context, a localized semi-group approach that integrates frequency localization with semi-group decay estimates for the KFP equation.
\end{rem}

Next, we study the uniform stability of the Cauchy problem \eqref{G1.4}.

\begin{thm}[Uniform stability]\label{Th2}
Let $F_0(x,v)=M+\sqrt{M}f_0(x,v)\geq 0$ and $G_0(x,v)=M+\sqrt{M}g_0(x,v)\geq 0$. Let $0<\varepsilon<\frac{1}{2}$ and let $N\geq 3$ be an integer. Also, let $f$ and $g$ be a pair of solutions to the Cauchy problem \eqref{G1.4} with initial data $f_0(x,v)$ and $g_0(x,v)$, respectively, where the initial data $f_0$ and $g_0$, the external force $E(t,x)$, and the inhomogeneous term $S = \frac{S_0}{\sqrt{M}}$ satisfy \eqref{G1.5} in Theorem \ref{Th1}. Then, we have the following:
\begin{itemize}
  \item[(i)] If the initial data $f_0$ and $g_0$ satisfy
\begin{align}\notag
f_0 - g_0\in L_v^2(\dot B_{2,\infty}^{s_0}),
\end{align}
with $s_0\in\big(-\frac{1}{2},\frac{1}{2}\big]$, then we have the following time-weighted uniform stability estimate:
\begin{align}\label{G1.8}
\sup_{t\geq 0}(1 + t)^{\frac{s - s_0}{2}}\|(f - g)(t)\|_{L_v^2(\dot B_{2,\infty}^{s}\cap\dot H^{N-1})}\leq C_1\|f_0 - g_0\|_{L_v^2(\dot B_{2,\infty}^{s_0}\cap\dot H^{N-1})},
\end{align}
for any $s\in\big[-\frac{1}{2}+\varepsilon,\frac{3}{2}-\varepsilon \big]$ with $s\geq s_0$.
  \item[(ii)] 
If the initial data $f_0$ and $g_0$ satisfy
\begin{align}\notag
f_0 - g_0\in L_v^2(\dot B_{2,\infty}^{s}),
\end{align}  
with $s\in \big(\frac{1}{2},\frac{3}{2}-\varepsilon\big]$, then we have the following uniform stability estimate
\begin{align}\label{DLNZ-1.10}
 \sup_{t\geq0}
\| ( f-g)(t)\|_
{L_v^2(\dot B^s_{2,\infty}\cap\dot H^{N-1})}
\leq
C_1\|f_0-g_0\|_
{L_v^2(\dot B^s_{2,\infty}\cap\dot H^{N-1})}.
\end{align}
\end{itemize}
Here, in both \eqref{G1.8} and \eqref{DLNZ-1.10}, $C_1$ is a positive constant independent of $t$. 
\end{thm}

\begin{rem}
Unlike \cite{Deguchi-2025} on the uniform stability of the Navier-Stokes-Fourier system, we cannot handle the case for $s\in\big [-\frac{3}{2}+\varepsilon,\frac{3}{2}-\varepsilon\big]$. Since we need to overcome the weighted $v$ and the derivative $\nabla_{v}$ from $\frac{1}{2}E\cdot v\widetilde f$ and $E\cdot\nabla_{v}\widetilde f$, we employed different methods in our proof. Due to technical limitations (see \eqref{techinique}), we can only obtain the result for $s\in\big [-\frac{1}{2}+\varepsilon,\frac{3}{2}-\varepsilon\big]$.
\end{rem}

Since we have established the global classical solutions to the Cauchy problem \eqref{G1.4}, we can further conduct the result on the   time-periodic problem \eqref{G1.2}. For the sake of brevity, let $F_T = M+\sqrt{M}f_{T}$. With the help of \eqref{G1.3}, we rewrite the problem \eqref{G1.2} as
 \begin{equation}\label{G1.9}
\left\{\begin{aligned}
& \partial_t f_{T}+ v\cdot \nabla_{x} f_{T}-\mathcal{L}f_{T}\\
&\quad= E\cdot \Big(v\sqrt{M}-\nabla_{v}f_{T}+\frac{1}{2}vf_{T}\Big)+S \quad \text{in}\quad \mathbb R\times \mathbb R_{x}^3\times \mathbb R_{v}^3,\\
&f_{T}(t,x,v)=f_{T}(t+T,x,v)\quad\, \text{for}\quad (t,x,v)\in \mathbb R\times \mathbb R_x^3\times \mathbb R_v^3,  
 \end{aligned}
 \right.
\end{equation} 
with $S(t,x,v)=S(t+T,x,v)$ and $E(t,x)=E(t+T,x)$. Below, we state the existence and stability of the time-periodic problem \eqref{G1.9} as follows.

\begin{thm}[Existence and stability of the time-periodic solutions]\label{Th1.3}
Let $T > 0$ and let $N\geq 3$ be an integer. There exists a small $\sigma > 0$ such that if the time-periodic inhomogeneous source term $S = \frac{S_0}{\sqrt{M}}$ and the time-periodic force $E$, both of the same period $T > 0$, satisfy
\begin{align*}
\|E (t)\|_{\mathcal{C} (\mathbb{R}; \dot B_{2,\infty}^{-\frac{3}{2}}\cap \dot H^N) }+\|S (t)\|_{\mathcal{C} (\mathbb{R};L_v^2(\dot B_{2,\infty}^{-\frac{3}{2}}\cap \dot H^N) )}\leq \sigma,
\end{align*}
then we have the following  results:
\begin{itemize}
    \item (Existence of time-periodic solutions) The problem \eqref{G1.9} admits a unique time-periodic solution $f_T$ of the period $T$ satisfying
    \begin{align}\label{Th1.3.ad1}
    \|f_T\|_{\mathcal{C}(\mathbb{R};L_v^2(\dot B^{\frac{1}{2}}_{2,\infty}\cap \dot H^N))}  \leq C_2\sigma.
    \end{align}
    \item (Stability of time-periodic solutions) Consider the Cauchy problem \eqref{G1.4} and suppose that the  initial datum $f_0(x,v)$ satisfy
    \begin{align*}
    \|f_0\|_{L_v^2(\dot B_{2,\infty}^{\frac{1}{2}}\cap \dot H^N)}\leq \sigma,
    \end{align*}
   for  $\sigma>0$ further sufficiently small and
    \begin{align*}
    f_0 - f_T(0) \in L^p,
    \end{align*}
    for some $\frac{3}{2}< p\leq 2$. Then we have the time decay estimate:
    \begin{align}\label{G1.15}
    \|(f - f_T)(t)\|_{L_v^2(\dot B^{s}_{2,1})}\leq C_2 (1 + t)^{-\frac{s}{2}-\frac{3}{2}(\frac{1}{p}-\frac{1}{2})}\|(f - f_T)(0)\|_{L_v^2( L^p\cap\dot H^N)},
    \end{align}
    for any $s\in\big(-\frac{1}{2}+\varepsilon,\frac{3}{2}-\varepsilon\big)$ with $\frac{s}{2}+\frac{3}{2}\big(\frac{1}{p}-\frac{1}{2}\big)>0$ and  $0<\varepsilon<\frac{1}{2}$, where $f$ is the unique global solution of \eqref{G1.4}  satisfying
    \begin{align*}
    \|f(t)\|_{\mathcal{C} ([0,\infty);L_v^2(\dot B_{2,\infty}^{\frac{1}{2}}\cap \dot H^N))}\leq C_2\sigma.
    \end{align*}
     Here, $C_2$ is a positive constant independent of $t$. 
\end{itemize}
    
\end{thm}

\begin{rem}
By leveraging the embedding $\dot B^{s}_{2,1}(\mathbb{R}^3)\hookrightarrow  \dot H^s(\mathbb{R}^3)$ and the decay estimate \eqref{G1.15}, we are able to derive the following optimal time-decay estimate in the $L_v^2(\dot H^s)$-norm for $(f - f_{T})(t)$ in the sense of Matsumura and Nishida \cite{MN-PJASAMS-1979}:
\begin{align*}
    \|(f - f_T)(t)\|_{L_v^2(\dot H^{s}) }\lesssim (1 + t)^{-\frac{s}{2}-{\frac{3}{2}}({\frac{1}{p}}-{\frac{1}{2}})}\|(f - f_T)(0)\|_{L_v^2( L^p\cap\dot H^N)},
\end{align*}
for $\frac{3}{2}<p \leq 2$  and   $s\in \big(-{\frac{1}{2}}+\varepsilon,{\frac{3}{2}}-\varepsilon \big)$ such that $ \frac{s}{2}+\frac{3}{2} \big(\frac{1}{p}-\frac{1}{2} \big)>0$.
\end{rem}

As an immediate consequence of Theorem \ref{Th1.3} established above, whenever $E$ and $S$ in the problem \eqref{G1.4} are independent of time, that is $E=E(x)$ and $S=S(x,v)$, we can further give the following   result for the stationary problem concerning the KFP equation:
\begin{align}\label{G1.11}
  v\cdot\nabla_x \bar{f}- \mathcal{L}\bar{f} =   E\cdot \Big(v\sqrt{M} -  \nabla_{v} \bar{f}+\frac{1}{2}  v \bar{f}\Big)+ S, \quad (x,v)\in   \mathbb R^3_{x}\times \mathbb R^3_{v}    ,
\end{align}
with the far-field boundary condition $\bar{f}(x,v)=\frac{\bar{F}(x,v) - M}{\sqrt{M}}\to 0$ as $|x|\to\infty$, for any $v\in\mathbb{R}_{v}^3$, where corresponding to \eqref{G1.1},  $\bar F(x,v)=M+\sqrt{M}\bar f (x,v)$ satisfies 
\begin{align*}
v \cdot \nabla_x \bar{F} + E \cdot \nabla_v \bar{F} = \nabla_v \cdot (\nabla_v \bar{F} + v \bar{F}) + S_0.  
\end{align*}

\begin{thm}[Existence and stability of the  stationary solutions] \label{cor1.4}
Let $N\geq 3$ be an integer. There exists a small constant $\delta^{*}>0$ such that if the external force $E(x)$ and the inhomogeneous source term $S(x,v)$ satisfy
\begin{align*} 
\|E\|_{\dot B^{-\frac{3}{2}}_{2,\infty}\cap \dot H^N}+\|S\|_{L_v^2(\dot B^{-\frac{3}{2}}_{2,\infty}\cap \dot H^N)}\leq   \delta^{*},
\end{align*}
then the stationary problem \eqref{G1.11} has a unique solution $\bar f (x,v)$  such that $\bar F(x,v)=M+\sqrt{M}\bar f (x,v)\geq 0$ and
\begin{align*}
\|\bar f\|_{L_v^2(\dot B^{\frac{1}{2}}_{2,\infty}\cap\dot H^N)}\leq C_3  \delta^*.
\end{align*}
If the initial datum $f_0(x,v)$ of the Cauchy problem \eqref{G1.4} further satisfy that $F_0=M+\sqrt{M}f_0(x,v)\geq 0$ and
\begin{align*}  
\|f_0\|_{L_v^2(\dot B^{\frac{1}{2}}_{2,\infty}\cap \dot H^N)} \leq \delta^{*},
\end{align*}
with $\delta^{*}>0$ further sufficiently small and if it holds
\begin{align*}    
f_0 - \bar f (0) \in L^p,
\end{align*}    
for some $\frac{3}{2}< p\leq 2$, then we have the time decay estimate:    
\begin{align*}    
\|(f - \bar f)(t)\|_{L_v^2(\dot H^{s})}\leq C_3 (1 + t)^{-\frac{s}{2}-\frac{3}{2}(\frac{1}{p}-\frac{1}{2})}\|(f - \bar f )(0)\|_{L_v^2( L^p\cap\dot H^N)},\quad \forall \,t\geq 0,
\end{align*}    
for any $s\in\big(-\frac{1}{2}+\varepsilon,\frac{3}{2}-\varepsilon\big)$ with $\frac{s}{2}+\frac{3}{2}\big(\frac{1}{p}-\frac{1}{2}\big)>0$ and  $0<\varepsilon<\frac{1}{2}$, where $f$ is the unique global solution of \eqref{G1.4}   satisfying   
\begin{align*}    
\|f(t)\|_{\mathcal{C} ([0,\infty);L_v^2(\dot B_{2,\infty}^{\frac{1}{2}}\cap \dot H^N))}  \leq C_3\delta^{*}.
\end{align*}
Here, $C_3>0$ is positive constant independent of $t$.

\end{thm}

\begin{rem}
An early classic reference connecting periodic and stationary problems is a seminal piece of work  Valli \cite{Valli-1983}, who investigated periodic and stationary solutions to the compressible isentropic Navier-Stokes equations by utilizing  a stability method.
\end{rem}

\begin{rem}  
It should be emphasized that the external force field under consideration  need not  be a potential field. Indeed, as pointed out in \cite{DN-2026}, for $m > 2$ and $0 <  \epsilon  \ll 1$, a non-potential stationary force field may take the form of
\begin{align*}
 E_\epsilon(x)=\frac{\epsilon(-x_2,x_1,0)}{\langle x\rangle^{2m}},\qquad x\in\mathbb R^3,   
\end{align*}
which is rotational. More precisely, one can compute that
\begin{align*}
 \nabla_x\times E_\epsilon\not\equiv0,   
\end{align*}
and  hence $E_\epsilon$ turns out to be non-potential. Moreover,
\begin{align*}
E_\epsilon\in\dot B^{-3/2}_{2,\infty}(\mathbb R^3)\cap \dot H^N(\mathbb R^3),\qquad\|E_\epsilon\|_{\dot B^{-3/2}_{2,\infty}\cap \dot H^N}\leq C_{m,N} \epsilon .    
\end{align*}
Therefore, by choosing $ \epsilon $ small enough, $E_\epsilon$ meets the force assumption in our main results. This indicates that our results are applicable to rotational, non-conservative external forces.
\end{rem}

\subsection{The strategies in our proofs}
The difficulty in proving Theorem \ref{Th1} arises from the nonlinear terms that contain weighted $v$ or the derivative $\nabla_{v}$, such as $\frac{1}{2}E\cdot v f$ and $E\cdot \nabla_{v}f$ in \eqref{G1.4}. When directly applying the semi-group theory, it is extremely difficult to handle these terms because the decay of our semi-group is related to the physical space rather than the velocity space. To overcome this difficulty, inspired by \cite{CDM-krm-2011} for the Vlasov-Fokker-Planck-Euler system, we decompose the nonlinear term $S_{f}:=-E\cdot\nabla_{v}f+\frac{1}{2}E\cdot v f$ as follows (see Remark \ref{remB.1} below):
\begin{align*}
S_{f} =\nabla_{v}\cdot G-\frac{1}{2}v\cdot G+E\cdot av\sqrt{M},
\end{align*}
where $G=-E\{\mathbf{I}-\mathbf{P}_0\}f$.
By leveraging the coercivity of $-\mathcal{L}$ and the property of $\mathbf{P}$ 
given by \begin{align*}\nabla_{v}\cdot G-\frac{1}{2}v\cdot G \perp {\rm Rang}\,\mathbf{P},    
\end{align*} 
we are able to prove (see Lemma \ref{LB.2}) that
\begin{align*} 
\|\widehat{f}(t,\xi,v)\|_{L_{v}^2}^2\lesssim   e^{-\lambda_{2}t\min\{1,|\xi|^2\}}  {}\|\widehat{f_0}(\xi,v)\|_{L_v^2}^2+\int_0^t  e^{-\lambda_{2}(t-\tau)\min\{1,|\xi|^2\}}\big(\|\widehat{ G}(\tau)\|_{L_v^2}^2+\|\nu^{-\frac{1}{2}}\widehat{ h}(\tau)\|_{L_v^2}^2  \big) {\rm d}\tau.
\end{align*}
This implies that $v$ and $\nabla_{v}$ can be eliminated via the Fourier analysis techniques. Consequently, we recognize that for the KFP equation \eqref{G1.4}, we don't need to estimate the mixed derivatives or the $v$-weighted energy estimates of $f(t,x,v)$ to enclose our {\it uniform a priori estimates}. 

To derive the {\it a priori estimates} of $f(t,x,v)$, we decompose it into low- and high-frequency components as follows (abbreviate $f(t,x,v)$ as $f(t)$):
\begin{align*}
f(t)=f_{L}(t)+f_{H}(t).
\end{align*}
By the Duhamel principle, we denote the low-frequency part of the solution $f_{L}(t)$ as
\begin{align*}
f_{L}(t)=&\, e_{L}^{t\mathcal B}f_{0}+\int_0^t e_L^{(t-\tau)\mathcal{B}} \Big[ E\cdot\Big( \frac{1}{2}vf-\nabla_{v}f \Big) (\tau)+E\cdot v\sqrt{M}(\tau)+ S(\tau)    \Big]{\rm d}\tau\nonumber\\
=&\, e_{L}^{t\mathcal B}f_{0}+ \int_0^t e_L^{(t-\tau)\mathcal{B}} \big(   E\cdot(1+a) v\sqrt{M}+ S     \big)(\tau){\rm d}\tau+\int_0^t e_L^{(t-\tau)\mathcal{B}} \Big( \nabla_{v}\cdot G-\frac{1}{2}v\cdot G\Big)(\tau){\rm d}\tau \nonumber\\
=&\, e_{L}^{t\mathcal B}f_{0}+  \mathcal{I}_1+\mathcal{I}_2.
\end{align*}
Although both $\mathcal I_1$ and $\mathcal I_2$ require the use of duality argument \eqref{duality} for processing, the methods are different. On the one hand, for the source term $S$ in $\mathcal{I}_1$, we need to prove the following fact:
\begin{align}\label{T1}
\bigg\|\int_0^t  e_L^{(t - \tau)\mathcal{B}} S(\tau)    {\rm d}\tau     \bigg\|_{L_v^2(\dot B_{2,\infty}^{{\frac{1}{2}}})}\lesssim \sup_{0\leq t\leq T_1}\|S(t)\|_{L_v^2(\dot B_{2,\infty}^{-{\frac{3}{2}}})}.
\end{align}
The proof of \eqref{T1} is non-trivial in semi-group theory. The above inequality is equivalent to exchanging the second-order spatial derivative for the condition of integrability. However, the second-order spatial derivative is approximately equal to the rate of $(1 + t)^{-1}$, and it does not satisfy the $L^1$ integrability condition. Therefore, our function framework in the low-frequency for $\dot B_{2,\infty}^{1/2}$ is critical. This idea is inspired by the works \cite{Deguchi-2024-Math Ann,Deguchi-2025} regarding the compressible Navier-Stokes system. For the term $\mathcal{I}_2$, noting that the estimate of $G$ in Lemma \ref{LB.3} is square-integrable, thus, we need a precise Littlewood-Paley decomposition to prove the following inequality (see the proof in Lemma \ref{L3.1}):
\begin{align}\label{T2}
\bigg\|\int_0^t  e_L^{(t - \tau)\mathcal{B}} \Big( \nabla_{v}\cdot G-\frac{1}{2}v\cdot G      \Big)(\tau){\rm d}\tau     \bigg\|_{L_v^2(\dot B_{2,\infty}^{{\frac{1}{2}}})}\lesssim \sup_{0\leq t\leq T_1}\|G(t)\|_{L_v^2(\dot B_{2,\infty}^{-{\frac{1}{2}}})},
\end{align}
which is quite different from $\mathcal{I}_1$ for which only Fourier analysis is enough. Moreover,   \eqref{T2} demonstrates that the function $G$ exhibits a faster rate of the first spatial derivative compared to   \eqref{T1}.
For the high-frequency part of $f$, we mainly employ the delicate energy method and restrict it in the framework of $L_v^2(\dot H^1\cap\dot H^N)$ with $N\geq 3$. Consequently, we establish the following Lyapunov-type inequality (see \eqref{G3.23}): 
\begin{align*} 
\frac{{\rm d}}{{\rm d}t}  \|f(t)\|_{L_v^2(\dot H^1\cap \dot H^N)}^2+ \kappa_5 \|f(t)\|_{L_v^2(\dot H^1\cap \dot H^N)}^2 \lesssim \sup_{{0\leq t\leq T_1}}\|E(t)\|_{ H^N}^2  +\sup_{{0\leq t\leq T_1}}\|S(t)\|_{L_v^2(H^N)}^2+\|f_0\|_{L_v^2(\dot B^\frac{1}{2}_{2,\infty})}^2 , 
\end{align*}
for any $0\leq t\leq T_1$. By applying Gr\"{o}nwall's inequality and incorporating the low-frequency estimate \eqref{G3.2}, we obtain the following uniform estimate:
\begin{align*}
 \sup_{{0\leq t\leq T_1}}\|f(t)\|_{L_v^2(\dot B^{\frac{1}{2}}_{2,\infty}\cap \dot H^N)} 
 \lesssim&\, \sup_{0\leq t\leq T_1} \Big\{\|E(t)\|_{  \dot B^{-\frac{3}{2}}_{2,\infty}\cap\dot H^N}+\|S(t)\|_{L_v^2( \dot B^{-\frac{3}{2}}_{2,\infty}\cap\dot H^N)}  \Big\}+\|f_0\|_{L_v^2(\dot B^{\frac{1}{2}}_{2,\infty}\cap \dot H^N)},
\end{align*}
which, together with the Continuity argument, leads to the global strong solution $f(t,x,v)$ of the Cauchy problem \eqref{G1.4}.

In the proof of Theorem \ref{Th2}, we first consider the following {difference} equation (see \eqref{G4.1}):
\begin{align}\label{T3}
  \partial_t \widetilde f+v\cdot\nabla \widetilde f-\mathcal{L} \widetilde f=&\, E\cdot \nabla_{v} \widetilde f+\frac{1}{2}E\cdot v\widetilde f\nonumber\\
=:&\, \nabla_v \cdot \widetilde G -\frac{1}{2}v\cdot \widetilde G+E\cdot \widetilde av \sqrt{M}.   
 \end{align}
Next, we establish a time-weighted estimate of \eqref{T3}. At low frequencies, we make full use of the semi-group theory, and at high frequencies, we continue to employ the precise energy method. It is crucial to note that even though we manage to overcome the $v$-weighted growth and the first-order velocity derivative $\nabla_{v}$, the range of the decay rates {is} reduced as a result. The decay rates we have obtained are in $L_v^2(\dot B_{2,1}^{s})$ with $s \in \big(-\frac{1}{2}+\varepsilon, \frac{3}{2}-\varepsilon\big)$, and this range is different from that in \cite{Deguchi-2025} where $s \in \big(-\frac{3}{2}+\varepsilon, \frac{3}{2}-\varepsilon\big)$.
In Lemma \ref{L4.1}, since $L_v^2(\dot B^{\frac{1}{2}}_{2,\infty})$ is the critical space, we need to discuss $s$ for three cases: $s\in\big[-\frac{1}{2}+\varepsilon,\frac{1}{2}\big)$, $s = \frac{1}{2}$, and $s\in \big(\frac{1}{2},\frac{3}{2}-\varepsilon\big]$, and consequently derive the estimate \eqref{G4.2} for the low-frequency part. With the estimate \eqref{G4.11} in Lemma \ref{L4.2} for high frequency, we immediately prove the uniform stability of the Cauchy problem \eqref{G1.4}.

To prove the existence and stability of the time-periodic solutions of the problem \eqref{G1.9}, we mainly utilize the classical Serrin’s method \cite{Serrin-1959-ARMA} (see also \cite{Mp-1991-Nonlinearity,DL-2015-AMSSB,ukai} for further developments). We establish the Cauchy sequence $\{f^*(nT)\}_{n\geq 1}$ in the hybrid homogeneous space $L_v^2(\dot B_{2,\infty}^{1} \cap \dot H^{N - 1})$. By the Banach completeness theorem and Fatou's lemma, we find the limit $f_{\infty}^*$. Then, we assume that $f_{T}(t)$ is a global solution of the problem \eqref{G1.4}, {subject to the condition that  $f_{T}(0)=f_{\infty}^*$.} {With this in hand}, we can prove that $f_T (t)$ is a time-periodic solution of the time-periodic problem \eqref{G1.9} through the {\it local-in-time uniqueness} argument. By the interpolation inequality and the embedding $L^p(\mathbb{R}^3) \hookrightarrow \dot B_{2,\infty}^{-\frac{3}{2}(\frac{1}{p}-\frac{1}{2})}(\mathbb{R}^3) $, we establish the time-asymptotic stability of $f$ for $p\in \big(\frac{3}{2},2\big]$. Hence, we complete the proof of Theorem \ref{Th1.3}. 

Finally, to prove the existence and stability of the stationary solutions of KFP equation \eqref{G1.11}, if we further assume that the external $E(x)$ and $S(x,v)$ are independent of $t$, by defining $f^\tau_{T}(t,x,v):=f_{T}(t+\tau,x,v)$, we can prove that $f^\tau_{T}(t,x,v)$ is $T$-periodic for any $\tau\in \mathbb R$. Combining the uniqueness of $f_{T}$ in Theorem \ref{Th1.3}, we further have $f_{T}^\tau = f_{T}$ for any $\tau\in \mathbb R$. By selecting $t = 0$, we consequently obtain the unique stationary solution $\bar f(x,v)=f_{T}(0,x,v)$ to the stationary problem \eqref{G1.11}.

\subsection{Outline of this paper}
The rest of this paper is organized as follows. In Section \ref{sec2}, we introduce some notations, the Littlewood-Paley decomposition, the macro-micro decomposition, and the properties of $\mathcal{L}$, and then present the time-decay estimates of the semi-group. In Section \ref{sec3}, we respectively establish the {\it uniform a priori estimates} of $f(t,x,v)$ at low and high frequencies, and further prove the global well-posedness of the Cauchy problem  \eqref{G1.4}. In Section \ref{sec4}, by using the semi-group theory and the delicate energy method at low and high frequencies respectively, we establish the uniform stability of the Cauchy problem \eqref{G1.4}. In Section \ref{sec5}, by utilizing the classical Serrin's approach, we prove the existence and stability of the time-periodic solution to the periodic problem \eqref{G1.9}. Finally, in Section \ref{sec6}, based on Theorem \ref{Th1.3}, we obtain the existence and stability of stationary solutions of the KFP equation \eqref{G1.11}.

\section{Preliminaries}\label{sec2}

In this section, we introduce some notations, the homogeneous Besov spaces, the macro-micro decomposition,  the properties of the operator $\mathcal{L}$, and the time decay estimates of semi-group, which are frequently utilized throughout this paper.

\subsection{Notations}
The symbol $C$ represents a general positive constant that is independent of the time $t$ and may vary from one line to another. The notation $A\backsim B$ denotes the equivalence relation $\frac{A}{C}\leq B\leq CA$, where $C$ is a  positive constant. The notation $A\lesssim B$ signifies that $A\leq CB$ for some constant $C>0$.
We shall use $\langle\cdot,\cdot\rangle$ and $(\cdot,\cdot)$ to denote the $L^2$ inner product in $\mathbb{R}^3_v$ and $\mathbb{R}_x^3$, respectively, along with the corresponding $L^2$ norms $|\cdot|_2$ and $\|\cdot\|_{L^2}$.
Furthermore, we denote the $L^2$ inner product in $\mathbb{R}^3_x\times \mathbb{R}_{v}^3$ by $\langle\cdot,\cdot\rangle_{x,v}$, along with its corresponding $L^2$ norm $\|\cdot\|_{L_{x,v}^2}$. For simplicity, we use $L^2$ and $L^2_{v}$ to represent $L^2(\mathbb R_x^3)$ and $L^2(\mathbb R_v^3)$, respectively. We define $L_v^2(L^p)=L^2(\mathbb R^3_v;L^p(\mathbb R_x^3))$ as the Lebesgue mixed velocity-spatial space, which is equipped with the norm $\|\cdot\|_{L_v^2(L^p)}$ for $1\leq p\leq \infty$.
Similarly, we adopt the notations $L_v^2(H^s )$, $L_v^2(\dot H^s)$, and $L_v^2(\dot B^s_{p,r})$ for $s\in \mathbb{R}$ and $1 \leq p, r \leq \infty$.

As in \cite{CDM-krm-2011}, we  define the velocity function $\nu(v):=1 + |v|^2$ and denote the corresponding norms $|\cdot|_{\nu}$ and $\|\cdot\|_{\nu}$ as follows:
\begin{align*}
|g|_{\nu}^2:= \int _{\mathbb{R}^3} (|\nabla_v g(v)|^2+\nu(v)|g(v)|^2){\rm d}v ,\quad \text{for}\quad g=g(v),   
\end{align*}
and
\begin{align*}
\|g\|_{\nu}^2:= \int_{\mathbb{R}^3_x}\!\int _{ \mathbb{R}^3_v} (|\nabla_v g(x,v)|^2+\nu(v)|g(x,v)|^2){\rm d}v{\rm d}x ,\quad \text{for}\quad g=g(x,v).
\end{align*}
For a multi-index $\alpha = (\alpha_1, \alpha_2, \alpha_3)$, we precisely define $\partial^\alpha = \partial^\alpha_x = \partial^{\alpha_1}_{x_1} \partial^{\alpha_2}_{x_2} \partial^{\alpha_3}_{x_3}$. The length of $\alpha$ is  defined as $|\alpha| = \alpha_1 + \alpha_2 + \alpha_3$. For the sake of brevity, we write $\partial_i$ instead of $\partial_{x_i}$ for each $i = 1, 2, 3$. Additionally, we use the symbols $\Delta$ and $\nabla$ to denote $\Delta_x$ and $\nabla_x$ respectively.

For a Banach space $Z$ and an interval $I\subset\mathbb{R}$, we define $\|(g,h)\|_{Z}$ as $\|g\|_{Z}+\|h\|_{Z}$, and $\mathcal{C}(I; Z)$ as the set of continuous functions on $I$ whose values lie in $Z$. For any time $\mathcal{T} > 0$ and $1 \leq p \leq \infty$, $L^p(0,\mathcal{T};Z)$ denotes the space of measurable functions $f: [0,\mathcal{T}] \to Z$ such that the mapping $t \mapsto \|f(t)\|_{Z}$ belongs to $L^p(0,\mathcal{T})$,  equipped with the norm $\|\cdot\|_{L^p(0,\mathcal{T};Z)} := \|\cdot\|_{L_{\mathcal{T}}^p(Z)}$.
We denote $\mathcal{S}$ as the set of all Schwartz functions on $\mathbb{R}^3_{x}$, and $\mathcal{S}^\prime$ as the set of all tempered distributions on $\mathbb{R}^3_{x}$. For a function $g(t,x,v) \in \mathcal{S}(\mathbb{R}^3_{x})$, its Fourier transform with respect to $x$ is defined as
\begin{align*}
\widehat{g}(t,\xi,v)=\mathcal{F}g(t,\xi,v) = \int_{\mathbb{R}^3} e^{-{\rm i} x\cdot\xi} g(t,x,v) \mathrm{d}x, \quad \text{where }\quad x\cdot \xi = \sum_{j = 1}^3 x_{j}\xi_{j},  
\end{align*}
for all $\xi \in \mathbb{R}^3$. Here, ${\rm i} = \sqrt{-1}$ represents the imaginary unit. Besides, $\mathcal{F}^{-1}$ denotes the inverse Fourier transform.

\subsection{The Littlewood-Paley decomposition and homogeneous Besov spaces}
We briefly introduce the Littlewood–Paley decomposition and the homogeneous Besov spaces,  for more details, interested reader can refer to \cite[Chapter 2]{BCD-Book-2011}.
Let $\chi(\xi)$ be a smooth, radially symmetric, non-increasing function whose support is within the ball $B(0,\frac{4}{3})$ and which satisfies $\chi(\xi)\equiv 1$ in the ball $B(0,\frac{3}{4})$. 
Define $\varphi(\xi):=\chi(\frac{\xi}{2})-\chi(\xi)$. Then, $\varphi$ satisfies the following properties:
\begin{align*}
\sum_{j \in \mathbb{Z}} \varphi(2^{-j}\xi) = 1 \quad \text{for all }\quad \xi \neq 0, \quad \text{and} \quad \mathrm{supp}\,\varphi \subset \Big\{ \xi \in \mathbb{R}^3 \Big| \frac{3}{4} \leq |\xi| \leq \frac{8}{3} \Big\}.
\end{align*}
For any $j\in \mathbb{Z}$, we define the following homogeneous dyadic blocks $\dot{\Delta}_j$ as follows:
\begin{align*}
\dot{\Delta}_j g := \mathcal{F}^{-1}\left( \varphi(2^{-j}\cdot)\mathcal{F}(g) \right) = 2^{3j} h(2^j\cdot) \ast g,  
\end{align*}
where $ h := \mathcal{F}^{-1}\varphi$.

Let $\mathcal{P}$ denote the set of all polynomials on $\mathbb{R}^3$. Define the space of tempered distributions on $\mathbb{R}^3$ modulo polynomials as follows:
\begin{align*}
\mathcal{S}'_h (\mathbb R^3):= \mathcal{S}' / \mathcal{P}  (\mathbb R^3) . 
\end{align*}
Then, for any $g \in \mathcal{S}'_h (\mathbb R^3)$, the homogeneous Littlewood-Paley decomposition holds:
\begin{align*}
g= \sum_{j \in \mathbb{Z}} \dot{\Delta}_j g,
\end{align*}
where the dyadic blocks satisfy
\begin{align*}
\dot{\Delta}_j \dot{\Delta}_k g = 0 \quad \text{whenever }\quad |j - k| \ge 2.
\end{align*}
For any $j\in \mathbb{Z}$, we further define the lower-frequency cut-off operator $\dot{S}_j$ as
\begin{align*}
\dot{S}_j g=\sum_{j^\prime< j}\dot{\Delta}_{j^\prime}g, \quad \text{for any }  g\in \mathcal{S}^{\prime}(\mathbb{R}^3).
\end{align*}
  For any $g\in \mathcal{S}^\prime(\mathbb{R}^3)$, we also denote 
\begin{align*}
g_{L} :=\dot{S}_{j_0}g,\quad \text{and}\quad g_{H} := g - \dot{S}_{j_0}g,    \quad \text{for some}\quad j_0\in \mathbb{Z}.
\end{align*}
 
With the help of those dyadic blocks, we now present the following definition of homogeneous Besov spaces.
\begin{defn}
For any $ s \in \mathbb{R} $ and $ 1 \leq p, r\leq \infty $, the homogeneous Besov space $ \dot{B}_{p,r}^s$ is defined as
\begin{align*}
\dot B_{p,r}^s:=\big\{ g\in \mathcal{S}^\prime_h \,\big{|}\, \|g\|_{\dot{B}_{p,r}^s}:=\big\|\{2^{ks}\|\dot\Delta_k g\|_{L^p}\}_{k\in\mathbb Z}\big\|_{l^r}<\infty\big\}. 
\end{align*}
In particular, when $p = r = 2$, it holds that $\dot B^{s}_{2,2}=\dot H^s$.
\end{defn}

We also need   a class of mixed velocity-spatial homogeneous Besov spaces, which can be described as follows.
\begin{defn}
For any $ s \in \mathbb{R} $ and $ 1 \leq p, r\leq \infty $, we adopt the notation $L_v^2(\dot B_{p,r}^s)$ to represent the space $L^2( \mathbb R^3_{v} ; \dot B_{p,r}^s(\mathbb R^3_{x}))$,  defined as 
\begin{align*}
 L_v^2(\dot B_{p,r}^s):=\bigg\{ g\in L_v^2(\mathcal{S}^\prime_h) \,\big{|}\, \|g\|_{L_v^2(\dot{B}_{p,r}^s)}:=\bigg(\int_{\mathbb R^3} \|g\|_{\dot B_{p,r}^s}^2{\rm d}v\bigg)^\frac{1}{2} <\infty\bigg\}.    
\end{align*}

\end{defn}

The first lemma concerns Bernstein's inequalities.

\begin{lem}\label{L2.1}{\rm(\!\!\cite[Chapter 2]{BCD-Book-2011})}
Let $0 < r < R$, $1\leq p\leq r\leq \infty$, and $m\in \mathbb{N}$. Define the ball $\mathcal{B}=\{\xi\in\mathbb{R}^{3}: |\xi|\leq  R\}$ and the annulus $\mathcal{C}=\{\xi\in\mathbb{R}^{3}: \lambda r\leq |\xi|\leq \lambda  R\}$. For any $ g\in L^p$ and $\lambda > 0$, it holds that
\begin{equation*} \left\{
\begin{aligned}
{\rm{supp}}\, \mathcal{F}(g) \subset& \lambda \mathcal{B} \Rightarrow \|D^{m}g\|_{L^r}\lesssim\lambda ^{m+d(\frac{1}{p}-\frac{1}{r})}\|g\|_{L^p}, \nonumber\\
{\rm{supp}}\, \mathcal{F}(g) \subset& \lambda  \mathcal{C} \Rightarrow \lambda^{m}\|g\|_{L^{p}}\lesssim\|D^{m}g\|_{L^{p}}\lesssim \lambda^{m}\|g\|_{L^{p}}.
\end{aligned}\right.
\end{equation*}
 \end{lem}

By virtue of Lemma \ref{L2.1}, we present the following fundamental properties of Besov spaces.

\begin{lem}{\rm(\!\!\cite[Chapter 2]{BCD-Book-2011})}\label{L2.2}
The following properties hold:
\begin{itemize}
\item{} For $s\in\mathbb{R}$, $1\leq p_{1}\leq p_{2}\leq \infty$ and $1\leq r_{1}\leq r_{2}\leq \infty$, it holds that
\begin{equation}\notag
\begin{aligned}
\dot{B}^{s}_{p_{1},r_{1}}\hookrightarrow \dot{B}^{s-3 (\frac{1}{p_1} -\frac{1}{p_2})}_{p_{2},r_{2}}.
\end{aligned}
\end{equation}
\item{} For $1\leq p\leq r\leq\infty$, we have the following chain of continuous embeddings:
\begin{equation}\nonumber
\begin{aligned}
  \dot{B}^{0}_{p,1}\hookrightarrow L^{p}\hookrightarrow \dot{B}^{0}_{p,\infty}\hookrightarrow \dot B_{r,\infty}^{\varsigma},\quad \varsigma=-3\Big(\frac{1}{p}-\frac{1}{r}\Big).
\end{aligned}
\end{equation}
In addition, $\dot{B}^{\frac{3}{p}}_{p,1}$ is continuously embedded in the set of continuous functions decaying to 0 at infinity.

\item{} Let $s\in \mathbb{R}$ and $1\leq r\leq \infty$. Let $r^\prime$ denote the conjugate exponent of $r$. Then, the following duality estimates are valid:
\begin{align} \label{duality}
\langle g,h\rangle\lesssim \|g\|_{\dot B_{2,r}^s} \|h\|_{\dot B_{2,r}^{-s}}, \quad \text{and} \quad \|g\|_{\dot B_{2,r}^s}\lesssim \sup_{\phi\in \mathcal{S}}\langle g,\phi\rangle,
\end{align}
where the supremum is taken over all $\phi$ that satisfy $\|\phi\|_{\dot B_{2,r^\prime}^{-s}}\leq 1$ and $0 \notin {\rm Supp}\mathcal{F}\phi$.

\item{}  The following real interpolation properties are satisfied for $1\leq p,r\leq \infty$, $s_1<s_2$, and $\theta\in (0,1)$:
\begin{align}\notag
\|g\|_{\dot B_{p,r}^{\theta s_1+(1-\theta)s_2}}\lesssim\|g\|_{\dot B_{p,1}^{\theta s_1+(1-\theta)s_2}}\lesssim \frac{1}{\theta(1-\theta)(s_2-s_1)}\|g\|_{\dot B_{p,\infty}^{s_1}}^{\theta} \|g\|_{\dot B_{p,\infty}^{s_2}}^{1-\theta}  .  
\end{align}

\item{}  For any $\zeta>0$, we have the embeddings:
\begin{equation}\nonumber
\begin{aligned}
H^{s+\zeta}\hookrightarrow \dot{B}^{s}_{2,1}\hookrightarrow \dot{H}^{s}.
\end{aligned}
\end{equation}
\item{}
Let $\Lambda^{\sigma}$ be defined by $\Lambda^{\sigma}:=(-\Delta )^{\frac{\sigma}{2}}g:=\mathcal{F}^{-1}\big{(} |\xi|^{\sigma}\mathcal{F}(g) \big{)}$ for $\sigma\in \mathbb{R}$ and $g\in{\mathcal S}^{'}_{h}(\mathbb{R}^3)$, then $\Lambda^{\sigma}$ is an isomorphism from $\dot{B}^{s}_{p,r}$ to $\dot{B}^{s-\sigma}_{p,r}$.
\item{} Let $1\leq p_{1},p_{2},r_{1},r_{2}\leq \infty$, $s_{1}\in\mathbb{R}$ and $s_{2}\in\mathbb{R}$ satisfy
\begin{align*}
 s_{2}<\frac{3}{2}\quad\text{\text{or}}\quad s_{2}=\frac{3}{2}~\text{and}~r_{2}=1.
\end{align*}
   Then the  space $\dot{B}^{s_{1}}_{2,r_{1}}\cap \dot{B}^{s_{2}}_{2,r_{2}}$ endowed with the norm $\|\cdot \|_{\dot{B}^{s_{1}}_{2,r_{1}}}+\|\cdot\|_{\dot{B}^{s_{2}}_{2,r_{2}}}$ is a Banach space and has the weak compact and Fatou properties$:$ If $g_{n}$ is a uniformly bounded sequence of $\dot{B}^{s_{1}}_{2,r_{1}}\cap \dot{B}^{s_{2}}_{2,r_{2}}$, then an element $g$ of $\dot{B}^{s_{1}}_{2,r_{1}}\!\cap \dot{B}^{s_{2}}_{2,r_{2}}$ and a subsequence $g_{n_{k}}$ exist such that $g_{n_{k}}\rightarrow g$ in $\mathcal{S}'$ and
    \begin{align*} 
    \begin{aligned}
    \|g\|_{\dot{B}^{s_{1}}_{2,r_{1}}\cap \dot{B}^{s_{2}}_{2,r_{2}}}\lesssim \liminf_{n_{k}\rightarrow \infty} \|g_{n_{k}}\|_{\dot{B}^{s_{1}}_{2,r_{1}}\cap \dot{B}^{s_{2}}_{2,r_{2}}}.
    \end{aligned}
    \end{align*}
\end{itemize}
\end{lem}

In order to control the nonlinear terms in the KFP equation \eqref{G1.3}, we require the following product estimates lemmas.
\begin{lem}\label{L2.3}{\rm(\!\!\cite[Chapter 2]{BCD-Book-2011})}
The following statements hold:
\begin{itemize}
\item{} Let $s>0$, $1\leq  r\leq \infty$.  Then $\dot B^{s}_{2,r}\cap L^\infty$ is an algebra and
 \begin{align*} 
\|gh\|_{\dot B_{2,r}^s}\lesssim \|g\|_{L^\infty}\|h\|_{\dot B_{2,r}^s}+\|h\|_{L^\infty}\|g\|_{\dot B_{2,r}^s},    
 \end{align*}
which, in particular, when $s = \frac{3}{2}$ and $r = 1$, leads to
\begin{align*}
\|gh\|_{\dot B_{2,1}^{\frac{3}{2}}}\lesssim \|g\|_{\dot B_{2,1}^{\frac{3}{2}}}\|h\|_{\dot B_{2,1}^{\frac{3}{2 }}} .    
\end{align*}

\item{} Let the real numbers $s_1$, $s_2$  and $r$ fulfill
\begin{align*}
 1\leq  q,r_1,r_2\leq\infty, \quad s_1<\frac{3}{2},\quad s_2<\frac{3}{2},\quad s_1+s_2>0,\quad \frac{1}{r}=\frac{1}{r_1}+\frac{1}{r_2}.
\end{align*}
  Then, it holds  that
 \begin{align} \label{newinterpolation}
\|gh\|_{\dot B_{2,r}^{s_1+s_2-\frac{3}{2}}}\lesssim \|g\|_{ \dot B_{2,r_1}^{s_1 }}\|h\|_{\dot B_{2,r_2}^{s_2}} .    
 \end{align}

\item{} Let  the real numbers $s_1$ and $s_2$  fulfill
\begin{align*}
 s_1\leq \frac{3}{2},\quad s_2< \frac{3}{2},\quad s_1+s_2\geq 0.   
\end{align*}
Then, we have
 \begin{align} \label{inter}
\|gh\|_{\dot B_{2,\infty}^{s_1+s_2-\frac{3}{2}}}\lesssim \|g\|_{ \dot B_{2,1}^{s_1 }}\|h\|_{\dot B_{2,\infty}^{s_2}}.    
 \end{align}
    \end{itemize}
\end{lem}

Finally, we provide the Gagliardo–Nirenberg and Sobolev interpolation  inequalities.
\begin{lem} [{\!\cite[Lemma 2.1]{CDM-krm-2011}} and {\cite[Lemmas 2.1--2.2]{Dk-MZ-1992}}]\label{LA.1}   
For any $g,h\in H^N(\mathbb R^3) $ with $N\geq 3$ and any multi-index $\alpha$  with $1\leq|\alpha|\leq N$, it holds that
\begin{align*}
\|g\|_{L^{\infty } } \lesssim&\, \|\nabla  g\|_{L^{2} }^{\frac{1}{2}}
\|\nabla^{2}g\|_{L^{2}  }^{\frac{1}{2}}, \\
\|gh\|_{H^{N-1}  }  \lesssim&\, \|g\|_{H^{N-1}  }\|\nabla  h\|_{H^{N-1}  }, \\
\|\partial^{\alpha}(gh)\|_{L^{2}  }
\lesssim&\,\|\nabla  g\|_{H^{N-1}  }\|\nabla  h\|_{H^{N-1}  },\\
\|g\|_{L^6 } \lesssim&\,  \|\nabla  g\|_{L^2 }\lesssim  \|g\|_{H^1  },\\
  \|g\|_{L^q  } \lesssim&\, \|g\|_{H^1 }, \quad 2\leq q\leq 6.
\end{align*}
\end{lem}

\subsection{Macro-micro decomposition  and some properties of $\mathcal{L}$}
In this subsection, we initially state the macro-micro decomposition, which was originally introduced in \cite{Gy-IUMJ-2004} for the Boltzmann equation and further developed in \cite{DFT-2010-CMP} for a KFP equation. It allows us to express $f(t, x, v)$ as the sum of the macroscopic component $\mathbf{P}f$ and the microscopic part $\{\mathbf{I}-\mathbf{P}\}f$, namely,
\begin{align}\label{G2.5}
f=\mathbf{P}f+\{\mathbf{I}-\mathbf{P}\}f.    
\end{align}
Here, $\mathbf{P}$ represents the orthogonal projection operator, which is defined as follows: 
\begin{equation*}  \mathbf{P}:  L_{v}^{2}(\mathbb{R}^3)\rightarrow \mathcal{N}:=\mathrm{Span}\{\sqrt{M}, v_1 \sqrt{M},  v_2\sqrt{M},v_3\sqrt{M}\}.
\end{equation*} 
This operator can be decomposed into
\begin{align*}  \mathbf{P} := \mathbf{P}_{0} \oplus \mathbf{P}_{1},  
\end{align*}
where $\mathbf{P}_{0}f := a \sqrt{M}$ and $\mathbf{P}_{1}f := b \cdot v \sqrt{M}$.
Moreover, $a = a^f$ and $b = b^f$ represent the  macroscopic components of $\mathbf{P}f$, and they are defined by 
\begin{align*}  a = a^f(t, x) = \int_{\mathbb{R}^3} \sqrt{M} f(t, x, v) \, \mathrm{d}v, \quad  b = b^f(t, x) = \int_{\mathbb{R}^3} v \sqrt{M} f(t, x, v) \, \mathrm{d}v.  \end{align*}
Notice that according to the definitions of $a$ and $b$,   \begin{align*}
\mathbf{P}f \perp \{\mathbf{I}-\mathbf{P}\}f   
\end{align*} 
holds in $L_v^2$ for any $(t,x)$.

As stated in \cite{AMTU-CPDE-2001}, the linearized Fokker-Planck operator $\mathcal{L}$ demonstrates certain dissipative properties. Similar to  \cite{DFT-2010-CMP,CDM-krm-2011},   the following  properties of $\mathcal{L}$ hold. For the sake of brevity, we omit the proofs here.
\begin{lem}\label{L2.4}
The linearized Fokker-Planck operator $\mathcal{L}$ is self-adjoint and has the following properties:
\begin{itemize}
    \item  $\mathcal{L}$ is a non-positive operator satisfying
    \begin{align*}
    \langle \mathcal{L} f, f \rangle=-\int_{\mathbb R^3} \bigg| \nabla_{v} \bigg(   \frac{f}{\sqrt{M}}     \bigg)
    \bigg |^2 M{\rm d}v  \leq 0,
    \end{align*}
and
\begin{align*}
 \ker \mathcal{L}= {\rm Span}\big\{ \sqrt{M} \big\},  \quad {\rm Rang} \mathcal{L}={\rm Span} \big\{ \sqrt{M} \big\}^\perp .
\end{align*}
    \item 
There is a constant $\bar\lambda >0$   such that
\begin{align*}
-\langle \mathcal{L}f, f\rangle \geq \bar\lambda |\{\mathbf{I}-\mathbf{P}_0\}f|_{\nu}^2.
\end{align*}

\item
There is a constant  $\lambda_0>0$ such that
\begin{align}\label{G2.6}
-\langle \mathcal{L}f, f\rangle \geq  \lambda_0 |\{\mathbf{I}-\mathbf{P}\}f|_{\nu}^2+b^2.
\end{align}
\end{itemize}
\end{lem}

\subsection{Decay estimates of semi-group}
To capture the solution behavior of the spectral properties of the KFP equation, we need to rewrite \eqref{G1.3} in the form of an evolution equation:
\begin{align}\label{B.1}
\partial_t f=\mathcal{B}f +E\cdot\Big(v\sqrt{M}-\nabla_{v}f+\frac{1}{2}vf    \Big)+S(t), \quad (t,x,v)\in \mathbb R^{+}\times \mathbb R^3_{x}\times \mathbb R^3_{v},
\end{align}
where $\mathcal B$ is the   operator given by
\begin{align}\label{def.B.add}
\mathcal{B}=-v\cdot\nabla_x+\mathcal{L}.    
\end{align}
By ignoring the force and source terms in \eqref{B.1}, following \cite{LY-AA-2015}, we consider the following linear Cauchy problem for the KFP equation:
 \begin{equation}\label{B.2}
\left\{\begin{aligned}
&  \partial_t f=\mathcal{B}f, \quad (t,x,v)\in \mathbb R\times \mathbb R^3_{x}\times \mathbb R^3_{v},\\ 
& f(0,x,v)=f_0(x,v), 
\quad \quad\,   (x,v)\in    \mathbb R^3_{x}\times \mathbb R^3_{v},
 \end{aligned}
 \right.
\end{equation} 
and have the following result.

\begin{lem}\label{LB.1}{\rm(\!\!\cite{LY-AA-2015})}
For any $t \geq 0$ and $\xi\in\mathbb{R}^3$, the solution  $f(t,x,v)$ to \eqref{B.2} satisfies
\begin{align*}
\|\widehat{f}(t,\xi,v)\|_{L_{v}^2}\lesssim   e^{{-\lambda_{1}t\min\{1,|\xi|^2\}}}  {}\|\widehat{f_0}(\xi,v)\|_{L_v^2},
\end{align*}
where $\lambda_1 > 0$ is a uniform constant.
\end{lem}

Notice that Lemma \ref{LB.1} cannot be applied directly to overcome the difficulties arising from nonlinear terms that contain weighted $v$ or the derivative $\nabla_{v}$, such as $\frac{1}{2}E\cdot v f$ and $E\cdot \nabla_{v}f$.
To circumvent the   difficulties arising from the weighted $v$ and derivative $\nabla_{v}$, inspired by    the Fourier analysis techniques developed in \cite[Theorem 3.1]{CDM-krm-2011}, we consider the following  Cauchy problem of the linearized KFP equation:
 \begin{equation}\label{B.3}
\left\{\begin{aligned}
&  \partial_t f=\mathcal{B}f+S_{f}, \quad (t,x,v)\in \mathbb R\times \mathbb R^3_{x}\times \mathbb R^3_{v},\\ 
& f(0,x,v)=f_0(x,v), 
\quad \quad\,   (x,v)\in    \mathbb R^3_{x}\times \mathbb R^3_{v},
 \end{aligned}
 \right.
\end{equation} 
where the additional inhomogeneous source $S_{f}$ takes the form
\begin{align*}
 S_{f}:=\nabla_{v}\cdot G-\frac{1}{2}v\cdot G+h.   
\end{align*}
Here, $ G = G(t, x, v) \in \mathbb{R}^3 $ and $ h = h(t, x, v) \in \mathbb{R} $ satisfy  
\begin{align*}  
\mathbf{P}_0 G = 0 \quad \text{and} \quad \mathbf{P} h = 0,  
\end{align*} 
under which one can verify that $\mathbf{P}S_{f}=0$.
According to Duhamel’s principle, the solution  of \eqref{B.3}  can be expressed as  
\begin{align} \label{B.4} 
f(t) = e^{t\mathcal{B}}f_0 + \int_0^t e^{t\mathcal{B}}S_f(t-s)\, {\rm d}s , 
\end{align}  

\begin{rem}\label{remB.1}
If we set $S_{f}= -E\cdot\nabla_{v}f+\frac{1}{2}E\cdot v f$, then we have the following decomposition of $S_{f}$:
\begin{align*}
S_{f}=&\,- E\cdot\nabla_{v}\{\mathbf{I}-\mathbf{P}_0\}f+\frac{1}{2}E\cdot v   \{\mathbf{I}-\mathbf{P}_0\}f- E\cdot\nabla_{v}  \mathbf{P}_0 f+\frac{1}{2}E\cdot v     \mathbf{P}_0 f\nonumber\\
=&\,\nabla_{v}\cdot G-\frac{1}{2}v\cdot G+E\cdot a v\sqrt{M},
\end{align*}
where 
\begin{align*}
G:=-E\{\mathbf{I}-\mathbf{P}_0\}f.    
\end{align*}
Thus, with the estimates of $G$, we are able to handle the nonlinear terms $E\cdot\nabla_{v}f $ and $ \frac{1}{2}v\cdot f$.
\end{rem}

Motivated by the hyperbolic/hyperbolic–parabolic argument \cite{SK-HMJ-1985}, we are able to establish sharp point-wise decay estimates of the solution to the   problem \eqref{B.3}.

\begin{lem}\label{LB.2}
For any $t \geq 0$ and $\xi\in\mathbb{R}^3$, the solution  $f(t,x,v)$ to \eqref{B.3} satisfies
\begin{align}\label{B.5}
\|\widehat{f}(t,\xi,v)\|_{L_{v}^2}^2\lesssim&\,   e^{-\lambda_{2}t\min\{1,|\xi|^2\}}  {}\|\widehat{f_0}(\xi,v)\|_{L_v^2}^2\nonumber\\
&+\int_0^t  e^{-\lambda_{2}(t-\tau)\min\{1,|\xi|^2\}}\big(\|\widehat{ G}(\tau)\|_{L_v^2}^2+\|\nu^{-\frac{1}{2}}\widehat{ h}(\tau)\|_{L_v^2}^2  \big) {\rm d}\tau,
\end{align}
where $\lambda_2 > 0$ is a uniform constant.
\end{lem}

\begin{rem}
The sharp point-wise decay estimates of solutions also assist us in comprehending that the solutions of   \eqref{B.3} display the diffusion effect at low frequencies and the spectrum gap at high frequencies. This mechanism is termed   hypocoercivity in \cite{Villani-MAMS-2009,D-2011-Nonlinearity}.
\end{rem}

\begin{proof}[Proof of Lemma \ref{LB.2}]
By applying the Fourier transform with respect to $x$ to the equation in \eqref{B.3}, one has
\begin{align}\label{B.6}
\partial_{t}\widehat{f}+{\rm i}v\cdot \xi \widehat{f} -\mathcal{L}\widehat{f}=\nabla_v\cdot\widehat{G}-\frac{1}{2}v\cdot\widehat{G}+\widehat{h}.
\end{align}
Taking the inner product of \eqref{B.6} with $\overline{\widehat{f}}$ and integrating over $v$, the real part of the resulting expression yields
\begin{align*}
&\frac{{\rm d}}{{\rm d}t} \|\widehat{f}\,\|_{L_v^2}^2+2|\widehat{b}|^2+2\mathbf{Re}\int_{\mathbb{R}^3} 
(-\mathcal{L}\{\mathbf{I}-\mathbf{P}\}\widehat{f}\,|\{\mathbf{I}-\mathbf{P}\}\widehat{f}\,\mathrm{d}v   )\nonumber\\
 &\quad =2\mathbf{Re}\int_{\mathbb{R}^3}\big(\nabla_{v}\cdot \widehat{G}-\frac{1}{2}v\cdot\widehat{G}|
\{\mathbf{I}-\mathbf{P}\}\widehat{f}\,\big)\mathrm{d}v
+2\mathbf{Re}\int_{\mathbb{R}^3}(\widehat{h}|\{\mathbf{I}-\mathbf{P}\}\widehat{f}\,)\mathrm{d}v,
\end{align*}
where we have used the following facts: 
\begin{align*}
\mathbf{P}h=0,\quad \nabla_{v}\cdot  G-\frac{1}{2}v\cdot G\perp {\rm Rang}\mathbf{P}.
\end{align*}
Leveraging the coercivity of $-\mathcal{L}$,   applying \eqref{G2.6} in Lemma \ref{L2.4},  integration by parts, and  Cauchy-Schwarz's inequality, we further  get
\begin{align}\label{B.7}
\frac{{\rm d}}{{\rm d}t} \|\widehat{f}\,\|_{L_v^2}^2 +2 |\widehat{b}|^2 + {\lambda_0} |\{\mathbf{I}-\mathbf{P}\}\widehat{f}\,|_{\nu}^2  
\lesssim \|\widehat{G}\|_{L^2_v}^2 + \|\nu^{-\frac{1}{2}}\widehat{h}\|_{L^2_v}^2.
\end{align}
Since we already have the dissipation estimate for $b$,  we only need to further consider the corresponding dissipation estimate for $a$ as follows. Due to $\mathbf{P}S_f = 0$, it follows from \eqref{B.3}$_1$ that
\begin{equation}\label{B.8}
\left\{\begin{aligned}
&\partial_{t} {a}+ \nabla \cdot {b}=0,\\
&\partial_{t}  {b}_i+  \partial_{i}{a} +\sum_{j=1}^3\partial_j\Gamma_{ij}(\{\mathbf{I}-\mathbf{P}\} {f}\,)  + {b}_i=0,
 \end{aligned}
 \right.
\end{equation}
where 
\begin{align} \label{gammaa}
\Gamma_{ij}(\cdot):={\big\langle (v_iv_j-1)\sqrt{M},\,\cdot\,   \big\rangle} 
\end{align}
is the moment functional. Note that the equations \eqref{B.8}$_1$ and \eqref{B.8}$_2$ can be obtained by respectively multiplying
\eqref{B.3} by $\sqrt{M}$ and $v_i\sqrt{M}$ with $1\leq i \leq3$ and then integrating the results over $\mathbb R^3_{v}$.

Taking the Fourier transform in $x$ to \eqref{B.8} yields
\begin{equation}\label{B.9}
\left\{\begin{aligned}
&\partial_{t} \widehat {a}+  {\rm i}\xi \cdot \widehat  {b}=0,\\
&\partial_{t}  \widehat{b}_i+   {\rm i}\xi_{i}\widehat{a} +\sum_{j=1}^3 {\rm i}  \xi_j\Gamma_{ij}(\{\mathbf{I}-\mathbf{P}\} \widehat{f}\,)  + \widehat{b}_i=0.
 \end{aligned}
 \right.
\end{equation}
By applying the basic energy method argument for \eqref{B.9},  we   deduce that
\begin{align}\label{B.10}
\frac{{\rm d}}{{\rm d}t} \mathbf{Re}\mathcal{E}_1(\widehat{f}\,)+\frac{ \lambda_3 |\xi|^2}{1+|\xi|^2}  |\widehat{a
}|^2 \lesssim  \|\{\mathbf{I}-\mathbf{P}\}\widehat{f}\,\|_{L_v^2}^2+| \widehat{b}|^2,   
\end{align}
for some constant $\lambda_3>0$, where 
\begin{align*}
\mathcal{E}_1(\widehat{f}\,):= -\frac{1}{1+|\xi|^2}(\widehat{a}|{\rm i}\xi\cdot\widehat{b}).   
\end{align*}

Define the   energy functional $\mathcal{E}_{\mathcal{F}}(\widehat{f} )$ by
\begin{align*}
\mathcal{E}_{\mathcal{F}}(\widehat{f} ):=\|f\|_{L^2_v}^2  +\tau_1\mathbf{Re}\mathcal{E}_1(\widehat{f}\,)  ,
\end{align*}
where $0<\tau_1 \ll1$ is a sufficiently small constant. By Young's inequality, we have
\begin{align*}
    |\mathcal{E}_1(\widehat{f}\,)|\leq \frac{|\xi|^2}{2(|\xi|^2+1)} (|\widehat{a}|^2+|\widehat{b}|^2) \leq \frac{1}{2} \|\widehat{f}\,\|_{L_v^2}^2,
\end{align*}
which implies that
\begin{align}\label{B.11}
\mathcal{E}_{\mathcal{F}}(\widehat{f} )\backsim \|\widehat{f}\,\|_{L_v^2}^2.
\end{align}
The sum of Equation \eqref{B.7} and $\tau_1\times $   \eqref{B.10} leads to
\begin{align*} 
\partial_t\mathcal{E}_{\mathcal{F}}(\widehat{f} )+\frac{\lambda_{4}|\xi|^2}{1+|\xi|^2}\mathcal{E}_{\mathcal{F}}(\widehat{f} )\lesssim  \|\widehat{G}\|_{L^2_v}^2+\|\nu^{-\frac{1}{2}}\widehat{h}\|_{L^2_v}^2 , 
\end{align*}
for some constant $\lambda_{4}>0$.

According to Gr\"{o}nwall's inequality, it follows that
\begin{align}\label{B.12}
 \mathcal{E}_{\mathcal{F}}(\widehat{f} )\leq e^{-\frac{\lambda_{4}|\xi|^2}{1+|\xi|^2}}  \mathcal{E}_{\mathcal{F}}(\widehat{f}_0)+\int_{0}^te^{-\frac{\lambda_{4}|\xi|^2}{1+|\xi|^2}(t-\tau)}\big(\|\widehat{G}(\tau)\|_{L^2_v}^2+\|\nu^{-\frac{1}{2}}\widehat{h}(\tau)\|_{L^2_v}^2\big)\mathrm{d}\tau.
\end{align}
As shown in \cite{Ks-1983}, with the help of \eqref{B.11},    the desired estimate \eqref{B.5} can be derived by applying \eqref{B.12} and then setting the inhomogeneous source $S _f = 0$ and the initial data $f_0 = 0$, respectively. For brevity, the detailed proof is omitted here.
\end{proof}

Recall \eqref{def.B.add} and \eqref{B.4}. Thanks to the time-decay estimates  in Lemma \ref{LB.2}, and following the approach of \cite[Lemma 3.3]{Deguchi-2025}, we can further establish the subsequent $L^2$-type decay rates in both the low-frequency and high-frequency regions in the homogeneous Besov space as follows. 

\begin{lem}\label{LB.3}
Let $t\geq 0$ and $1\leq r\leq \infty$. Let $s,s_0\in \mathbb{R}$ such that $s_0\leq s$. Then, there exists $j_0\in \mathbb{Z}$ such that, for any $1\leq r\leq\infty$ and $f_0\in L_v^2(\dot B_{2,r}^{s_0})$, the following inequality holds:
\begin{align}\label{B.13}
\|\dot S_{j_0}e^{t\mathcal{B}} f_0\|_{L_v^2(\dot B_{2,r}^s)}\lesssim (1+t)^{-\frac{s - s_0}{2}} \|f_0\|_{L_v^2(\dot B_{2,r}^{s_0})}.
\end{align}
For any $1\leq r\leq\infty$ and $f_0\in L_v^2(\dot B_{2,r}^{s})$, we have
\begin{align}\label{B.14}
\|({\rm Id} - \dot S_{j_0})e^{t\mathcal{B}} f_0\|_{L_v^2(\dot B_{2,r}^s)}\lesssim e^{-\lambda_5t} \|f_0\|_{L_v^2(\dot B_{2,r}^{s})},
\end{align}
where $\lambda_5>0$ is a constant.
For any $1\leq r\leq\infty$ and $f_0\in L_v^2(\dot B_{2,r}^{s-\frac{2}{r}})$, the following inequality is valid:
\begin{align}\label{B.15}
\|\dot S_{j_0}e^{t\mathcal{B}}f_0\|_{L^r(L_v^2(\dot B_{2,r}^s))}\lesssim \|f_0\|_{L_v^2(\dot B_{2,r}^{s-\frac{2}{r}})}.
\end{align}
The estimates \eqref{B.13}, \eqref{B.14} and \eqref{B.15} still hold when $\mathcal{B}$ is replaced by its adjoint operator $\mathcal B^*=v\cdot\nabla_{x}+\mathcal{L}$.

In addition, for the inhomogeneous part,   it holds that
\begin{align}\label{B.23}
\bigg\|\int_0^t \dot S_{j_0}e^{(t-\tau)\mathcal{B}} S_{f}(\tau){\rm d}\tau\bigg\|_{L_v^2(\dot B_{2,r}^s)}^2 \lesssim \int_0^t (1+t-\tau)^{-{(s - s_0)}} \Big(\|G (\tau)\|_{L_v^2(\dot B_{2,r}^{s_0})}^2+\|\nu^{-\frac{1}{2}}h(\tau) \|_{L_v^2(\dot B_{2,r}^{s_0})}^2 \Big) {\rm d}\tau,
\end{align}
for any $2\leq r\leq \infty$ and $G(\tau),\nu^{-\frac{1}{2}}h(\tau) \in L_v^2(\dot B_{2,r}^{s_0})$, 
and
\begin{align}\label{B.24}
\bigg\|\int_0^t \dot ({\rm Id}-S_{j_0})e^{(t-\tau)\mathcal{B}} S_{f}(\tau){\rm d}\tau\bigg\|_{L_v^2(\dot B_{2,r}^s)}^2 \lesssim \int_0^t  e^{-\lambda_5(t-\tau)} \Big(\|G (\tau)\|_{L_v^2(\dot B_{2,r}^{s})}^2+\|\nu^{-\frac{1}{2}}h(\tau) \|_{L_v^2(\dot B_{2,r}^{s})}^2 \Big) {\rm d}\tau,
\end{align}
for any $2\leq r\leq \infty$ and $G(\tau),\nu^{-\frac{1}{2}}h(\tau) \in L_v^2(\dot B_{2,r}^{s })$. 
\end{lem}

\begin{proof}
The proofs of \eqref{B.13}--\eqref{B.15} in Lemma \ref{LB.3} can be obtained in the same way as those in 
\cite[Lemma 3.3]{Deguchi-2025} and \cite[Lemma 4.1]{Deguchi-2024-Math Ann} on compressible Navier-Stokes equations. For brevity, we omit them here and   interested {readers}  refer to \cite{Deguchi-2025,Deguchi-2024-Math Ann} for more   details.

Now we give the proof of \eqref{B.23} and \eqref{B.24}. We divide it into two parts for low and high frequencies respectively.
For the case of low frequency, we set
\begin{align*}
\mathcal{R}_L(t):=\int_0^t \dot S_{j_0}e^{(t-\tau)\mathcal B}S_f(\tau){\rm d}\tau.    
\end{align*}
From Lemma \ref{LB.2}, we have
\begin{align*}
\|\dot\Delta_j \mathcal{R}_L(t)\|_{L^2_{x,v}}^2
\lesssim
\int_0^t e^{-c2^{2j}(t-\tau)}
\Big(
\|\dot\Delta_j G(\tau)\|_{L^2_{x,v}}^2
+
\|\nu^{-\frac12}\dot\Delta_j h(\tau)\|_{L^2_{x,v}}^2
\Big){\rm d}\tau,  
\end{align*}
for each $j \leq j_0$ and a certain constant $c > 0$.
Multiplying both sides by $2^{2js}$, we obtain
\begin{align*}
2^{2js}\|\dot\Delta_j \mathcal{R}_L(t)\|_{L^2_{x,v}}^2
\lesssim&\,
\int_0^t
2^{2j(s-s_0)}e^{-c2^{2j}(t-\tau)}
\Big(
2^{2js_0}\|\dot\Delta_j G(\tau)\|_{L^2_{x,v}}^2+
2^{2js_0}\|\nu^{-\frac12}\dot\Delta_j h(\tau)\|_{L^2_{x,v}}^2
\Big){\rm d}\tau \nonumber\\
\lesssim&\, \int_0^t (1+t-\tau)^{-(s-s_0)}
\Big(
2^{2js_0}\|\dot\Delta_j G(\tau)\|_{L^2_{x,v}}^2
+
2^{2js_0}\|\nu^{-\frac12}\dot\Delta_j h(\tau)\|_{L^2_{x,v}}^2
\Big){\rm d}\tau,
\end{align*}
where we have used the following standard low-frequency bound:
\begin{equation*} 
2^{2j(s-s_0)}e^{-c2^{2j}(t-\tau)}
\lesssim
(1+t-\tau)^{-(s-s_0)},\quad \forall\, j\le j_0 .
\end{equation*}

Now, we set
\begin{align*}
A_j(t):=&\,2^{2js}\|\dot\Delta_j \mathcal R_L(t)\|_{L^2_{x,v}}^2,\\
B_j(\tau):=&\,2^{2js_0}\|\dot\Delta_j G(\tau)\|_{L^2_{x,v}}^2, \\H_j(\tau):=&\,2^{2js_0}\|\nu^{-\frac12}\dot\Delta_j h(\tau)\|_{L^2_{x,v}}^2 .   
\end{align*}
Then, for all $j\le j_0$, it holds that
\begin{align*}
A_j(t)
\lesssim
\int_0^t (1+t-\tau)^{-(s-s_0)}\big(B_j(\tau)+H_j(\tau)\big){\rm d}\tau .   
\end{align*}
 
Since $r\geq2$, we have $\frac{r}{2}\geq1$, and thus $\ell^{\frac{r}{2}}$ is a Banach space. By applying Minkowski's inequality in $\ell^{\frac{r}{2}}$, we  infer that
\begin{align*}
 \|(A_j(t))_{j\le j_0} \|_{\ell^{\frac{r}{2}}}
\lesssim&\,
\int_0^t (1+t-\tau)^{-(s-s_0)}
\|(B_j(\tau)+H_j(\tau))_{j\le j_0}\|_{\ell^{\frac{r}{2}}} {\rm d}\tau
\\
\lesssim&\,
\int_0^t (1+t-\tau)^{-(s-s_0)}
\Big(
\|(B_j(\tau))_{j\le j_0}\|_{\ell^{\frac{r}{2}}}
+
\|(H_j(\tau))_{j\le j_0}\|_{\ell^{\frac{r}{2}}}
\Big){\rm d}\tau .
\end{align*}
Observe that
\begin{align*}
\|(A_j(t))_{j\le j_0}\|_{\ell^{\frac{r}{2}}}
=
\Big\|\big(2^{js}\|\dot\Delta_j \mathcal{R}_L(t)\|_{L^2_{x,v}}\big)_{j\le j_0}\Big\|_{\ell^r}^{\,2}
=
\|\mathcal R_L(t)\|_{L_v^2(\dot B_{2,r}^{s})}^{2},    
\end{align*}
and similarly,
\begin{align*}
\|(B_j(\tau))_{j\le j_0}\|_{\ell^{\frac{r}{2}}}
\leq&\, 
\|G(\tau)\|_{L_v^2(\dot B_{2,r}^{s_0})}^{2},\\
\|(H_j(\tau))_{j\le j_0}\|_{\ell^{\frac{r}{2}}}
\leq&\, 
\|\nu^{-\frac12}h(\tau)\|_{L_v^2(\dot B_{2,r}^{s_0})}^{2}.    
\end{align*} 
Combining the above estimates yields \eqref{B.23}.

For the case of high frequency, we set
\begin{align*}
\mathcal R_H(t):=\int_0^t ({\rm Id}-\dot S_{j_0})e^{(t-\tau)\mathcal B}S_f(\tau) {\rm d}\tau .    
\end{align*}
By Lemma \ref{LB.2}, for every $j\ge j_0-1$, it holds 
\begin{align*}
\|\dot\Delta_j \mathcal R_H(t)\|_{L^2_{x,v}}^2
\lesssim
\int_0^t e^{-\lambda_2 (t-\tau)}
\Big(
\|\dot\Delta_j G(\tau)\|_{L^2_{x,v}}^2
+
\|\nu^{-\frac12}\dot\Delta_j h(\tau)\|_{L^2_{x,v}}^2
\Big){\rm d}\tau .    
\end{align*}
Multiplying both sides of the above inequality by $2^{2js}$, we obtain
\begin{align*}
2^{2js}\|\dot\Delta_j \mathcal R_H(t)\|_{L^2_{x,v}}^2
\lesssim
\int_0^t e^{-\lambda_2 (t-\tau)}
\Big(
2^{2js}\|\dot\Delta_j G(\tau)\|_{L^2_{x,v}}^2
+
2^{2js}\|\nu^{-\frac12}\dot\Delta_j h(\tau)\|_{L^2_{x,v}}^2
\Big){\rm d}\tau .
\end{align*}
Define
\begin{align*}
\widetilde A_j(t):=&\,2^{2js}\|\dot\Delta_j \mathcal R_H(t)\|_{L^2_{x,v}}^2    \\
\widetilde B_j(\tau):=&\,2^{2js}\|\dot\Delta_j G(\tau)\|_{L^2_{x,v}}^2, \\
\widetilde H_j(\tau):=&\,2^{2js}\|\nu^{-\frac12}\dot\Delta_j h(\tau)\|_{L^2_{x,v}}^2 .
\end{align*}
Then, for $j\ge j_0-1$, it yields
\begin{align*}
\widetilde A_j(t)
\lesssim
\int_0^t e^{-\lambda_2 (t-\tau)}
\big(\widetilde B_j(\tau)+\widetilde H_j(\tau)\big){\rm d}\tau .    
\end{align*}

Since $r\ge2$, we  again apply Minkowski's inequality in $\ell^{\frac{r}{2}}$-norm to deduce
\begin{align*}
\|(\widetilde A_j(t))_{j\ge j_0-1}\|_{\ell^{\frac{r}{2}}}
\lesssim&\,
\int_0^t e^{-\lambda_2 (t-\tau)}
\|(\widetilde B_j(\tau)+\widetilde H_j(\tau))_{j\ge j_0-1}\|_{\ell^{\frac{r}{2}}}{\rm d}\tau
\\
\lesssim&\,
\int_0^t e^{-\lambda_2 (t-\tau)}
\Big(
\|(\widetilde B_j(\tau))_{j\ge j_0-1}\|_{\ell^{\frac{r}{2}}}
+
\|(\widetilde H_j(\tau))_{j\ge j_0-1}\|_{\ell^{\frac{r}{2}}}
\Big){\rm d}\tau .
\end{align*}
By the definition of the   norm in the Besov space $ {L_v^2(\dot B_{2,r}^{s})} $, we have
\begin{align*}
\|(\widetilde A_j(t))_{j\ge j_0-1}\|_{\ell^{\frac{r}{2}}}
=
\|\mathcal R_H(t)\|_{L_v^2(\dot B_{2,r}^{s})}^{2},    
\end{align*}
while
\begin{align*}
\|(\widetilde B_j(\tau))_{j\ge j_0-1}\|_{\ell^{\frac{r}{2}}}
\leq&\, 
\|G(\tau)\|_{L_v^2(\dot B_{2,r}^{s})}^{2}, \\
\|(\widetilde H_j(\tau))_{j\ge j_0-1}\|_{\ell^{\frac{r}{2}}}
\leq&\,
\|\nu^{-\frac12}h(\tau)\|_{L_v^2(\dot B_{2,r}^{s})}^{2}.
\end{align*}
Therefore,
\begin{align*}
\|\mathcal R_H(t)\|_{L_v^2(\dot B_{2,r}^{s})}^{2}
\lesssim
\int_0^t e^{-\lambda_2 (t-\tau)}
\Big(
\|G(\tau)\|_{L_v^2(\dot B_{2,r}^{s})}^{2}
+
\|\nu^{-\frac12}h(\tau)\|_{L_v^2(\dot B_{2,r}^{s})}^{2}
\Big){\rm d}\tau ,    
\end{align*}
which completes the proof of \eqref{B.24}.
\end{proof}

\begin{rem}
The range of values of $r$ for the inhomogeneous part is $2 \leq r \leq \infty$, while the range of values of $r$ for the homogeneous part is $1 \leq r \leq \infty$. This difference arises because the inhomogeneous part is square-integrable.
\end{rem}

\section{Global existence of classical solution to the KFP equation}\label{sec3}
In this section, we prove the existence and uniqueness of the global solution to the Cauchy problem \eqref{G1.4}. To this end, we aim to derive the {\it uniform a priori estimates} of $f(t,x,v)$. To begin with, we assume that $f(t,x,v)$ is a classical solution to the Cauchy problem \eqref{G1.4} defined on the time interval $0 \leq t < T_1$, where $T_1 > 0$ can be arbitrary. It satisfies the following inequality:
\begin{align}\label{G3.1}
\sup_{{0\leq t\leq T_1}} \left\{\|f(t)\|_{L_v^2(\dot B_{2,\infty}^{\frac{1}{2}}\cap\dot H^N)}+\|E(t)\|_{\dot B_{2,\infty}^{-\frac{3}{2}}\cap \dot H^N}+\|S (t)\|_{L_v^2(\dot B_{2,\infty}^{-\frac{3}{2}}\cap \dot H^N)}\right\}\leq \delta,
\end{align}
where $N\geq3$ is an integer and $0 < \delta < 1$ is a sufficiently small generic constant.

We decompose $f(t,x,v)$ into low-frequency and high-frequency components, such that 
\begin{align*}
f(t,x,v) = f_L(t,x,v) + f_H(t,x,v),    
\end{align*}
where $f_L(t,x,v) = \dot{S}_{j_0}f(t,x,v)$. First, we present an estimate of $f_L(t,x,v)$ in the homogeneous Besov space $\dot{B}_{2,\infty}^{\frac{1}{2}}$.

\subsection{Estimates of the solution at low frequencies}
\begin{lem}\label{L3.1}
For the classical solution $f$ to the problem   \eqref{G1.4}, it holds that
\begin{align}\label{G3.2}
 \|f_{L}(t)\|_{L_v^2(\dot B_{2,\infty}^{\frac{1 }{2}})}  \lesssim&\, \sup_{{0\leq t\leq T_1}}\|E(t)\|_{\dot B_{2,\infty}^{-\frac{3}{2}}\cap\dot B^{-\frac{1}{2}}_{2,\infty}} + \sup_{{0\leq t\leq T_1}}\|S(t)\|_{L_v^2(\dot B_{2,\infty}^{-\frac{3}{2}})}+\|f
_0\|_{L_v^2(\dot B_{2,\infty}^{\frac{1 }{2}})},
\end{align}
for any $0\leq t \leq T_1$.
\end{lem}
\begin{proof}

Applying the operator $\dot{S}_{j_0}$ to \eqref{G1.4} and using the Duhamel principle yields  
\begin{align}\label{G3.3}
f_{L}(t)= e_{L}^{t\mathcal B}f_{0}+\int_0^t e_L^{(t-\tau)\mathcal{B}} \Big[ E\cdot\Big( \frac{1}{2}vf-\nabla_{v}f \Big) (\tau)+E\cdot v\sqrt{M}(\tau)+ S(\tau)    \Big]{\rm d}\tau, 
\end{align}
where $e_{L}^{t\mathcal{B}}:=\dot S_{j_0}e^{t\mathcal{B}}$ and
\begin{align*}
 S_{f}:= E\cdot\Big( \frac{1}{2}vf-\nabla_{v}f \Big) =\nabla_{v}\cdot G-\frac{1}{2}v\cdot G+E\cdot av \sqrt{M} ,
\end{align*}
with
\begin{align*}
G:= -E\{\mathbf{I}-\mathbf{P}_0\}f .   
\end{align*}
Taking the $\dot B_{2,\infty}^{\frac{1}{2}}$-norm of \eqref{G3.3} and applying the triangle inequality, we obtain
\begin{align} \label{G3.4}
\|f_{L}(t)\|_{L_v^2(\dot B_{2,\infty}^{\frac{1}{2}})}\lesssim&\, \|e_{L}^{t\mathcal{B}} f_0\| _{L_v^2(\dot B_{2,\infty}^{\frac{1}{2}})} +\bigg\|\int_0^t  e_L^{(t-\tau)\mathcal{B}} \big(   E\cdot(1+a) v\sqrt{M}+ S     \big)(\tau){\rm d}\tau     \bigg\|_{L_v^2(\dot B_{2,\infty}^{\frac{1}{2}})}  \nonumber\\
&+\bigg\|\int_0^t  e_L^{(t-\tau)\mathcal{B}} \Big(  \nabla_{v}\cdot G-\frac{1}{2}v\cdot G      \Big)(\tau){\rm d}\tau     \bigg\|_{L_v^2(\dot B_{2,\infty}^{\frac{1}{2}})} \nonumber\\
\lesssim&\,\|f_0\|_{L_v^2(\dot B_{2,\infty}^{\frac{1}{2}})}+ \bigg\|\int_0^t  e_L^{(t-\tau)\mathcal{B}} \big(   E\cdot(1+a) v\sqrt{M}+ S     \big)(\tau){\rm d}\tau     \bigg\|_{L_v^2(\dot B_{2,\infty}^{\frac{1}{2}})}\nonumber\\
&+\bigg\|\int_0^t  e_L^{(t-\tau)\mathcal{B}} \Big(  \nabla_{v}\cdot G-\frac{1}{2}v\cdot G      \Big)(\tau){\rm d}\tau     \bigg\|_{L_v^2(\dot B_{2,\infty}^{\frac{1}{2}})}\nonumber\\
\equiv:&\, \|f_0\|_{L_v^2(\dot B_{2,\infty}^{\frac{1}{2}})}+I_1+I_2.
\end{align}
For any test function $\phi\in L_v^2(\mathcal{S})$, owing to   \eqref{duality} in  Lemma \ref{L2.2}, we have
\begin{align}\label{G3.5}
&\bigg\langle \int_0^t e^{(t-\tau)\mathcal{B}}_{L} \big(   E\cdot(1+a) v\sqrt{M}+ S     \big)(\tau) {\rm d}\tau, \phi      \bigg\rangle_{x,v}  \nonumber\\
=&\, \int_0^t \big \langle  e^{(t-\tau)\mathcal{B}}_{L} \big(   E\cdot(1+a) v\sqrt{M}+ S     \big)(\tau), \phi      \big \rangle_{x,v}{\rm d}\tau \nonumber\\
=&\, \int_0^t \big \langle  \big(   E\cdot(1+a) v\sqrt{M}+ S     \big)(\tau) , e^{(t-\tau)\mathcal{B}^*}_{L} \phi      \big \rangle_{x,v}{\rm d}\tau \nonumber\\
\lesssim&\, \int_0^t \big( \big|\big(   E\cdot(1+a) v\sqrt{M}+ S     \big) (\tau)\big|_{2}  ,\big|e_{L}^{(t-\tau)\mathcal{B}^*}\phi\big|_{2}  \big)  {\rm d}\tau \nonumber\\
\lesssim&\, \sup_{0\leq t\leq T_1}\big\|\big(   E\cdot(1+a) v\sqrt{M}+ S     \big)(t)\big\|_{L_v^2(\dot B_{2,\infty}^{-\frac{3}{2}})} \int_0^{\infty}\| e_{L}^{\tau \mathcal{B}^*}\phi\|_{L_v^2(\dot B_{2,1}^{\frac{3}{2}})} {\rm d}\tau \nonumber\\
\lesssim &\, \sup_{0\leq t\leq T_1}\big\|\big(   E\cdot(1+a) v\sqrt{M}+ S     \big)(t)\big\|_{L_v^2(\dot B_{2,\infty}^{-\frac{3}{2}})}  \|\phi\|_{L_v^2(\dot B_{2,1}^{-\frac{1}{2}})},
\end{align}
where we have employed the low-frequency estimate from \eqref{B.15} that
\begin{align*}
  \int_0^{\infty}\| e_{L}^{\tau \mathcal{B}^*}\phi\|_{L_v^2(\dot B_{2,1}^{\frac{3}{2}})} {\rm d}\tau\lesssim    \|\phi\|_{L_v^2(\dot B_{2,1}^{-\frac{1}{2}})}. 
\end{align*}
By applying the duality argument and  the inequalities  \eqref{G3.5} and \eqref{newinterpolation}, we find that
\begin{align}\label{G3.6}
I_1\lesssim&\, \sup_{0\leq t\leq T_1}\big\|\big(   E\cdot(1+a) v\sqrt{M}+ S     \big)(t)\big\|_{L_v^2(\dot B_{2,\infty}^{-\frac{3}{2}})}  \nonumber\\
\lesssim&\, \sup_{0\leq t\leq T_1} \|S(t)\|_{L_v^2(\dot B^{-\frac{3}{2}}_{2,\infty})}+ \sup_{0\leq t\leq T_1} \|E(t)\|_{ \dot B^{-\frac{3}{2}}_{2,\infty} }+\sup_{0\leq t\leq T_1} \|E(t)\|_{ \dot B^{-\frac{3}{2}}_{2,\infty} } \|a\|_{L_v^2(\dot B^{\frac{3}{2}}_{2,1})}\nonumber\\
\lesssim&\, \sup_{0\leq t\leq T_1} \|S(t)\|_{L_v^2(\dot B^{-\frac{3}{2}}_{2,\infty})}+ \Big(1+\|f\|_{L_v^2(\dot B^{\frac{1}{2}}_{2,\infty}\cap \dot H^{N})}\Big)\sup_{0\leq t\leq T_1} \|E(t)\|_{ \dot B^{-\frac{3}{2}}_{2,\infty} } \nonumber\\
\lesssim&\,  \sup_{0\leq t\leq T_1} \|S(t)\|_{L_v^2(\dot B^{-\frac{3}{2}}_{2,\infty})}+ \sup_{0\leq t\leq T_1} \|E(t)\|_{ \dot B^{-\frac{3}{2}}_{2,\infty} }.
\end{align}

To deal with the term $I_2$, we shall use Lemma \ref{LB.2} to overcome the difficulties arising from the weighted $v$ and the derivative $\nabla_{v}$. 
Define
\begin{align*}
\mathfrak{F}(t,x,v)
:=
\int_0^t \dot S_{j_0}e^{(t-\tau)\mathcal B}_{L}
\Big(\nabla_v\cdot G-\frac12 v\cdot G\Big)(\tau){\rm d}\tau.
\end{align*}
It follows from Lemma \ref{LB.2} that for every dyadic block $\dot\Delta_j$ with $j\le j_0$, the following low-frequency square-type estimate holds:
\begin{equation}\label{DLNGG3.7}
\|\dot\Delta_j \mathfrak{F}(t)\|_{L^2_{x,v}}^2
\lesssim
\int_0^t e^{-c2^{2j}(t-\tau)}
\|\dot\Delta_j G(\tau)\|_{L^2_{x,v}}^2{\rm d}\tau,
\end{equation}
for $j\leq j_0$.
Then, for any $\phi\in L_v^2(\dot B_{2,1}^{-\frac{1}{2}})$, we claim that
\begin{equation}\label{DLNGG3.8}
\big|
\left\langle \mathfrak{F}(t),\phi\right\rangle_{x,v}
\big|
\lesssim
\sup_{0\le \tau\le T_1}
\|G(\tau)\|_{L_v^2(\dot B_{2,\infty}^{-\frac12})}
\,
\|\phi\|_{L_v^2(\dot B_{2,1}^{-\frac12})}.
\end{equation}
 
In fact, by the definition of the homogeneous Besov norm, we have
\begin{equation*} 
\|\dot\Delta_j G(\tau)\|_{L^2_{x,v}}
\leq \sup_{0\leq t \leq T_1}\|G(t)\|_{L_v^2(\dot B_{2,\infty}^{-\frac12})}2^{\frac{j}{2}},
\end{equation*}
which together with \eqref{DLNGG3.7} gives
\begin{align*}
\|\dot\Delta_j \mathfrak{F}(t)\|_{L^2_{x,v}}^2
\lesssim&\,
\sup_{0\leq t \leq T_1}\|G(t)\|_{L_v^2(\dot B_{2,\infty}^{-\frac12})}^2 2^j
\int_0^t e^{-c2^{2j}(t-\tau)}{\rm d}\tau \nonumber\\
\lesssim&\, \sup_{0\leq t \leq T_1}\|G(t)\|_{L_v^2(\dot B_{2,\infty}^{-\frac12})}^2 2^{-j},
\end{align*}
for all $j\leq j_0$.
Thus, we have
\begin{equation}\label{DLNGG3.9}
\|\dot\Delta_j \mathfrak{F}(t)\|_{L^2_{x,v}}
\lesssim
\sup_{0\leq t \leq T_1}\|G(t)\|_{L_v^2(\dot B_{2,\infty}^{-\frac12})}2^{-\frac{j}{2}},
\end{equation}
for all $j\leq j_0$.
Multiplying \eqref{DLNGG3.9} by $2^{\frac{j}{2}
}$ and taking the supremum over $j\leq j_0$, we obtain
\begin{equation}\label{DLNGG3.10}
\|\mathfrak{F}_{L}(t)\|_{L_v^2(\dot B_{2,\infty}^{\frac{1}{2}})}
=
\sup_{j\leq j_0} 2^{\frac{j}{2}}\|\dot\Delta_j \mathfrak{F}(t)\|_{L^2_{x,v}}
\lesssim \sup_{0\leq t \leq T_1}\|G(t)\|_{L_v^2(\dot B_{2,\infty}^{-\frac12})} .
\end{equation}

We now shift our focus to the duality argument. Owing to the almost orthogonality of the homogeneous dyadic decomposition, it holds that
\begin{align*}
\langle \mathfrak{F}(t),\phi\rangle_{x,v}=\sum_{j\in\mathbb Z}\langle \dot\Delta_j \mathfrak{F}(t),\widetilde{\dot\Delta}_j\phi\rangle_{x,v},    
\end{align*}
with
\begin{align*}
 \widetilde{\dot\Delta}_j:=\dot\Delta_{j - 1}+\dot\Delta_j+\dot\Delta_{j + 1}.   
\end{align*}
Since $\mathfrak{F}(t)$ is localized in the low-frequency region by $\dot S_{j_0}$, only the indices $j\leq j_0 + N$ contribute, where $N$ is a universal integer that depends only on the dyadic partition. Therefore, we have
\begin{align*}
\left|
\langle \mathfrak{F}_{L}(t),\phi\rangle_{x,v}
\right|
\le
\sum_{j\le j_0+N}
\|\dot\Delta_j \mathfrak{F}(t)\|_{L^2_{x,v}}
\,
\|\widetilde{\dot\Delta}_j\phi\|_{L^2_{x,v}}.    
\end{align*}
Insert the factors $2^{\frac{j}{2}}$ and $2^{-\frac{j}{2}}$ to obtain
\begin{align*}
\left|
\langle \mathfrak{F}_{L}(t),\phi\rangle_{x,v}
\right|
\leq&\,
\sum_{j\le j_0+N}
\Big(2^{\frac{j}{2}}\|\dot\Delta_j \mathfrak{F}(t)\|_{L^2_{x,v}}\Big)
\Big(2^{-\frac{j}{2}}\|\widetilde{\dot\Delta}_j\phi\|_{L^2_{x,v}}\Big)\nonumber\\
\leq&\, \Big(\sup_{j\le j_0+N}2^{\frac{j}{2}}\|\dot\Delta_j \mathfrak{F}(t)\|_{L^2_{x,v}}\Big)
\sum_{j\le j_0+N}2^{-\frac{j}{2}}\|\widetilde{\dot\Delta}_j\phi\|_{L^2_{x,v}}.
\end{align*}
Substituting \eqref{DLNGG3.10} into the above inequality, one gets
\begin{align*}
\left|
\langle \mathfrak{F}_{L}(t),\phi\rangle_{x,v}
\right|
\lesssim&\,
\sup_{0\leq t \leq T_1}\|G(t)\|_{L_v^2(\dot B_{2,\infty}^{-\frac12})}
\sum_{j\le j_0+N}2^{-\frac{j}{2}}\|\widetilde{\dot\Delta}_j\phi\|_{L^2_{x,v}} \nonumber\\
\lesssim&\, \sup_{0\leq t \leq T_1}\|G(t)\|_{L_v^2(\dot B_{2,\infty}^{-\frac12})}
\sum_{j\le j_0+N}2^{-\frac{j}{2}}\big(\| {\dot\Delta}_{j-1}\phi\|_{L^2_{x,v}} +\| {\dot\Delta}_{j}\phi\|_{L^2_{x,v}}+\| {\dot\Delta}_{j+1}\phi\|_{L^2_{x,v}}\big) \nonumber\\
\lesssim&\, \sup_{0\leq t \leq T_1}\|G(t)\|_{L_v^2(\dot B_{2,\infty}^{-\frac12})}\sum_{j\in\mathbb Z}2^{-\frac{j}{2}}\|\dot\Delta_j\phi\|_{L^2_{x,v}}\nonumber\\
\lesssim&\, \sup_{0\leq t \leq T_1}\|G(t)\|_{L_v^2(\dot B_{2,\infty}^{-\frac12})}\|\phi\|_{L_v^2(\dot B_{2,1}^{-\frac12})}.
\end{align*}
Consequently, the inequality \eqref{DLNGG3.8} holds.

By combining \eqref{DLNGG3.8} with the duality argument, we have
\begin{align}\label{G3.8}
 I_2\lesssim     \sup_{0\leq t\leq T_1}\|G(t)\|_{L_v^2(\dot B_{2,\infty}^{-\frac{1}{2}})} 
 \lesssim&\,  \sup_{0\leq t\leq T_1}\|E(t)f(t)\|_{L_v^2(\dot B_{2,\infty}^{-\frac{1}{2}})}\nonumber\\
 \lesssim &\,  \sup_{0\leq t\leq T_1}\|E(t) \|_{ \dot B_{2,\infty}^{-\frac{1}{2}}} \|f(t)\|_{L_v^2(\dot B^{\frac{3}{2}}_{2,1})}\nonumber\\
 \lesssim&\, \delta  \sup_{0\leq t\leq T_1}\|E(t) \|_{ \dot B_{2,\infty}^{-\frac{1}{2}}}.
\end{align}
Substituting the estimates \eqref{G3.6} and \eqref{G3.8} into \eqref{G3.4}, 
we   obtain the desired estimate \eqref{G3.2}.    
\end{proof}

\subsection{Estimates of the solution at high frequencies}
Since the estimate of $f_{L}(t,x,v)$ has been established, we now focus on the high frequency estimate of $f(t,x,v)$.

\begin{lem}\label{L3.2}
For the classical solution $f$ to the problem   \eqref{G1.4}, there exists a positive constant  $\kappa_1 >0$ such that  
\begin{align}\label{G3.9}
&\frac{{\rm d}}{{\rm d}t}\|  f(t)\|_{L_v^2(\dot H^1\cap \dot H^N)}^2+\ka_1 \Big(\|b\|_{\dot H^1\cap \dot H^N}^2+\sum_{1\leq k\leq N}\|\nabla^k\{\mathbf{I}-\mathbf{P}\}f(t)\|_{\nu}^2\Big)     \nonumber\\
&\quad\leq \eta  \|  a(t)\|_{\dot H^1\cap \dot H^N}^2+C\sup_{{0\leq t\leq T_1}}\|E(t)\|_{ H^N}^2  +C\sup_{{0\leq t\leq T_1}}\|S(t)\|_{L_v^2(H^N)}^2 , 
\end{align}
for any $0\leq t \leq T_1$ and any $T_1>0$. Here and below, $0 < \eta < 1$ is a sufficiently small constant.
\end{lem}

\begin{proof}
Taking the differentiation $\partial^\alpha$ with $1\leq |\alpha|\leq N$ of \eqref{G1.4}$_1$, then multiplying the resulting identity by $\partial^\alpha f$ and integrating over $\mathbb R_x^3\times\mathbb R^3_v$ yields
\begin{align} \label{G3.10}
& \frac{{\rm d}}{{\rm d}t}\|\partial^\alpha f\|_{L_{x,v}^2}^2+2\kappa_2 \big(\|\partial^\alpha \{\mathbf{I}-\mathbf{P}\}f\|_{\nu}^2+2\|\partial^\alpha b\|_{L^2}^2\big)\nonumber\\
\lesssim&\,  |\langle v\sqrt{M}\cdot \partial^\alpha E,\partial^\alpha f\rangle_{x,v} |+|\langle\partial^\alpha S,\partial^\alpha f\rangle_{x,v} |+  |\langle \partial^\alpha[E\cdot vf],\partial^\alpha f\rangle _{x,v}|+|\langle \partial^\alpha [E\cdot\nabla_v f],\partial^\alpha f  \rangle_{x,v}  |\nonumber\\
\equiv:&\, \sum_{i=1}^4 J_{i},
\end{align}
for some constant $\kappa_{2}>0$.

Applying H\"{o}lder’s and Young’s inequalities along with the macro-micro decomposition \eqref{G2.5} yields
\begin{align}\label{G3.11}
J_1+J_2\leq&\,  C\|f\|_{L_v^2(\dot H^1\cap \dot H^N)} \Big(  \sup_{0\leq t\leq T_1}\|E(t)\|_{\dot H^1\cap \dot H^N}+   \sup_{0\leq t\leq T_1}\|S(t)\|_{L_v^2(\dot H^1\cap \dot H^N)}    \Big)   \nonumber\\
\leq&\, \eta\Big(\|(a,b)(t)\|_{  {\dot H^1\cap \dot H^N }} ^2 +\sum_{1\leq k\leq N}\|\nabla^k\{\mathbf{I}-\mathbf{P}\}f(t)\|_{\nu}^2\Big)+C\sup_{ 0\leq t\leq T_1}\|E(t)\|_{ H^N}^2 \nonumber\\
&+C\sup_{0\leq t\leq T_1}\|S(t)\|_{L_v^2(H^N)}^2.
\end{align}
On the other hand, by using the macro-micro decomposition \eqref{G2.5} and Lemma \ref{LA.1}, we get
\begin{align}\label{G3.12}
J_3+J_4\leq&\, C    \sup_{0\leq t\leq T_1}\|\nabla (E  f)(t)\|_{L_v^2(H^{N-1})} \|\nabla (v\{\mathbf{I}-\mathbf{P}\}f+v\mathbf{P}f)(t)\|_{L_v^2(H^{N-1})}\nonumber\\
&+ C    \sup_{0\leq t\leq T_1}\|\nabla (E  f)(t)\|_{L_v^2(H^{N-1})} \|\nabla (\nabla_{v}\{\mathbf{I}-\mathbf{P}\}f+\nabla_{v}\mathbf{P}f)(t)\|_{L_v^2(H^{N-1})}\nonumber\\
\leq&\,\eta\Big(\|(a,b)(t)\|_{ { \dot H^1\cap \dot H^N } } ^2 +\sum_{1\leq k\leq N}\|\nabla^k\{\mathbf{I}-\mathbf{P}\}f(t)\|_{\nu}^2\Big) +C\sup_{0\leq t\leq T_1} \|E(t)\|_{\dot H^1\cap H^N}^2 \|f(t)\|_{L_v^2(\dot H^1\cap\dot H^N)}^2 \nonumber\\
\leq&\, \eta\Big(\|(a,b)(t)\|_{ { \dot H^1\cap \dot H^N }} ^2 +\sum_{1\leq k\leq N}\|\nabla^k\{\mathbf{I}-\mathbf{P}\}f(t)\|_{\nu}^2\Big) +C\sup_{0\leq t\leq T_1} \|E(t)\|_{ H^N}^2 .
\end{align}
Inserting \eqref{G3.11} and \eqref{G3.12} into \eqref{G3.10} and then taking the sum over $1\leq|\alpha|\leq N$, we obtain \eqref{G3.9} and hence complete the proof of Lemma \ref{L3.2}.
\end{proof}

Next, we provide the estimate of $\|  a(t)\|_{\dot H^2\cap \dot H^N}$. For this purpose, motivated by \cite{Ks-1983}, we   first  study the following coupled hyperbolic-parabolic system satisfied by $a$ and $b$, which are the coefficient functions of macroscopic part $\mathbf{P}f$:
\begin{equation}\label{G3.13}
\left\{\begin{aligned}
&\partial_{t}a+\nabla \cdot {b} =\int_{\mathbb R^3_{v}}\sqrt{M}S\,{\rm d}v,\\
&\partial_{t} b_i+b_i+\partial_{i} a+\sum_{j=1}^3\partial_j\Gamma_{ij}(\{\mathbf{I}-\mathbf{P}\}f) =E_i(1+a)+\int_{\mathbb R^3_{v}}v\sqrt{M}S\,{\rm d}v,    
 \end{aligned}
 \right.
\end{equation}
for $1\leq i,j\leq 3$, where $\Gamma_{ij}(\cdot)$ is defined by \eqref{gammaa}.

\begin{lem}\label{L3.3}
For the classical solution $f$ to the problem   \eqref{G1.4}, there exists a positive constant  $\kappa_3 >0$ such that  
\begin{align}\label{G3.14}
&-\frac{{\rm d}}{{\rm d}t}\sum_{1\leq|\alpha|\leq N-1}\int_{\mathbb{R}^3}\partial^\alpha a\partial^\alpha \nabla\cdot b\,{\rm d}x+\kappa_3\| a\|_{\dot H^2\cap \dot H^N}^2\nonumber\\
&\quad\lesssim
\|b\|_{\dot H^1\cap\dot H^N}^2+\sum_{1\leq k\leq N}\|\nabla^k\{\mathbf{I}-\mathbf{P}\}f\|_{L_{x,v}^2}^2+\sup_{{0\leq t\leq T_1}} \|E(t)\|_{ H^N }^2+\sup_{{0\leq t\leq T_1}} \|S(t)\|_{L_v^2(H^N)}^2, 
\end{align}
for any $0\leq t \leq T_1$. 
\end{lem}

\begin{proof}
Applying $\partial^\alpha$ with $1\leq|\alpha|\leq N - 1$ to the second equation in \eqref{G3.13}, and then multiplying the resulting identity by $\partial_i a$, one gets
\begin{align}\label{G3.15}
   \|\partial^\alpha \nabla  a \|_{L^2}^2 
= \,&\sum_{i=1}^3 \int_{\mathbb{R}^3} \partial^\alpha \partial_i a \partial^\alpha\Big[-\partial_t b_i-b_i-\sum_{j=1}^3 \partial_j \Gamma_{ij}(\{\mathbf{I}-\mathbf{P}\}f)+E_i(1+a) \Big] \mathrm{d} x \nonumber\\
&+\sum_{i=1}^3 \int_{\mathbb{R}^3_{x}}\!\int_{\mathbb R^3_{v}}v \partial^\alpha \partial_i a \sqrt{M}\partial^\alpha S\,  {\rm d}v\mathrm{d} x \nonumber\\
=&-\frac{\mathrm{d}}{\mathrm{d} t} \sum_{i=1}^3 \int_{\mathbb{R}^3} \partial^\alpha \partial_i a \partial^\alpha b_i \mathrm{d} x+\sum_{i=1}^3 \int_{\mathbb{R}^3} \partial^\alpha \partial_i \partial_t a \partial^\alpha b_i\,\mathrm{d} x \nonumber\\
& -\sum_{i=1}^3 \int_{\mathbb{R}^3} \partial^\alpha \partial_i a \partial^\alpha\Big[b_i+\sum_{j=1}^3 \partial_j\Gamma_{ij}(\{\mathbf{I}-\mathbf{P}\} f)+E_i(1+a)\Big] \mathrm{d} x \nonumber\\
&+\sum_{i=1}^3 \int_{\mathbb{R}^3_{x}}\!\int_{\mathbb R^3_{v}}v \partial^\alpha \partial_i a \sqrt{M}\partial^\alpha S\,  {\rm d}v\mathrm{d} x \nonumber\\
\equiv:&\, -\frac{\mathrm{d}}{\mathrm{d} t} \sum_{i=1}^3 \int_{\mathbb{R}^3} \partial^\alpha \partial_i a \partial^\alpha b_i\, \mathrm{d} x+K_1+K_2+K_3 . 
\end{align}
It follows from the first part of \eqref{G3.13} and integration by parts that
\begin{align}\label{G3.16}
  -\frac{\mathrm{d}}{\mathrm{d} t} \sum_{i=1}^3 \int_{\mathbb{R}^3} \partial^\alpha \partial_i a \partial^\alpha b_i \mathrm{d} x = {\frac{{\rm d}}{{\rm d}t}}\int_{\mathbb{R}^3}  \partial^\alpha a\partial^\alpha \nabla\cdot b\,  {\rm d}x.   
\end{align}
Similarly, we have
\begin{align}\label{G3.17}
K_1=&\,     - \int_{\mathbb{R}^3}\partial^\alpha\partial_t a\partial^\alpha \nabla\cdot b\,{\rm d}x\nonumber\\
=&\,\|\partial^\alpha\nabla\cdot b\|_{L^2}^2-\int_{\mathbb R^3_{x}}\!\int_{\mathbb R^3_{v}} \partial^\alpha\nabla\cdot b\sqrt{M}\partial^\alpha S\, {\rm d}v{\rm d}x \nonumber\\
\lesssim&\, \|b\|_{\dot H^1\cap\dot H^N}^2+ \sup_{0\leq t\leq T_1}\|S(t)\|_{L_v^2(\dot H^1\cap \dot H^{N-1})}^2,
\end{align}
and
\begin{align}\label{G3.18}
K_3\leq&\, C\|a\|_{\dot H^2\cap \dot H^N} \sup_{0\leq t\leq T_1}\|S(t)\|  _{L_v^2(\dot H^1\cap \dot H^{N-1})}  \nonumber\\
\leq&\, \eta \|a\|_{\dot H^2\cap \dot H^N}^2+C \sup_{0\leq t\leq T_1}\|S(t)\|  _{L_v^2(\dot H^1\cap \dot H^{N-1})} ^2.
\end{align}

Finally, for the remaining term $K_2$, since the moment functional appearing in $\Gamma_{ij}$ can absorb any velocity derivative and any velocity weight, by applying Lemma \ref{LA.1}, we derive
\begin{align}\label{G3.19}
 K_2\leq&\, \eta \|a\|_{\dot H^2\cap \dot H^N}^2+C\|b\|_{\dot H^1\cap \dot H^N}^2+C   \sum_{1\leq k\leq N}\|\nabla^k\{\mathbf{I}-\mathbf{P}\}f\|_{L_{x,v}^2}^2+C \sup_{0\leq t\leq T_1 } \|E(t)\|_{\dot H^1\cap \dot H^N}^2 \nonumber\\
 &+ C \sup_{0\leq t\leq T_1 } \|E(t)\|_{\dot H^1\cap \dot H^N}  \|a\|_{\dot H^1\cap \dot H^N} \nonumber\\
 \leq &\,  \eta \|a\|_{\dot H^2\cap \dot H^N}^2+C\|b\|_{\dot H^1\cap \dot H^N}^2+C   \sum_{1\leq k\leq N}\|\nabla^k\{\mathbf{I}-\mathbf{P}\}f\|_{L_{x,v}^2}^2+C \sup_{0\leq t\leq T_1 } \|E(t)\|_{  H^N}^2  .
\end{align}
Plugging the estimates \eqref{G3.16}, \eqref{G3.17}, \eqref{G3.18} and \eqref{G3.19} into \eqref{G3.15}, we end up with \eqref{G3.14} and then complete the proof of Lemma \ref{L3.3}.
\end{proof}

\subsection{Proof of Theorem \ref{Th1}}
With Lemmas \ref{L3.1}–\ref{L3.3} established, we are now in a postion to prove the global existence of classical solutions to the Cauchy problem \eqref{G1.4}.
\begin{proof}[Proof of Theorem \ref{Th1}]
Define a temporal energy functional $\mathcal{E}^H(f(t))$ in the high-frequency regime and the corresponding dissipation rate $\mathcal{D}^H(f(t))$ as
follows:
\begin{align} \label{meab}
\mathcal{E}^H(f(t)):=&\,\|f(t)\|_{L_v^2(\dot H^1\cap \dot H^N)}^2 -\tau_2 \sum_{1\leq|\alpha|\leq N-1} \int_{\mathbb{R}^3}\partial^\alpha a\partial^\alpha \nabla\cdot b\,{\rm d}x, \\\label{meac}
\mathcal{D}^H(f(t)):=&\, \|  a(t)\|_{\dot H^2\cap \dot H^N}^2+\|b(t)\|_{\dot H^1\cap \dot H^N}^2+\sum_{1\leq k \leq N}\|\nabla^k\{\mathbf{I}-\mathbf{P}\}f(t)\|_{\nu}^2 ,
\end{align}
where $0 < \tau_2 \ll 1$ is a sufficiently small constant.
Then, for \eqref{meab}, by Young's inequality, we also have
\begin{align}\label{add.efequiv}
\mathcal{E}^H(f(t))\backsim  \|f(t)\|_{L_v^2(\dot H^1\cap \dot H^N)}^2,
\end{align}
uniformly for all $0\leq t\leq T_1$. Thus, by combining \eqref{G3.9} and \eqref{G3.14}, we  infer that
\begin{align}\label{G3.22}
\frac{{\rm d}}{{\rm d}t} \mathcal{E}^H(f(t))+ \kappa_4  \mathcal{D}^H(f(t))\leq \eta \|a(t)\|_{\dot H^1}^2+C\sup_{{0\leq t\leq T_1}}\|E(t)\|_{ H^N}^2  +C\sup_{{0\leq t\leq T_1}}\|S(t)\|_{L_v^2(H^N)}^2 , 
\end{align}
for some constant $\kappa_4 > 0$ and any $0\leq t\leq T_1$ with $T_1>0$. Adding the term $\kappa_{4}\|f_{L}(t)\|_{L_v^2(\dot B_{2,\infty}^{\frac{1}{2}})}^2$ to both sides of  \eqref{G3.22}, using   \eqref{G3.2} and recalling \eqref{meac}, we arrive at
\begin{align}\label{G3.23}
\frac{{\rm d}}{{\rm d}t} \mathcal{E}^H(f(t))+ \kappa_5 \|f(t)\|_{L_v^2(\dot H^1\cap \dot H^N)}^2 \lesssim \sup_{{0\leq t\leq T_1}}\|E(t)\|_{ H^N}^2  +\sup_{{0\leq t\leq T_1}}\|S(t)\|_{L_v^2(H^N)}^2+\|f_0\|_{L_v^2(\dot B^\frac{1}{2}_{2,\infty})}^2 , 
\end{align}
for some constant $\kappa_5>0$. Here, we have utilized the following inequality: 
\begin{align*}
 \|f(t)\|_{L_v^2(\dot H^1\cap \dot H^N)}^2\leq&\, \|(a,b)(t)\|_{ \dot H^1\cap \dot H^N }^2+\sum_{1\leq k\leq N}\|\nabla^k\{\mathbf{I}-\mathbf{P}\}f(t)\|_{L_{x,v}^2}^2  \\
 \leq&\, \|f_{L}(t)\|_{L_v^2(\dot B_{2,\infty}^{\frac{1}{2}})}^2+ \mathcal{D}^H(f(t)). 
\end{align*}
Recall \eqref{add.efequiv}. By applying Gr\"{o}nwall’s inequality to \eqref{G3.23}, we  further get
\begin{align*}
 \|f(t)\|_{L_v^2(\dot H^1\cap \dot H^N)}^2\lesssim &\, \sup_{0\leq t\leq T_1}\Big\{\|E(t)\|_{ H^N}^2+\|S(t)\|_{L_v^2(H^N)}^2 \Big\}\int_0^t e^{-\kappa_{5}\tau}{\rm d}\tau\nonumber\\
& +\|f_0\|_{L_v^2(\dot B^\frac{1}{2}_{2,\infty})}^2\int_0^t e^{-\kappa_{5}\tau}{\rm d}\tau +e^{-\ka_{5} t} \|f_0\|_{L_v^2(\dot H^1\cap \dot H^N)}^2\nonumber\\
 \lesssim&\, \sup_{0\leq t\leq T_1} \Big\{\|E(t)\|_{ H^N}^2+\|S(t)\|_{L_v^2(H^N)}^2  \Big\}+\|f_0\|_{L_v^2(\dot B^{\frac{1}{2}}_{2,\infty}\cap \dot H^N)}^2,
\end{align*}
which implies that
\begin{align}\label{G3.24}
 \|f(t)\|_{L_v^2(\dot H^1\cap \dot H^N)} 
 \lesssim&\, \sup_{0\leq t\leq T_1} \Big\{\|E(t)\|_{ H^N}+\|S(t)\|_{L_v^2(H^N)}  \Big\}+\|f_0\|_{L_v^2(\dot B^{\frac{1}{2}}_{2,\infty}\cap \dot H^N)}.
\end{align}
Adding \eqref{G3.2} and \eqref{G3.24}, we arrive at
\begin{align}\label{G3.25}
 \|f(t)\|_{L_v^2(\dot B^{\frac{1}{2}}_{2,\infty}\cap \dot H^N)} 
 \lesssim&\, \sup_{0\leq t\leq T_1} \Big\{\|E(t)\|_{  \dot B^{-\frac{3}{2}}_{2,\infty}\cap\dot H^N}+\|S(t)\|_{L_v^2( \dot B^{-\frac{3}{2}}_{2,\infty}\cap\dot H^N)}  \Big\}+\|f_0\|_{L_v^2(\dot B^{\frac{1}{2}}_{2,\infty}\cap \dot H^N)},
\end{align}
for any $0\leq t \leq T_1$.
In addition, a priori estimates \eqref{G3.1} can be justified by choosing
\begin{align*}
\|E(t) \|_{  \dot B_{2,\infty}^{-\frac{3}{2}}\cap \dot H^N   }+\|S(t) \|_{ L_v^2(\dot B_{2,\infty}^{-\frac{3}{2}}\cap \dot H^N)  } + \|f_0\|_{L_v^2(\dot B_{2,\infty}^{\frac{1}{2}}\cap \dot H^N)}
\end{align*}
to be sufficiently small.

The proof for the local existence of classical solutions is standard. For the sake of brevity, we omit it here. Interested readers can refer to \cite{DFT-2010-CMP} for a more detailed proof on a simpler KFP equation
than \eqref{G1.4}$_1$. Subsequently, the global existence of solutions to \eqref{G1.4} is derived from the obtained {\it global a priori estimates} and the continuity argument. Moreover, \eqref{G1.6} holds for all $t\geq 0$ in terms of  \eqref{G3.25}. Finally, by applying the maximum principle (cf. \cite{Gy-IUMJ-2004}), we obtain
\begin{align*}  F(t,x,v) = M+\sqrt{M}f (t,x,v)\geq 0. \end{align*} 
Thus, we complete the proof of Theorem \ref{Th1}.
\end{proof}

\section{Uniform stability of the KFP equation}\label{sec4}

In this section, we move forward to studying the uniform stability of global solutions to the  Cauchy problem \eqref{G1.4} in order to prove Theorem \ref{Th2}.  For the purpose, let $f$ and $g$ be a pair of solutions to the Cauchy problem \eqref{G1.4} with initial data $f_0(x,v)$ and $g_0(x,v)$, respectively. Suppose that the initial data $f_0$ and $g_0$, the external force $E(t,x)$, and the inhomogeneous term $S(t,x,v) $ satisfy \eqref{G1.5} in Theorem \ref{Th1}. 
We first present the difference equation between two solutions $f$ and $g$ by denoting $\widetilde f$ as  
\begin{align*} 
\widetilde f  =  f - g. 
\end{align*}
From \eqref{G1.4}, it can be inferred that
 \begin{equation}\label{G4.1}
\left\{\begin{aligned}
&\partial_t \widetilde f+v\cdot\nabla \widetilde f-\mathcal{L} \widetilde f=\, E\cdot \nabla_{v} \widetilde f+\frac{1}{2}E\cdot v\widetilde f 
=:  \nabla_v \cdot \widetilde G -\frac{1}{2}v\cdot \widetilde G+\widetilde aE\cdot v \sqrt{M},\\
&\widetilde f(t,x,v)|_{t=0}=\widetilde f_0(x,v)= f_0 (x,v)-g_0 (x,v),
 \end{aligned}
 \right.
\end{equation} 
where $\widetilde  G$ is given by 
\begin{align*}
 \widetilde  G:=-E\{\mathbf{I}-\mathbf{P}_0\}\widetilde f.   
\end{align*}
Note that $\{\mathbf{I}-\mathbf{P}\} (\nabla_v \cdot \widetilde G -\frac{1}{2}v\cdot \widetilde G)=0$.

Below, we shall establish the time-weighted estimate in hybrid homogeneous Besov spaces $L_v^2(\dot B^{\frac{1}{2}}_{2,\infty}\cap\dot H^{N - 1})$ for $\widetilde f(t,x,v)$ with $N\geq 3$ at low frequencies and high frequencies respectively.

\subsection{Time-weighted estimates at low and high frequencies }
\begin{lem} \label{L4.1}
Under the assumptions of Theorem \ref{Th2}, for any $s\in\big[-\frac{1}{2}+\varepsilon,\frac{3}{2}-\varepsilon\big]$ with $s\geq s_0$ and $s_0\in\big(-\frac{1}{2},\frac{1}{2}\big]$,  we have
\begin{align}\label{G4.2}
\sup_{t\geq 0} (1+t)^{\frac{s-s_0}{2}}\|\widetilde f\|_{L_v^2(\dot B_{2,\infty}^s)} \lesssim \delta_0 {\widetilde{\mathcal{D}}_{\varepsilon}}+\|\widetilde f_0\|_{L_v^2(\dot B_{2,\infty}^{s_0}\cap \dot B_{2,\infty}^{\frac{3}{2}-\varepsilon})},   
\end{align}
where  $\widetilde {\mathcal{D}}_{\varepsilon}$ is given by
\begin{align*}
{\widetilde {\mathcal{D}}_{\varepsilon}}:=\sup_{s_1\leq \bar s\leq \frac{3}{2}-\varepsilon} \sup_{t\geq 0}\,(1+t)^{\frac{\bar s-s_0}{2}} \|\widetilde f (t)\|_{L_v^2(\dot B_{2,\infty}^{\bar s}\cap \dot H^{N-1})},\quad s_1=\max\{0,{s_0}\},   
\end{align*}
and $\delta_0$ is defined by \eqref{G1.5}.
\end{lem}

\begin{proof}
Through the Duhamel principle, we can express the solution of \eqref{G4.1} as
\begin{align}\label{G4.3}
\widetilde f(t)=e^{t\mathcal B} \widetilde f_0+\int_0^t e^{\tau\mathcal{B}}\widetilde S_{f}(t-\tau){\rm d}\tau,  
\end{align}
where
\begin{align*}
 \widetilde S_{f}:=   \nabla_v \cdot \widetilde G -\frac{1}{2}v\cdot \widetilde G+E\cdot \widetilde av \sqrt{M}.
\end{align*}
We first consider the case in which $s\in \big[-\frac{1}{2}+\varepsilon,\frac{1}{2}\big)$. By taking the $\dot B^{s}_{2,\infty}$-norm of \eqref{G4.3} and leveraging \eqref{B.23}--\eqref{B.24} and the triangle inequality, one obtains
\begin{align}\label{G4.4}
\|\widetilde f(t)\|_{L_v^2(\dot B_{2,\infty}^s )}^2\lesssim&\, \|e^{t \mathcal{B}} \widetilde f_0\|_{L_v^2(\dot B_{2,\infty}^s)}^2+\bigg \|\int_0^t e_{L}^{\tau\mathcal{B}}\widetilde S_{f}(t-\tau){\rm d}\tau \bigg\|_{L_v^2(\dot B_{2,\infty}^s)}  ^2 +\bigg \|\int_0^t e_{H}^{\tau\mathcal{B}}\widetilde S_{f}(t-\tau){\rm d}\tau \bigg\|_{L_v^2(\dot B_{2,\infty}^s)}^2   \nonumber\\
\lesssim&\, (1+t)^{-{s-s_0}}\|\widetilde f_0\|_{L_v^2(\dot B_{2,\infty}^{s_0}\cap \dot B_{2,\infty}^s)}^2+ \int_0^t (1+\tau)^{-\frac{1}{2}-s}\|  \widetilde G(t-\tau)   \|_{L_v^2(\dot B_{2,\infty}^{-\frac{1}{2}})}^2{\rm d}\tau \nonumber\\
&+ \int_0^t (1+\tau)^{-\frac{1}{2}-s}\|E\widetilde a(t-\tau)\|_{ \dot B_{2,\infty}^{-\frac{1}{2}} }^2{\rm d}\tau + \int_0^t  e^{-\lambda_{2}\tau} \|  \widetilde G(t-\tau)   \|_{L_v^2(\dot B_{2,\infty}^{s})}^2 {\rm d}\tau\nonumber\\
&+ \int_0^t  e^{-\lambda_{2}\tau} \|E\widetilde a(t-\tau)\|_{ \dot B_{2,\infty}^{s} }^2 {\rm d}\tau \nonumber\\
\lesssim&\,(1+t)^{- {s-s_0} }\|\widetilde f_0\|_{L_v^2(\dot B_{2,\infty}^{s_0}\cap \dot B_{2,\infty}^s)}^2+ \int_0^t (1+\tau)^{-\frac{1}{2}-s}\|  \widetilde G(t-\tau)   \|_{L_v^2(\dot B_{2,\infty}^{-\frac{1}{2}}\cap\dot B^{\frac{1}{2}}_{2,\infty})}^2{\rm d}\tau  \nonumber\\
&+ \int_0^t (1+\tau)^{-\frac{1}{2}-s}\|  \widetilde E\widetilde a(t-\tau)   \|_{{\dot B_{2,\infty}^{-\frac{1}{2}}\cap\dot B^{\frac{1}{2}}_{2,\infty}}}^2{\rm d}\tau\nonumber\\
\equiv:&\, (1+t)^{- {s-s_0} }\|\widetilde f_0\|_{L_v^2(\dot B_{2,\infty}^{s_0}\cap \dot B_{2,\infty}^s)}^2+\widetilde I_{1}+\widetilde I_{2}.
\end{align}
By the inequalities \eqref{newinterpolation} and \eqref{inter}, one has
\begin{align}\label{G4.5}
\widetilde I_1+\widetilde I_2\lesssim &\,\sup_{t\geq 0} \|E(t) \|_{\dot B^{\frac{1}{2}}_{2,\infty}\cap \dot B_{2,\infty}^{\frac{1}{2}+\varepsilon} }^2  \int_0^t (1+\tau)^{-\frac{1}{2}- {s}}   \|\widetilde f(t-\tau)\|_{L_v^2(\dot B^{\frac{1}{2}}_{2,\infty}\cap \dot B_{2,\infty}^{\frac{3}{2}-\varepsilon} )}^2{\rm d}\tau \nonumber\\
&+\sup_{t\geq 0} \|E(t) \|_{\dot B^{\frac{1}{2}}_{2,\infty}\cap \dot B_{2,\infty}^{\frac{1}{2}+\varepsilon} }^2  \int_0^t (1+\tau)^{-\frac{1}{2}-s}   \|\widetilde a(t-\tau)\|_{ \dot B^{\frac{1}{2}}_{2,\infty}\cap \dot B_{2,\infty}^{\frac{3}{2}-\varepsilon}  }^2{\rm d}\tau \nonumber\\
\lesssim&\, {\widetilde{\mathcal{D}}_{\varepsilon}^2} \sup_{t\geq 0} \|E(t)\|_{H^N} ^2  \int_0^t (1+\tau)^{-\frac{1}{2}-s} (1+t-\tau)^{-\frac{1}{2}+{s_0}} {\rm d}\tau \nonumber\\
\lesssim&\, { \delta_0^2 {\widetilde{\mathcal{D}}_{\varepsilon}^2}(1+t)^{- {(s-s_0)} }},
\end{align}
where we have made use of the fact that for $s_0>-\frac{1}{2}$ and $s<\frac{1}{2}$, which guarantees $s_0-\frac{1}{2}>-1$ and $-s-\frac{1}{2}>-1$. As a result, we have the following inequality:
\begin{align} \label{techinique}
 \int_0^t (1+\tau)^{-\frac{1}{2}-s} (1+t-\tau)^{-\frac{1}{2}+ {s_0} } {\rm d}\tau\lesssim&\, \bigg( \int_0^{\frac{t}{2}}+\int_\frac{t}{2}^t\bigg)  (1+\tau)^{-\frac{1}{2}- {s} } (1+t-\tau)^{-\frac{1}{2}+ {s_0} } {\rm d}\tau\nonumber\\
 \lesssim&\, (1+t)^{-\frac{1}{2}+ {s_0} } \int_0^\frac{t}{2}(1+\tau)^{-\frac{1}{2}- {s} }{\rm d}\tau\nonumber\\
 &+(1+t)^{-\frac{1}{2}- {s} }\int_{\frac{t}{2}}^t (1+t-\tau)^{-\frac{1}{2}+ {s_0} } {\rm d}\tau\nonumber\\
 \lesssim&\, (1+t)^{-  {(s-s_0)} }.
\end{align}

Substituting \eqref{G4.5} into \eqref{G4.4} results in
\begin{align}\label{G4.6}
\|\widetilde f(t)\|_{L_v^2(\dot B_{2,\infty}^s )}\lesssim   (1+t)^{-\frac{s-s_0}{2}}\|\widetilde f_0\|_{L_v^2(\dot B_{2,\infty}^{s_0}\cap \dot B_{2,\infty}^{\frac{1}{2}})}+    \delta_0 {\widetilde{\mathcal{D}}_{\varepsilon}}(1+t)^{-\frac{s-s_0}{2}},
\end{align}
 for any $s\in \big [-\frac{1}{2}+\varepsilon,\frac{1}{2}\big)$.

For the case $s = \frac{1}{2}$, we decompose $\widetilde f(t,x,v)$ into $\widetilde f_{L}(t,x,v)+\widetilde f_{H}(t,x,v)$. By applying the same argument for $f_{L}(t,x,v)$ as in Lemma \ref{L3.1} and using \eqref{B.24} for $f_{H}(t,x,v)$, we have
\begin{align}\label{G4.7}
\|\widetilde f(t)\|_{L_v^2(\dot B_{2,\infty}^{\frac{1}{2}} )}\lesssim&\, \|e^{t \mathcal{B}} \widetilde f_0\|_{L_v^2(\dot B_{2,\infty}^{\frac{1}{2}})}+\bigg \|\int_0^t e_{L}^{t\mathcal{B}}\widetilde S_{f}(t-\tau){\rm d}\tau \bigg\|_{L_v^2(\dot B_{2,\infty}^{\frac{1}{2}})}   +\bigg \|\int_0^t e_{H}^{t\mathcal{B}}\widetilde S_{f}(t-\tau){\rm d}\tau \bigg\|_{L_v^2(\dot B_{2,\infty}^{\frac{1}{2}})}   \nonumber\\    
\lesssim&\, \|\widetilde f_0\|_{L_v^2(\dot B_{2,\infty}^{\frac{1}{2}})}+\sup_{t\geq 0} \bigg\{\|\widetilde{G}(t)\|_{L_v^2(\dot B_{2,\infty}^{-\frac{1}{2}} )}+ \|\widetilde a {E}(t)\|_{ \dot B_{2,\infty}^{-\frac{1}{2}} }\bigg\}\nonumber\\
&+\sup_{t\geq 0} \bigg\{\|\widetilde{G}(t)\|_{L_v^2(\dot B_{2,\infty}^{ \frac{1}{2}} )} + \|\widetilde a{E}(t)\|_{ \dot B_{2,\infty}^{ \frac{1}{2}} } \bigg\}\bigg(\int_0^\infty e^{-\lambda_{2}t} {\rm d}t \bigg)^{\frac{1}{2}}\nonumber\\
\lesssim&\, (1+t)^{-\frac{s-s_0}{2}}\|\widetilde f_0\|_{L_v^2(\dot B_{2,\infty}^{s_0})}+\sup_{t\geq 0} \bigg\{\| {E}\widetilde f(t)\|_{L_v^2(\dot B_{2,\infty}^{ -\frac{1}{2}}  \cap\dot B_{2,\infty}^{ \frac{1}{2}} )} + \|\widetilde a {E}(t)\|_{\dot B_{2,\infty}^{ -\frac{1}{2}}  \cap\dot B_{2,\infty}^{ \frac{1}{2}}  } \bigg\}\nonumber\\
\lesssim&\, (1+t)^{-\frac{s-s_0}{2}}\|\widetilde f_0\|_{L_v^2(\dot B_{2,\infty}^{s_0})}+\sup_{t\geq 0}\|E(t)\|_{\dot B^{-\frac{3}{2}}_{2,\infty}\cap\dot H^N} {\sup_{t\geq 0}}\|\widetilde f(t)\|_{L_v^2(\dot B_{2,\infty}^{\frac{1}{2}})} \nonumber\\
\lesssim&\,  (1+t)^{-\frac{s-s_0}{2}}\|\widetilde f_0\|_{L_v^2(\dot B_{2,\infty}^{s_0})}+    \delta_0 {\widetilde{\mathcal{D}}_{\varepsilon}}(1+t)^{-\frac{s-s_0}{2}},
\end{align}
for $s=s_0=\frac{1}{2}$.

Finally, for the remaining case $s\in \big (\frac{1}{2},\frac{3}{2}-\varepsilon\big ]$, similar to \eqref{DLNGG3.8}, for any test function $\Phi\in L_v^2(\mathcal{S})$, we   conclude that
\begin{align*}
&\bigg|\bigg\langle \int_{{\frac{t}{2}}}^t e^{(t-\tau)\mathcal{B}}_{L} \widetilde S_{f}(\tau) {\rm d}\tau, \Phi      \bigg\rangle_{x,v} \bigg| \nonumber\\
\lesssim&\, {\sup_{\frac{t}{2}\leq \tau\leq t}}\bigg\{\|  E\widetilde{f}  (\tau) \|_{L_v^2(\dot B_{2,\infty}^{s-2+1})}+\|  \widetilde{a} E (\tau) \|_{ \dot B_{2,\infty}^{s-2+1} } \bigg\}  \sum_{j\leq j_0+N}2^{-sj}\|  \dot\Delta_{{j}}\Phi\|_{L_{x,v}^2 }  \nonumber\\
\lesssim &\,{\sup_{\frac{t}{2}\leq \tau\leq t}}\bigg\{\|   E\widetilde{f}  (\tau) \|_{L_v^2(\dot B_{2,\infty}^{s-1})}+\|  \widetilde{a} E (\tau) \|_{ \dot B_{2,\infty}^{s-1} } \bigg\}  \|\Phi\|_{L_v^2(\dot B^{-s}_{2,1})}.
\end{align*}
For $s=\frac12$ and $s_0<\frac12$, the corresponding
estimate follows from the argument of Lemma \ref{L3.1}
on $[\frac{t}{2},t]$. Together with the duality argument \eqref{duality},
the product estimate \eqref{newinterpolation},
and \eqref{B.23} and \eqref{G4.5}, lead to
\begin{align}\label{G4.8}
\bigg \|\int_0^t e_{L}^{{\tau}\mathcal{B}}\widetilde S_{f}(t-\tau){\rm d}\tau \bigg\|_{L_v^2(\dot B_{2,\infty}^{s})} \lesssim&\,
{
\bigg\|
\int_0^{\frac{t}{2}}e_L^{(t-\tau)\mathcal B}
\widetilde S_f(\tau){\rm d}\tau
\bigg\|_{L_v^2(\dot B_{2,\infty}^{s})}
+
\bigg\|
\int_{\frac{t}{2}}^t e_L^{(t-\tau)\mathcal B}
\widetilde S_f(\tau){\rm d}\tau
\bigg\|_{L_v^2(\dot B_{2,\infty}^{s})}
}
\nonumber\\
\lesssim&\,   {
\delta_0\widetilde{\mathcal D}_{\varepsilon}
(1+t)^{-\frac14-\frac{s}{2}}
\bigg(
\int_0^{\frac{t}{2}}(1+\tau)^{-\frac12+s_0}{\rm d}\tau
\bigg)^{\frac12}
}\nonumber\\
&+{ \sup_{ \frac{t}{2}\leq \tau\leq t}}\bigg\{\|   E\widetilde{f}  (\tau) \|_{L_v^2(\dot B_{2,\infty}^{s-1})}+\|  \widetilde{a} E (\tau) \|_{ \dot B_{2,\infty}^{s-1}} \bigg\}\nonumber\\
\lesssim&\, {
\delta_0\widetilde{\mathcal D}_{\varepsilon}
(1+t)^{-\frac{s-s_0}{2}}
+
\delta_0\widetilde{\mathcal D}_{\varepsilon}
\sup_{\frac{t}{2}\leq \tau\leq t}(1+\tau)^{-\frac{s-s_0}{2}}
}
\nonumber\\
\lesssim&\, \delta_0 { \widetilde{\mathcal{D}}_{\varepsilon} } (1+t)^{-\frac{s-s_0}{2}},
\end{align}
where we used $s_0>-\frac12$.
Besides, through \eqref{B.24} and direct calculation, we have
\begin{align}\label{G4.9}
\bigg \|\int_0^t e_{H}^{{\tau}\mathcal{B}}\widetilde S_{f}(t-\tau){\rm d}\tau \bigg\|_{L_v^2(\dot B_{2,\infty}^{s})}  \lesssim&\,
{
\bigg[
\int_0^t e^{-\lambda_2\tau}
\bigg(
\|E\widetilde f(t-\tau)\|_{L_v^2(\dot B_{2,\infty}^{s})}^2
+\|\widetilde aE(t-\tau)\|_{\dot B_{2,\infty}^{s}}^2
\bigg){\rm d}\tau
\bigg]^{\frac12}}
\nonumber\\
\lesssim&\,
\sup_{t\geq0}\|E(t)\|_{\dot B_{2,1}^{\frac32}}
\bigg(
\int_0^t e^{-\lambda_2\tau}
\|\widetilde f(t-\tau)\|_{L_v^2(\dot B_{2,\infty}^{s})}^2
{\rm d}\tau
\bigg)^{\frac12}
\nonumber\\
\lesssim&\,
{
\delta_0\widetilde{\mathcal D}_{\varepsilon}
\bigg(
\int_0^t e^{-\lambda_2\tau}
(1+t-\tau)^{-(s-s_0)}{\rm d}\tau
\bigg)^{\frac12}
}
\nonumber\\
\lesssim&\,
\delta_0{\widetilde{\mathcal D}_{\varepsilon}}
(1+t)^{-\frac{s-s_0}{2}}.
\end{align}

By combining the estimates \eqref{G4.8}--\eqref{G4.9}, we  further derive
\begin{align}\label{G4.10}
\|\widetilde f(t)\|_{L_v^2(\dot B_{2,\infty}^{s} )}\lesssim&\,    (1+t)^{-\frac{s-s_0}{2}}\|\widetilde f_0\|_{L_v^2(\dot B_{2,\infty}^{s_0}{\cap \dot B_{2,\infty}^{s})}}+    \delta_0 \widetilde{\mathcal{D}}_{\varepsilon} (1+t)^{-\frac{s-s_0}{2}},
\end{align}
for any $s\in \big (\frac{1}{2}, \frac{3}{2}-\varepsilon\big ]$.

Collecting all the estimates \eqref{G4.6}, \eqref{G4.7}, and \eqref{G4.10}, we consequently obtain \eqref{G4.2}.
\end{proof}

Next, we provide the estimate of $\widetilde f(t,x,v)$ in the  space $L_v^2(\dot H^2\cap \dot H^{N-1})$.
\begin{lem} \label{L4.2}
Under the assumptions of Theorem \ref{Th2}, for any $s\in\big[-\frac{1}{2}+\varepsilon,\frac{3}{2}-\varepsilon\big]$ with $s\geq s_0$ and $s_0\in\big(-\frac{1}{2},\frac{1}{2}\big]$,  we have
\begin{align}\label{G4.11}
\sup_{t\geq 0} (1+t)^{\frac{\frac{3}{2}-\varepsilon-s_0}{2}}\|\widetilde f\|_{L_v^2(\dot H^2\cap\dot H^{N-1})} \lesssim \delta_0 \widetilde{\mathcal{D}} +\|\widetilde f_0\|_{L_v^2(\dot H^2\cap \dot H^{N-1})},   
\end{align}
where    $\widetilde {\mathcal{D}}  $ is given by
\begin{align*}
\widetilde{\mathcal{D}} :=  \sup_{t\geq 0}\,(1+t)^{\frac{3}{4}-\frac{\varepsilon+s_0}{2}} \|\widetilde f (t)\|_{L_v^2(\dot  B^{\frac{3}{2}-\varepsilon}_{2,\infty} )},  
\end{align*}    
and $\delta_0$ is defined by \eqref{G1.5}.
\end{lem}

\begin{proof}
By applying the energy method to the equation \eqref{G4.1},  we have
\begin{align*} 
& \frac{{\rm d}}{{\rm d}t}\|  \widetilde f\|_{L_{v}^2(\dot H^2\cap \dot H^{N-1})}^2+2\kappa_6\sum_{2\leq k\leq N-1 } \|\nabla^k \{\mathbf{I}-\mathbf{P}\}\widetilde f\|_{\nu}^2+2\kappa_6\|  \widetilde b\|_{\dot H^2\cap\dot H^{N-1}}^2\nonumber\\
\lesssim&\,   \sum_{2\leq |\alpha|\leq N-1} |\langle \partial^\alpha[E\cdot v\widetilde f],\partial^\alpha \widetilde f\rangle _{x,v}|+\sum_{2\leq |\alpha|\leq N-1}|\langle \partial^\alpha [E\cdot\nabla_v \widetilde f],\partial^\alpha \widetilde f  \rangle_{x,v}  |\nonumber\\
\lesssim&\,  \sup_{t\geq 0} \|E\widetilde f(t)\|_{L_v^2(\dot H^2\cap \dot H^{N-1})}  \|(v\{\mathbf{I}-\mathbf{P}\}\widetilde f+v\mathbf{P}\widetilde f)(t)\|_{L_v^2( \dot H^2\cap \dot H^{N-1})} \nonumber\\
&+ \sup_{t\geq 0} \|E\widetilde f(t)\|_{L_v^2(\dot H^2\cap \dot H^{N-1})}  \|(\nabla_{v}\{\mathbf{I}-\mathbf{P}\}\widetilde f+\nabla_{v}\mathbf{P}\widetilde f)(t)\|_{L_v^2( \dot H^2\cap \dot H^{N-1})} \nonumber\\
\lesssim&\, \sup_{t\geq 0} \|E(t)\|_{\dot B^{-\frac{3}{2}}_{2,\infty}\cap\dot H^N} \|\widetilde f \|_{L_v^2(\dot B_{2,\infty}^{\frac{3}{2}-\varepsilon}\cap \dot H^{N-1})} \Big(\|(\widetilde a,\widetilde b)(t)\|_{ \dot H^2\cap \dot H^{N-1} }  +\sum_{2\leq k\leq N-1}\|\nabla^k\{\mathbf{I}-\mathbf{P}\}\widetilde f(t)\|_{\nu} \Big)\nonumber\\
\lesssim&\, \delta_0 \|\widetilde f \|_{L_v^2(\dot B_{2,\infty}^{\frac{3}{2}-\varepsilon}\cap \dot H^{N-1})}^2+\delta_0\Big(\|(\widetilde a,\widetilde b)(t)\|_{ \dot H^2\cap \dot H^{N-1} }^2  +\sum_{2\leq k\leq N-1}\|\nabla^k\{\mathbf{I}-\mathbf{P}\}\widetilde f(t)\|_{\nu}^2 \Big),
\end{align*}
for some constant $\kappa_6>0$,
which implies that
\begin{align}\label{G4.12}
& \frac{{\rm d}}{{\rm d}t}\|  \widetilde f\|_{L_{v}^2(\dot H^2\cap \dot H^{N-1})}^2+\kappa_7\sum_{2\leq k\leq N-1 } \|\nabla^k \{\mathbf{I}-\mathbf{P}\}\widetilde f\|_{\nu}^2+\kappa_7\|  \widetilde b\|_{\dot H^2\cap\dot H^{N-1}}^2   \nonumber\\
&\quad \lesssim \delta_0  \Big(\|\widetilde f \|_{L_v^2(\dot B_{2,\infty}^{\frac{3}{2}-\varepsilon} )}^2+ \|\widetilde a\|_{\dot H^{N-1}}^2\Big),
\end{align}
for some constant $\kappa_7>0$. 

To absorb the term $\|\widetilde a \|_{ \dot H^{N - 1}}^2$ on the right-hand side of \eqref{G4.12}, we begin by estimating the dissipation of $\widetilde a_{H}$. By applying the operator ${\rm Id}-\dot S_{j_0}$ to \eqref{G3.13}, we obtain 
\begin{equation}\label{G4.13}
\left\{\begin{aligned}
&\partial_{t}\widetilde a_{H}+\nabla \cdot \widetilde {b}_{H} =0,\\
&\partial_{t} \widetilde {b_i}_{H}+\widetilde{b_i}_{H}+{\partial_{i} \widetilde a}_{H}+\sum_{j=1}^3\partial_j\Gamma_{ij}(\{\mathbf{I}-\mathbf{P}\}\widetilde f_{H}) =(E_i\widetilde a)_{H},    
 \end{aligned}
 \right.
\end{equation}
for $1\leq i,j\leq 3$, where $\Gamma_{ij}(\cdot)$ is defined by \eqref{gammaa}.

Then, we define the following temporal functional
$\widetilde{\mathfrak{E}}_H(t)$:
\begin{align*}
\widetilde{\mathfrak{E}} _H(t) :=   -\int_{\mathbb{R}^3} \nabla^{N-2} \widetilde a_H \nabla^{N-2} \nabla\cdot \widetilde b_H \mathrm{d} x.
\end{align*} 
By Young's inequality, one gets
\begin{align*} 
|\widetilde{\mathfrak{E}}_H(t)| \lesssim  \|\nabla^{N-2} \widetilde a_H\|_{L^2}^2+ \|\nabla^{N-1} \widetilde b_H\|_{L^2}^2\lesssim  \|\widetilde f\|_{L_v^2(\dot H^{N-1})}^2.
\end{align*}
Similar to \eqref{G3.15}, from the second equation of \eqref{G4.13}, one has
\begin{align} \label{G4.14}
   \|\nabla^{N-1}    \widetilde a_{H} \|_{L^2}^2 
= \,&\sum_{i=1}^3 \int_{\mathbb{R}^3} \nabla^{N-2} \partial_i \widetilde a_{H} \nabla^{N-2}\Big[-\partial_t \widetilde {b_i}_{H}-\widetilde {b_i}_{H}-\sum_{j=1}^3 \partial_j \Gamma_{ij}(\{\mathbf{I}-\mathbf{P}\}\widetilde f_{H})+(E_i\widetilde a)_{H} \Big] \mathrm{d} x \nonumber\\
=&-\frac{\mathrm{d}}{\mathrm{d} t} \sum_{i=1}^3 \int_{\mathbb{R}^3} \nabla^{N-2} \partial_i \widetilde a_{H} \nabla^{N-2} \widetilde{b_i}_{H} \mathrm{d} x+\sum_{i=1}^3 \int_{\mathbb{R}^3} \nabla^{N-2}\partial_i \partial_t \widetilde a_{H} \nabla^{N-2} \widetilde {b_i}_{H}\,\mathrm{d} x \nonumber\\
& -\sum_{i=1}^3 \int_{\mathbb{R}^3} \nabla^{N-2} \partial_i \widetilde a_{H} \nabla^{N-2}\Big[\widetilde{b_i}_{H}+\sum_{j=1}^3 \partial_j\Gamma_{ij}(\{\mathbf{I}-\mathbf{P}\} \widetilde f_{H})+(E_i\widetilde a)_{H}\Big] \mathrm{d} x \nonumber\\
\equiv:&\, -\frac{\mathrm{d}}{\mathrm{d} t} \widetilde{\mathfrak{E}} _H(t) +\widetilde J_1+\widetilde J_2 . 
\end{align}

Notice that  $\partial_{t}\widetilde a_{H}$ can be replaced by using the mass conservation law in \eqref{G4.13} to obtain
\begin{align}\label{G4.15}
\widetilde J_1 = -\int_{\mathbb{R}^3}\nabla^{N-2}\partial_t\widetilde a_H\nabla^{N-2}\nabla\cdot \widetilde b_H{\rm d}x=\|\nabla^{N-2}\nabla\cdot\widetilde b_H\|_{L^2}^2. 
\end{align}
It follows from  Lemma \ref{LA.1} that
\begin{align}\label{G4.16}
\widetilde J_{2}\leq  &\, \frac{1}{2}\|\nabla^{N-1} \widetilde a_H\|_{L^2}^2+C\|\nabla^{N-2} \widetilde b_H\|_{L^2}^2+C\|\nabla^{N-1}\{\mathbf{I}-\mathbf{P}\}\widetilde f_H\|_{L_{x,v}^2}^2+C \|E_i\widetilde a\|_{\dot H^{N-2}}^2 \nonumber\\
\leq&\, \frac{1}{2}\|\nabla^{N-1} \widetilde a_H\|_{L^2}^2+C\|\widetilde b\|_{\dot H^{N-1}}^2+C\|\widetilde f\|_{L_v^2(\dot H^{N-1})}^2+C\sup_{t\geq 0}\|E(t)\|_{H^N}^2\|\widetilde a\|_{\dot B_{2,\infty}^{\frac{3}{2}-\varepsilon}\cap \dot H^{N-2}} ^2.    
\end{align}
Putting the estimates \eqref{G4.15}--\eqref{G4.16} into \eqref{G4.14} yields
\begin{align}\label{G4.17}
\frac{{\rm d}}{{\rm d}t}  \widetilde{\mathfrak{E}}_H(t)+\kappa_8  \|\widetilde a_{H}\|_{\dot H^{N-1}}^2\lesssim \|\widetilde f\|_{L_v^2(\dot B^{\frac{3}{2}-\varepsilon}_{2,\infty}\cap\dot H^{N-1})}^2,
\end{align}
for some constant $\kappa_8>0$. By adding $\eta^\prime\times$  \eqref{G4.17} to \eqref{G4.12}, where $\eta^\prime>0$ is sufficiently small, and then taking $\delta_0$ sufficiently small, we get
\begin{align}\label{G4.18}
& \frac{{\rm d}}{{\rm d}t}\Big(\|  \widetilde f\|_{L_{v}^2(\dot H^2\cap \dot H^{N-1})}^2{+\eta^{\prime}\widetilde{\mathfrak E}_H(t)\Big)}+\kappa_9\sum_{2\leq k\leq N-1 } \|\nabla^k \{\mathbf{I}-\mathbf{P}\}\widetilde f\|_{\nu}^2+\kappa_9\|\widetilde a_{H}\|_{\dot H^{N-1}}^2+\kappa_9\|  \widetilde b\|_{\dot H^2\cap\dot H^{N-1}}^2   \nonumber\\
&\quad \lesssim \delta_0 \|\widetilde f \|_{L_v^2(\dot B_{2,\infty}^{\frac{3}{2}-\varepsilon} )}^2,
\end{align}
for some constant $\kappa_9>0$.
By adding the term $\kappa_{9} \|\widetilde f\|_{L_v^2(\dot B_{2,\infty}^{\frac{3}{2}-\varepsilon})}^2$ to both sides of \eqref{G4.18}, we  derive the following Lyapunov inequality:
\begin{align*}
& \frac{{\rm d}}{{\rm d}t}\Big(\|  \widetilde f\|_{L_{v}^2(\dot H^2\cap \dot H^{N-1})}^2{+\eta^{\prime}\widetilde{\mathfrak E}_H(t)\Big)}+\kappa_{10}\Big(\|  \widetilde f\|_{L_{v}^2(\dot H^2\cap \dot H^{N-1})}^2{+\eta^{\prime}\widetilde{\mathfrak E}_H(t)\Big)}
 \lesssim  \|\widetilde f \|_{L_v^2(\dot B_{2,\infty}^{\frac{3}{2}-\varepsilon} )}^2\lesssim (1+t)^{\varepsilon+s_0-\frac{3}{2}}\widetilde{\mathcal{D}}^2,
\end{align*}    
for some constant $\kappa_{10}>0$. Thanks to the smallness of $\eta^\prime$, we have  
\begin{align*}
 \|  \widetilde f\|_{L_{v}^2(\dot H^2\cap \dot H^{N-1})}^2{+\eta^{\prime}\widetilde{\mathfrak E}_H(t)}\sim\|\widetilde f\|_{L_v^2(\dot H^2\cap\dot H^{N-1})}^2.
\end{align*}
Then, using Gr\"{o}nwall's inequality, we derive \eqref{G4.11}.
\end{proof}

\subsection{Proof of Theorem \ref{Th2}}
With Lemmas \ref{L4.1} and \ref{L4.2} established, we are now in a position to prove the uniform stability of the problem \eqref{G1.4}.

\begin{proof}[Proof of Theorem \ref{Th2}]
By combining \eqref{G4.2} and \eqref{G4.11}, we obtain
\begin{align*}
\sup_{t\geq 0}(1 + t)^{\frac{s - s_0}{2}}\|\widetilde f(t)\|_{L_v^2(\dot B_{2,\infty}^{s}\cap\dot H^{N-1})}\lesssim \delta_0\widetilde{\mathcal{D}}_{\varepsilon}+\delta_0\widetilde{\mathcal{D}}+  \|\widetilde f_0\|_{L_v^2(\dot B_{2,\infty}^{s_0}\cap \dot H^{N-1}  )},  
\end{align*}
which, together with the smallness of $\delta_0$ and the definitions of $\widetilde{\mathcal{D}}_{\varepsilon}$ and $\widetilde{\mathcal{D}}$, gives rise to
\begin{align*}
\sup_{t\geq 0}(1 + t)^{\frac{s - s_0}{2}}\|\widetilde f(t)\|_{L_v^2(\dot B_{2,\infty}^{s}\cap\dot H^{N-1})}\lesssim   \|\widetilde f_0\|_{L_v^2(\dot B_{2,\infty}^{s_0}\cap \dot H^{N-1}  )}. 
\end{align*}
Hence, \eqref{G1.8} is proved.

Finally, we focus on the uniform stability estimate \eqref{DLNZ-1.10}. Applying \eqref{G4.3} and Lemmas \ref{L2.3} and \ref{LB.3}, we immediately conclude that
\begin{align*}
&\|\widetilde f(t)\|_
{L_v^2(\dot B^s_{2,\infty}\cap\dot H^{N-1})}\\
\lesssim&\,
\|e^{t\mathcal B}\widetilde f_0\|_
{L_v^2(\dot B^s_{2,\infty}\cap\dot H^{N-1})}
+\bigg\|
\int_0^t e_L^{\tau\mathcal B}
\widetilde S_f(t-\tau){\rm d}\tau
\bigg\|_{L_v^2(\dot B^s_{2,\infty})}+\bigg\|
\int_0^t e_H^{\tau\mathcal B}
\widetilde S_f(t-\tau){\rm d}\tau
\bigg\|_{L_v^2(\dot H^{N-1})}\\
\lesssim&\,
\|\widetilde f_0\|_
{L_v^2(\dot B^s_{2,\infty}\cap\dot H^{N-1})}+\bigg[
\int_0^t(1+\tau)^{-\frac32}
\Big(
\|\widetilde G(t-\tau)\|_
{L_v^2(\dot B^{s-\frac32}_{2,\infty})}^2
+\|\widetilde aE(t-\tau)\|_
{\dot B^{s-\frac32}_{2,\infty}}^2
\Big){\rm d}\tau
\bigg]^{\frac12}\\
&+\bigg[
\int_0^t e^{-\lambda_5\tau}
\Big(
\|\widetilde G(t-\tau)\|_
{L_v^2(\dot H^{N-1})}^2
+\|\widetilde aE(t-\tau)\|_{\dot H^{N-1}}^2
\Big){\rm d}\tau
\bigg]^{\frac12}\\
\lesssim&\,
\|\widetilde f_0\|_
{L_v^2(\dot B^s_{2,\infty}\cap\dot H^{N-1})}+\sup_{t\geq0}\|E(t)\|_{\dot B^0_{2,\infty}}
\sup_{t\geq0}
\|\widetilde f(t)\|_{L_v^2(\dot B^s_{2,\infty})}
\bigg(\int_0^\infty(1+\tau)^{-\frac32}
 {\rm d}\tau\bigg)^{\frac12}\\
&+\sup_{t\geq0}
\|E(t)\|_{L^\infty\cap\dot H^{N-1}}
\sup_{t\geq0}
\|\widetilde f(t)\|_
{L_v^2(\dot B^s_{2,\infty}\cap\dot H^{N-1})}
\bigg(\int_0^\infty e^{-\lambda_5\tau}
 {\rm d}\tau\bigg)^{\frac12}\\
\lesssim&\,
\|\widetilde f_0\|_
{L_v^2(\dot B^s_{2,\infty}\cap\dot H^{N-1})}
+\sup_{t\geq0}
\|E(t)\|_{\dot B^{-\frac32}_{2,\infty}\cap\dot H^N}
\sup_{t\geq0}
\|\widetilde f(t)\|_
{L_v^2(\dot B^s_{2,\infty}\cap\dot H^{N-1})}\\
\lesssim&\,
\|\widetilde f_0\|_
{L_v^2(\dot B^s_{2,\infty}\cap\dot H^{N-1})}
+\delta_0\sup_{t\geq0}
\|\widetilde f(t)\|_
{L_v^2(\dot B^s_{2,\infty}\cap\dot H^{N-1})}.
\end{align*}
Taking the supremum over $t\geq0$ and using the
smallness of $\delta_0$, we conclude that
\begin{align*} 
\sup_{t\geq0}
\|\widetilde f(t)\|_
{L_v^2(\dot B^s_{2,\infty}\cap\dot H^{N-1})}
\lesssim
\|\widetilde f_0\|_
{L_v^2(\dot B^s_{2,\infty}\cap\dot H^{N-1})},
\end{align*}
for any $s\in(\frac12,\frac32-\varepsilon]$, which gives rise to 
\eqref{DLNZ-1.10}.
Thus, we complete the proof of Theorem \ref{Th2}.
\end{proof}

\section{Existence and stability of time-periodic solutions to the KFP equation}\label{sec5}
In this section, with the help of Theorem \ref{Th1} and Theorem \ref{Th2} concerning the Cauchy problem \eqref{G1.4}, we have got ready to prove Theorem \ref{Th1.3}, that is to establish the existence and stability of time-periodic solutions for the time-periodic problem \eqref{G1.9}. Our approach is built on Serrin’s  method (see \cite{Serrin-1959-ARMA,Mp-1991-Nonlinearity,DL-2015-AMSSB}), and we need to carefully adapt these ideas to the KFP equation.

\begin{proof}[Proof of Theorem \ref{Th1.3}]
First of all, in terms of Theorem \ref{Th1}, we let $f^*(t,x,v)$ be the unique global solution of \eqref{G1.4} with the initial data given by
\begin{align}\label{pf.add.id}
f^*(0,x,v)=f^*_0(x,v):=\frac{F_0^*(x,v)-M}{\sqrt{M}},  
\end{align}
where particularly we can choose $F_0^*(x,v)=M$ so that $f^*_0(x,v)\equiv 0$.
Then,  based on \eqref{G1.5} in Theorem \ref{Th1},  we infer that
\begin{align}\label{G5.1}
\|f^*(t)\|_{L_v^2(\dot B_{2,\infty}^{\frac{1}{2}}\cap \dot H^N)}
\lesssim    \|S (t)\|_{\mathcal{C} ([0,\infty);L_v^2(\dot B_{2,\infty}^{-\frac{3}{2}}\cap \dot H^N) )}+\|E (t)\|_{\mathcal{C} ([0,\infty); \dot B_{2,\infty}^{-\frac{3}{2}}\cap \dot H^N )}\lesssim \delta_0 , 
\end{align}
for any $t\geq 0$. Here, the solution  $f^*(t,x,v)\in \mathcal{C}([0,\infty);L_v^2(\dot B_{2,\infty}^{\frac{1}{2}}\cap \dot H^N))$, where $N\geq 3$ is an integer.

Let $m\geq k\geq 1$ be arbitrary integers. Note that $f^*(t+(m-k)T,x,v)$ is the global solution of \eqref{G1.4} with the initial data $f^*((m-k)T,x,v)$, because $E$ and $S$ are $T$-periodic in time and solutions of \eqref{G1.4} are unique, where $T>0$ is the given period. Therefore, in terms of Theorem \ref{Th2}, it follows from \eqref{G1.8}, \eqref{pf.add.id} and \eqref{G5.1} that for any $t\geq 0$,
\begin{align}\label{G5.2}
& \|f^*(t+(m-k)T)-f^*(t)\|_{L_v^2(\dot B_{2,\infty}^{1}\cap \dot H^{N-1})}\notag \\
& \lesssim   \|f^*((m-k)T)-f^*(0)\|_{L_v^2(\dot B_{2,\infty}^{\frac{1}{2}}\cap \dot H^{N-1})} (1+t)^{-\frac{1}{4}} \nonumber\\
& \lesssim \delta_0(1+t)^{-\frac{1}{4}}.
\end{align}
In particular, taking $t=kT$ in \eqref{G5.2}, we have 
$$
\|f^*(mT)-f^*(kT)\|_{L_v^2(\dot B_{2,\infty}^{1}\cap \dot H^{N-1})}\lesssim \delta_0(1+kT)^{-\frac{1}{4}} \to 0 \text{ as } k\to \infty.
$$
Note that according to \eqref{G5.1}, we have
\begin{align*}
\|f^*(mT)\|_{L_v^2(\dot B_{2,\infty}^{1}\cap\dot H^{N-1})}\lesssim   \|f^*(mT)\|_{L_v^2(\dot B_{2,\infty}^{\frac{1}{2}}\cap\dot H^{N-1})}\lesssim \delta_0.  
\end{align*}
Hence, it follows that $\{f^*(mT)\}_{m\geq 1}$ is a Cauchy sequence in $L_v^2(\dot B_{2,\infty}^{1} \cap \dot H^{N - 1})$. Thus, due to the Banach completeness theorem and Fatou's lemma, there exists a limit $f^{*}_{\infty}=f^{*}_{\infty}(x,v)$ that satisfies
\begin{align}\label{G5.3}
     f^*(mT) \rightarrow f^*_{\infty} \quad \text{in}\quad L_v^2(\dot B_{2,\infty}^{1}\cap\dot H^{N-1}),\quad \text{as}\quad m\rightarrow\infty,
\end{align}
and 
\begin{align}\label{G5.4}
\|f^*_{\infty}\|_{L_v^2(\dot B_{2,\infty}^{\frac{1}{2}}\cap\dot H^{N-1})}\leq   \liminf_{m\to\infty}  \|f^*(mT)\|_{L_v^2(\dot B_{2,\infty}^{\frac{1}{2}}\cap\dot H^{N-1})}\lesssim \delta_0,   
\end{align}
which implies $f^*_{\infty}\in L_v^2(\dot B_{2,\infty}^{\frac{1}{2}}\cap\dot H^{N-1})$.

We are now ready to establish the existence of solutions to the time-periodic problem \eqref{G1.9}. Let $f_{T}(t)$ be the unique global solution of the Cauchy problem \eqref{G1.4} with the initial data $f_{T}(0)=f_{\infty}^{*}$ satisfying \eqref{G5.4}.
From \eqref{G5.3} and in terms of Theorem \ref{Th2}, we have
\begin{align*} 
\|f_{T}(T)- f_{\infty}^*\|_{L_v^2(\dot B_{2,\infty}^{1}\cap\dot H^{N-1})}
=&\,
 \lim\limits_{n\rightarrow\infty} \|f_{T}(T)-f^*(T+(n-1)T)\|_{L_v^2(\dot B_{2,\infty}^{1}\cap\dot H^{N-1})}\nonumber\\ 
 \lesssim&\,  \lim\limits_{n\rightarrow\infty}\|f_{T}(0)-f^*((n-1)T)\|_{L_v^2(\dot B_{2,\infty}^{1 }\cap\dot H^{N-1})}\nonumber\\ 
=&\, \|f_{T}(0)-f_{\infty}^*\|_{L_v^2(\dot B_{2,\infty}^{1 }\cap\dot H^{N-1})}=0,
\end{align*}
which implies that $f_{T}(T)=f_{\infty}^*=f_T(0)$ in $L_v^2(\dot B_{2,\infty}^{1 }\cap\dot H^{N-1})$. With this fact together with Theorem \ref{Th2} again, we are able to obtain 
\begin{align}\label{pf.add.pseq}
f_T(t,x,v)=   f_T(t+T,x,v), \ \ \forall\,t\in \R.
\end{align}
Here, for $t\geq 0$, \eqref{pf.add.pseq} directly follows from Theorem \ref{Th2}, and for $t<0$, we can define $f_T(t,x,v)$ by solving the shifted Cauchy problem \eqref{G1.4} over $t\geq t_n:=-nT$ with $f_T(t_n)=f_{\infty}^*$ for $n=1,2,\dots$ in an iterative way, so that  \eqref{pf.add.pseq} follows by the same strategy.  
Thus, $f_{T}(t)$ is the unique time-periodic solution to problem \eqref{G1.9} and satisfies the estimate \eqref{Th1.3.ad1} in terms of Theorem \ref{Th1}. 

Finally, we turn to prove the time decay estimate \eqref{G1.15}.
By utilizing the embedding $L^p(\mathbb{R}^3) \hookrightarrow \dot B_{2,\infty}^{-3\left( \frac{1}{p}-\frac{1}{2}\right)}(\mathbb{R}^3)$  in Lemma \ref{L2.2} and the time-weighted estimate \eqref{G1.8}, we carry out the computation as follows:
\begin{align*}
    \|(f - f_T)(t)\|_{L_v^2(\dot B^{ s}_{2,1} )}  \lesssim&\,   \|(f - f_T)(t)\|_{L_v^2(\dot B_{2,\infty}^{-\frac{1}{2}+\varepsilon} )}^{\zeta } \|(f - f_T)(t)\|_{L_v^2(\dot B_{2,\infty}^{\frac{3}{2}-\varepsilon} )}^{1-\zeta}\nonumber\\
  \lesssim&\, (1 + t)^{-\frac{s}{2}-\frac{3}{2}(\frac{1}{p}-\frac{1}{2})}\|(f - f_T)(0)\|_{L_v^2( L^p\cap\dot H^N)},
\end{align*}
for any $s\in \big(-\frac{1}{2}+\varepsilon,\frac{3}{2}-\varepsilon \big)$, where $\zeta\in(0,1)$ satisfying
\begin{align*}
\Big(-\frac{1}{2}+\varepsilon\Big)\zeta + \Big(\frac{3}{2}-\varepsilon\Big)(1 - \zeta)=s.
\end{align*}
Therefore, the proof of Theorem \ref{Th1.3} has been completed.
\end{proof}

\section{Existence and stability of stationary solutions to the KFP equation}\label{sec6}
In this section, we are devoted to proving the existence and stability of stationary solutions to the KFP equation \eqref{G1.11} in case where the given force and source terms are time-independent. Indeed, to show Theorem \ref{cor1.4}, it suffices to obtain the following result.

\begin{thm}[Small-amplitude stationary solution obtained from Theorem \ref{Th1.3}]\label{Th6}
Let $N\ge 3$, and we define
\begin{align*}
 X:=L_v^2\big(\dot B^{\frac{1}{2}}_{2,\infty} \cap \dot H^N\big), \quad Y:=L_v^2\big(\dot B^{-\frac{3}{2}}_{2,\infty} \cap \dot H^N\big),\quad   Z:=  \dot B^{-\frac{3}{2}}_{2,\infty} \cap \dot H^N. 
\end{align*}
Assume that the external forces $E$ and $S$ are time-independent, i.e., 
\begin{align*}
E(t,x)\equiv \bar E(x), \quad S(t,x,v)\equiv \bar S(x,v),   
\end{align*}
with
\begin{align}\label{G6.1}
\|\bar S\|_Y\lesssim \delta_s,\quad \|\bar E\|_{Z}\lesssim\delta_{s}.   
\end{align}
Then the KFP equation \eqref{G1.11}  admits a unique   stationary solution $\bar f(x,v)$ satisfying
\begin{align}\label{independent}
\bar f\in     X,\quad \|\bar f\|_{X}\lesssim \delta_{s}.
\end{align}  
Here, $\delta_s>0$ is a uniformly sufficiently small constant.
\end{thm}
\begin{proof}
First, we fix an arbitrary period $T > 0$. Since the force and source terms are independent of time, that is,
$E(t,x) \equiv \bar{E}(x)$ and $S(t,x,v) \equiv \bar{S}(x,v)$,
we obtain
\begin{align*}
E(t + T,x) = \bar{E}(x) = E(t,x),\quad S(t + T,x,v) = \bar{S}(x,v) = S(t,x,v),
\end{align*}
for any $(t,x,v) \in \mathbb{R} \times \mathbb{R}_{x}^3 \times \mathbb{R}^3_{v}$.
Moreover,
\begin{align*}
\|E\|_{C(\mathbb R;Z)}
=
\sup_{t\in\mathbb R}\|E(t)\|_Z
=
\|\bar E\|_Z,\quad     \|S\|_{C(\mathbb R;Y)}
=
\sup_{t\in\mathbb R}\|S(t)\|_Y
=
\|\bar S\|_Y.
\end{align*}
Therefore, if \eqref{G6.1} holds, where $\delta_s$ is sufficiently small, then the assumptions of Theorem \ref{Th1.3} are satisfied.
Then, we apply Theorem \ref{Th1.3} to obtain a unique small-amplitude $T$-periodic solution, namely, there exists a unique time-periodic solution
$f_{T}$ such that
\begin{align}\label{add.s7.tres}
f_T\in C(\mathbb R;X),\quad 
\|f_{T}\|_{C(\mathbb R;X)}\lesssim \delta_s,    
\end{align}
and
\begin{align*}
  f_T (t+T,x,v)=f_T (t,x,v), 
\end{align*}
for any $(t,x,v)\in \mathbb R\times\mathbb R^3_{x} \times\mathbb R^3_{v}$.

Next, let's introduce an arbitrary time-shift of the periodic solution.
For any $\tau\in \mathbb R$, we define
\begin{align*}
f_T^\tau(t,x,v):=f_T(t+\tau,x,v) .    
\end{align*}
We claim that $f_T^\tau(t,x,v)$ 
is also a small-amplitude $T$-periodic solution to the KFP problem \eqref{G1.9}.  
Indeed, since $f_{T}$ solves
\begin{align}   \label{G6.2}
\partial_t f_{T} + v\cdot\nabla   f_{T} - \mathcal{L}f_{T} =   \bar E\cdot \Big(v\sqrt{M} -  \nabla_{v} f_{T}+\frac{1}{2}  v f_{T}\Big)+ \bar S,
\end{align} 
we obtain
\begin{align*}
 \partial_t f_T^\tau(t,x,v) =
\partial_t f_T(t+\tau,x,v),   
\end{align*} 
and
\begin{align*}
 &- v\cdot\nabla  f^\tau_{T} +\mathcal{L}f^\tau_{T} +  \bar E\cdot \Big(v\sqrt{M} -  \nabla_{v} f^\tau_{T}+\frac{1}{2}  v f^\tau_{T}\Big)+ \bar S\nonumber\\
& \quad=
 \,- v\cdot\nabla  f_{T}(t+\tau) +\mathcal{L}f_{T}(t+\tau) +  \bar E\cdot \Big(v\sqrt{M} -  \nabla_{v} f_{T}(t+\tau)+\frac{1}{2}  v f_{T}(t+\tau)\Big)+\bar S.   
\end{align*}
Hence, it holds that
\begin{align*}
 \partial_t f^\tau_{T} + v\cdot\nabla  f^\tau_{T} - \mathcal{L}f_{T}^\tau =  \bar E\cdot \Big(v\sqrt{M} -  \nabla_{v} f_{T}^\tau+\frac{1}{2}  v f_{T}^\tau\Big)+ \bar S.   
\end{align*}
Therefore, $f_T^\tau$ solves the same KFP equation \eqref{G6.2}
with the same force $\bar E(x)$ and  source $\bar S(x,v)$. 
Since $f_T $ is $T$-periodic, we have
\begin{align*}
f^\tau_{T}(t+T,x,v)=f_T (t+T+\tau,x,v)=f_T (t+\tau,x,v),    
\end{align*}
which yields
$f_T^\tau(t+T,x,v) = f_T^\tau (t,x,v)$.    
Thus, $f^\tau_{T} $ is also $T$-periodic.

By the definition of the norm for $C(\mathbb R;X)$, we have
\begin{align*} 
\|f_T^\tau \|_{C(\mathbb R;X)}=\|f_T\|_{C(\mathbb R;X)}\lesssim \delta_s,
\end{align*}
which implies that $f_T^\tau $ is also another small $T$-periodic solution belonging to the same class as $f_T $. Noticing the uniqueness of $f_{T}$ in Theorem \ref{Th1.3}, we then have
$f_T^\tau = f_T $ for any $ \tau\in \mathbb R $,
namely,
\begin{align}\label{G6.3}
 f_T (t+\tau,x,v)
=
f_T(t,x,v),
\end{align}
for any $t,\tau \in \mathbb R$ and $(x,v)\in \mathbb R_{x}^3\times \mathbb R_{v}^3$.
Taking $t=0$ in \eqref{G6.3}, one has
$f_T (\tau,x,v)
=f_T (0,x,v)$, for any   $(x,v)\in \mathbb R_{x}^3\times \mathbb R_{v}^3$ and
for any $\tau\in \mathbb R$.
Thus, the solution is independent of time. Let
$\bar f(x,v) := f_T (0,x,v)$.   
Then, it holds that 
$f_T(t,x,v) \equiv \bar f (x,v)$. 
As a result, we have shown that $\bar f (x,v)$ is a stationary solution.
Since the solution is independent of time, we have $\partial_t f_T  =0$.
Substituting this into the KFP probelm \eqref{G1.9} yields
\begin{align}\label{G6.4}
  v\cdot\nabla   \bar f - \mathcal{L}\bar f  =   \bar E\cdot \Big(v\sqrt{M} -  \nabla_{v} \bar f +\frac{1}{2}  v \bar f \Big)+ \bar S.
\end{align}
Thus, $\bar f(x,v)$ is a stationary solution of  \eqref{G1.11}, and \eqref{independent} directly follows from \eqref{add.s7.tres}. 
 
Finally, we show  the uniqueness of the stationary solution $\bar f (x,v)$.
Assume that $\bar f^{(1)}$ and $\bar f^{(2)}$ are two stationary solutions within the same small class to the equation \eqref{G6.4}. Since they are stationary, they can be regarded as time-periodic solutions for any $T > 0$ by setting
\begin{align*}
f^{(j)}(t,x,v):= \bar f^{(j)}(x,v), \quad j = 1,2.
\end{align*}
Then both are small $T$-periodic solutions of the same KFP equation with the same force $\bar E(x)$ and source $\bar S(x,v)$. According to the uniqueness stated in Theorem \ref{Th1.3}, we have $\bar f^{(1)}(x,v)=\bar f^{(2)}(x,v)$. Therefore, we complete the proof of Theorem \ref{Th6}.
\end{proof}

\medskip 
\noindent{\bf Acknowledgements:} 
The research of Renjun Duan was partially supported by the General Research Fund (Project No.~14303523) from RGC of Hong Kong, the grant from the National Natural Science Foundation of China (Project No.~12425109) and the Direct Grant (4053652) from CUHK. The research of both Fucai Li and Jinkai Ni were supported by NSFC (Grant No. 12331007). The research of Fucai Li was also supported by the ``333 Project" of Jiangsu Province. The research of Zhu Zhang was supported by the Early Career Scheme (No. 25303523) of the Hong Kong Research Grants Council (RGC).

\medskip 
\noindent\textbf{Conflict of interest.} The authors do not have any possible conflicts of interest.

\medskip 
\noindent\textbf{Data availability statement.}
 Data sharing is not applicable to this article as no data sets were generated or analyzed during the current study.

\bibliographystyle{plain}

\begin{thebibliography}{aaa}

\bibitem{AAMN-APDE-2024} D. Albritton,  S. Armstrong,  J.-C. Mourrat,  M. Novack, 
Variational methods for the kinetic Fokker-Planck equation,
{\it Anal. PDE} {\bf 17 (6)} (2024)  1953--2010.

\bibitem{AMTU-CPDE-2001} A. Arnold,  P. Markowich, G. Toscani,  A. Unterreiter,
On convex Sobolev inequalities and the rate of convergence to equilibrium for Fokker-Planck type equations,
{\it Comm. Partial Differential Equations} {\bf 26} (2001)  43--100.

\bibitem{ah} B. Avelin, M.  Hou,  Weak and Perron solutions for stationary Kramers-Fokker-Planck equations in bounded domains, {\it Potential Anal.} {\bf 64 (1)} (2026)  Paper No. 11.  

\bibitem{BCD-Book-2011} H. Bahouri, J. Chemin, R. Danchin,
{\it  Fourier Analysis and Nonlinear Partial Differential Equations,}
Grundlehren Math. Wiss, 343,
Springer, Heidelberg, 2011. 

\bibitem{BBJM-esaim-2005} C. Baranger, L. Boudin, P. Jabin, S. Mancini,
A modeling of biospray for the upper airways, CEMRACS 2004 -- {\it Mathematics and Applications to Biology and Medicine}, 41--47, ESAIM Proc.,  14, EDP Sciences, Les Ulis, 2005.

\bibitem{Bayraktar-Feng-Li-2024} E. Bayraktar, Q. Feng,   W. Li, 
Exponential entropy dissipation for weakly self-consistent Vlasov-Fokker-Planck equations,
{\it J. Nonlinear Sci.} {\bf 34 (1)} (2024)   Paper No. 7.

\bibitem{bzz} J.  Bedrossian, W.  Zhao, R.   Zi,  
Landau damping, collisionless limit, and stability threshold for the Vlasov-Poisson equation with nonlinear Fokker-Planck collisions, 
{\it Comm. Math. Phys.}  {\bf 406 (7)} (2025)  Paper No. 166.

\bibitem{bit}
 J. A. Bittencourt, \emph{Fundamentals of plasma physics}. Springer Science and Business Media, 1986.

\bibitem{Bolley-Guillin-Malrieu-2010}F. Bolley,  A. Guillin,  F. Malrieu,  
Trend to equilibrium and particle approximation for a weakly selfconsistent Vlasov-Fokker-Planck equation,
{\it M2AN Math. Model. Numer. Anal.} {\bf 44 (5)} (2010)  867--884.

 \bibitem{BG} M. Bostan, T.  Goudon,  High-electric-field limit for the Vlasov-Maxwell-Fokker-Planck system, {\it Ann. Inst. H. Poincar\'{e} C Anal. Non Lin\'{e}aire} {\bf 25 (6)} (2008)  1221--1251.  


\bibitem{Bouchut-Dolbeault-1995} F. Bouchut, J. Dolbeault,  
On long time asymptotics of the Vlasov-Fokker-Planck equation and of the Vlasov-Poisson-Fokker-Planck system with Coulombic and Newtonian potentials,
{\it Differential Integral Equations} {\bf 8 (3)} (1995)  487--514.

 \bibitem{BS-sima25} G.  Brigati, G.  Stoltz,   How to construct explicit decay rates for kinetic Fokker-Planck equations, {\it SIAM J. Math. Anal.} {\bf 57 (4)} (2025)   3587--3622.  

\bibitem{Ccq-CMS-2019} C. Cao, 
The kinetic Fokker-Planck equation with weak confinement force,
{\it Commun. Math. Sci.} {\bf 17 (8)} (2019)  2281--2308.

\bibitem{Ccq-JEE-2021}C. Cao,  
The kinetic Fokker-Planck equation with general force,
{\it J. Evol. Equ.} {\bf 21 (2)} (2021)  2293--2337.


\bibitem{cm} 
K. Carrapatoso, S.  Mischler, 
The kinetic Fokker-Planck equation in a domain: ultracontractivity, hypocoercivity, and long-time asymptotic behavior,
{\it Atti Accad. Naz. Lincei Rend. Lincei Mat. Appl.} {\bf 35 (4)} (2024) 643--680. 

\bibitem{CDM-krm-2011} 
J.-A. Carrillo, R.-J. Duan, A. Moussa, Global classical solutions close to equilibrium 
to the Vlasov-Fokker-Planck-Euler system, {\it Kinet. Relat. Models} {\bf 4 (1)} (2011) 227--258.

\bibitem{chae} M.  Chae,  The global classical solution of the Vlasov-Maxwell-Fokker-Planck system near Maxwellian, {\it Math. Models Methods Appl. Sci.} {\bf 21 (5)} (2011)  1007--1025.  


 \bibitem{cha} S. Chandrasekhar, Stochastic problems in physics and astronomy, {\it Reviews of Modern Physics.}
 {\bf 15  (1)}  (1943) 1--89. 

\bibitem{chav} P.-H. Chavanis, Generalized thermodynamics and Fokker-Planck equations: Applications to stellar dynamics and two-dimensional turbulence, {\it Phys. Rev. E} {\bf 68} (2003) 036108.

\bibitem{Choi-2016-nonlinearity} Y.-P. Choi,  
Global classical solutions of the Vlasov-Fokker-Planck equation with local alignment forces,
{\it Nonlinearity} {\bf 29 (7)} (2016)  1887--1916. 

\bibitem{Dk-MZ-1992} K. Deckelnick, 
Decay estimates for the compressible Navier-Stokes equations in unbounded domains,
{\it Math. Z.}   
{\bf 209 (1)} (1992) 115--130.

\bibitem{Degond-1986} P. Degond,  
Global existence of smooth solutions for the Vlasov-Fokker-Planck equation in 1 and 2 space dimensions,
{\it Ann. Sci. \'{E}cole Norm. Sup. (4)} {\bf 19 (4)}  (1986)  519--542.


\bibitem{Deguchi-2024-Math Ann} N. Deguchi,  
On the stability of stationary compressible Navier-Stokes flows in 3D,
{\it Math. Ann.} {\bf 390 (3)} (2024)  4361--4404.

\bibitem{Deguchi-2025} N. Deguchi, Stability analysis of time-periodic solutions to the Navier-Stokes-Fourier system in 3D whole space, arXiv:2601.00034v2.

\bibitem{doi}
 M. Doi, S. F. Edwards, {\it The Theory of Polymer Dynamics,} International Series of Monographs on Physics, Clarendon Press, 1988.  
 
\bibitem{do} J. Dolbeault,  
Free energy and solutions of the Vlasov-Poisson-Fokker-Planck system: external potential and confinement (large time behavior and steady states), 
{\it J. Math. Pures Appl. (9)} {\bf 78 (2)} (1999)  121--157. 

\bibitem{do-li} J. Dolbeault, X.  Li,  
{$\varphi$}-entropies: convexity, coercivity and hypocoercivity for Fokker-Planck and kinetic Fokker-Planck equations,  
{\it Math. Models Methods Appl. Sci.} {\bf 28 (13)} (2018)  2637--2666. 

 \bibitem{dms}   J. Dolbeault, C. Mouhot,  C. Schmeiser,  Hypocoercivity for linear kinetic equations conserving mass, {\it Trans. Amer. Math. Soc.,} {\bf 367 (6)} (2015) 3807--3828.   
    
 \bibitem{dre} K. Dressler,  
Stationary solutions of the Vlasov-Fokker-Planck equation,
{\it Math. Methods Appl. Sci.} {\bf 9 (2)} (1987)  169--176. 

\bibitem{D-2011-Nonlinearity} 
R.-J. Duan,  
Hypocoercivity of linear degenerately dissipative kinetic equations,
{\it Nonlinearity} {\bf 24 (8)} (2011)  2165--2189.

\bibitem{DFT-2010-CMP} 
R.-J. Duan, M. Fornasier, G. Toscani, 
A kinetic flocking model with diffusion,
{\it Comm. Math. Phys.} {\bf 300 (1)} (2010) 95--145.

\bibitem{DL-2013-KRM}  
R.-J. Duan,  S. Liu, 
Cauchy problem on the Vlasov-Fokker-Planck equation coupled with the compressible Euler equations through the friction force,
{\it Kinet. Relat. Models} {\bf 6(4)} (2013)  687--700.

\bibitem{DL-2015-AMSSB}
R.-J. Duan, S. Liu,  
Time-periodic solutions of the Vlasov-Poisson-Fokker-Planck system,
{\it Acta Math. Sci. Ser. B (Engl. Ed.)} {\bf 35 (4)} (2015)  876--886.

\bibitem{DLN-2026} 
R.-J. Duan, Y. Lei, J. Ni,   Global dynamics of the Boltzmann equation driven by a time-periodic source in $\mathbb{R}^3$, preprint, 2026.

\bibitem{DN-2026} 
R.-J. Duan, J. Ni, Three-dimensional time-periodic problem on the Boltzmann equation with external force, arXiv:2604.21339.

\bibitem{FQW-M3AS-2024} 
Z. Fang,  K. Qi,  H. Wen,  
The small Deborah number limit for the fluid-particle flows: incompressible case,
{\it Math. Models Methods Appl. Sci.} {\bf 34 (12)}  (2024)    2265--2304.

 \bibitem{ft-1} 
R. Farwig,  K.  Tsuda,  
The time periodic problem for the Navier-Stokes equations on half spaces with moving boundary: linear theory,
{\it J. Differential Equations} {\bf 411} (2024) 531--603. 

\bibitem{ft-2} R. Farwig,  K.  Tsuda,  
The time periodic problem for the Navier-Stokes equations on half spaces with moving boundary: nonlinear theory, {\it Math. Ann.} {\bf 391 (3)} (2025)  4791--4846. 


\bibitem{feic}
E. Feireisl, P. Gwiazda, A. \'{S}wierczewska-Gwiazda,  
Time periodic motion of temperature driven compressible fluids,  
{\it Math. Ann.} {\bf 387 (3--4)} (2023)  1603--1627. 

\bibitem{Frank-2005} T. D. Frank,  {\it Nonlinear Fokker–Planck Equations: Fundamentals and Applications} (Berlin: 
Springer)  2005.

\bibitem{feia}
E. Feireisl, S.  Matu\v{s}{u}-Ne\v{c}asov\'{a},  H. Petzeltov\'{a},   I. Stra\v{s}kraba, On the motion of a viscous compressible fluid driven by a time-periodic external force, {\it Arch. Ration. Mech. Anal.} {\bf 149 (1)} (1999)  69--96. 

 \bibitem{feib}
 E. Feireisl, P. B. Mucha, A. Novotn\'{y}, M.  Pokorn\'{y},  
 Time-periodic solutions to the full Navier-Stokes-Fourier system, 
{\it Arch. Ration. Mech. Anal.} {\bf 204 (3)} (2012)  745--786. 

\bibitem{gsz} R.  Glassey, J.  Schaeffer, Y.  Zheng,   Steady states of the Vlasov-Poisson-Fokker-Planck system, {\it J. Math. Anal. Appl.} {\bf 202 (3)} (1996)  1058--1075.  

\bibitem{GMM-MSM-2017}  M. P. Gualdani, S. Mischler,  C. Mouhot,  
Factorization of non-symmetric operators and exponential H-theorem,
{\it M\'{e}m. Soc. Math. Fr. (N.S.)}   {\bf  153} (2017).

\bibitem{Gy-IUMJ-2004} Y. Guo,  
The Boltzmann equation in the whole space,
{\it Indiana Univ. Math. J.} {\bf 53 (4)} (2004)  1081--1094.

\bibitem{Ks-1983}  S. Kawashima, {\it Systems of a Hyperbolic-Parabolic Composite Type, with Applications to the Equations of Magnetohydrodynamics,}
 PhD thesis, Kyoto University, 1983.

\bibitem{LMW-SIMA-2017} F. Li,  Y. Mu,  D. Wang,  
Strong solutions to the compressible Navier-Stokes-Vlasov-Fokker-Planck equations: global existence near the equilibrium and large time behavior,
{\it SIAM J. Math. Anal.} {\bf 49  (2)} (2017)  984--1026.

\bibitem{LLY-JSP-2022}  
H.-L. Li,  S. Liu, T. Yang, 
The Navier-Stokes-Vlasov-Fokker-Planck system in bounded domains,
{\it J. Stat. Phys.} {\bf 186 (3)} (2022) Paper No. 42.

\bibitem{LWW-ARMA-2022} H.-L. Li,  T. Wang,  Y. Wang,  
Wave phenomena to the three-dimensional fluid-particle model,
{\it Arch. Ration. Mech. Anal.} {\bf 243 (2)} (2022)  1019--1089.

\bibitem{lib}
R. L. Liboff, Kinetic theory: classical, quantum, and relativistic descriptions. Springer Science
and Business Media, 2003.

\bibitem{LY-AA-2015} L. Luo,  H. Yu,  
Spectrum analysis of the linear Fokker-Planck equation,
{\it Anal. Appl. (Singap.)}  {\bf 15 (3)} (2017)  313--331.

\bibitem{Mp-1991-Nonlinearity} P. Maremonti, 
Existence and stability of time-periodic solutions to the Navier-Stokes equations in the whole space,
{\it Nonlinearity} {\bf 4 (2)} (1991)  503--529.

\bibitem{MN-PJASAMS-1979} A. Matsumura,  T. Nishida,  
The initial value problem for the equations of motion of compressible viscous and heat-conductive fluids,
{\it  Proc. Japan Acad. Ser. A Math. Sci.} {\bf 55 (9)} (1979)  337--342.
 
\bibitem{mis-mou} S.  Mischler, C. Mouhot,   Exponential stability of slowly decaying solutions to the kinetic-Fokker-Planck equation, {\it Arch. Ration. Mech. Anal.} {\bf 221 (2)} (2016)  677--723.  

\bibitem{Serrin-1959-ARMA}
J. Serrin, 
A note on the existence of periodic solutions of the Navier-Stokes equations,
{\it Arch. Rational Mech. Anal.} {\bf 3} (1959) 120--122.

\bibitem{SK-HMJ-1985} Y. Shizuta, S. Kawashima, 
Systems of equations of hyperbolic-parabolic type with applications to the discrete Boltzmann equation,
{\it Hokkaido Math. J.} {\bf 14 (2)} (1985) 249--275.


 \bibitem{Ts} K. Tsuda,  On the existence and stability of time periodic solution to the compressible Navier-Stokes equation on the whole space, {\it Arch. Ration. Mech. Anal.} {\bf 219 (2)} (2016)  637--678.  

\bibitem{ukai} S. Ukai,  Time-periodic solutions of the Boltzmann equation,  
{\it Discrete Contin. Dyn. Syst.} {\bf 14 (3)} (2006) 579--596. 

\bibitem{ukaiy} S. Ukai, T.  Yang,  The {B}oltzmann equation in the space {$L^2\cap L^\infty_\beta$}: global and time-periodic solutions,  {\it Anal. Appl. (Singap.)} {\bf 4 (3)}(2006) 263--310. 

\bibitem{Valli-1983} A. Valli,  
Periodic and stationary solutions for compressible Navier-Stokes equations via a stability method,
{\it Ann. Scuola Norm. Sup. Pisa Cl. Sci. (4)} {\bf 10 (4)} (1983)  607--647.

\bibitem{Villani-MAMS-2009} C. Villani,  
Hypocoercivity,
{\it Mem. Amer. Math. Soc.} {\bf 202 (950)} (2009).

\bibitem{YY-SIMA-2010} T. Yang,  H. Yu,  
Global classical solutions for the Vlasov-Maxwell-Fokker-Planck system,
{\it SIAM J. Math. Anal.} {\bf 42 (1)} (2010)  459--488.

\end{thebibliography}

\end{document}